\documentclass[11pt,a4paper]{article}

\usepackage[english]{babel}
\usepackage{amsmath}
\usepackage{amssymb}
\usepackage{graphicx}        

\usepackage{cases}
\usepackage{amsfonts}
\usepackage{color}
\usepackage{verbatim}
\usepackage{caption}
\usepackage{subcaption}
\usepackage{fullpage}

\usepackage{array}

\usepackage{epsfig}
\usepackage{float}
\usepackage{url}

\usepackage{theorem}

\usepackage{fullpage}

\usepackage{hyperref}

\newcommand{\rosso}[1]{\color{black} #1} 

\usepackage{stackengine}
\newcommand{\tensorm}{
\setstackgap{S}{0.4ex}%
\mathrel{\Shortstack{{.} {.} {.}}}}

\newtheorem{rem}{Remark}[section]
\newtheorem{lem}{Lemma}[section]
\newtheorem{thm}{Theorem}[section]

\newenvironment{pf}{\noindent{\bf Proof. \/}\noindent%
}{\hfill\EndProofMarker}
\newcommand{\EndProofMarker}{$\Box$}
\newcommand{\diver}{\ensuremath{\operatorname{div}}}

\def\XXint#1#2#3{{\setbox0=\hbox{$#1{#2#3}{\int}$ }
\vcenter{\hbox{$#2#3$ }}\kern-.6\wd0}}

\begin{document}
  
\title{A Cahn--Hilliard phase field model coupled to an Allen--Cahn model of viscoelasticity at large strains}
\author{A. Agosti$^{\sharp}$\footnote{Corresponding author. E-mail address: {\tt abramo.agosti@unipv.it} \newline \textit{Email addresses:} {\tt abramo.agosti@unipv.it} (A. Agosti), {\tt pierluigi.colli@unipv.it} (P. Colli), {\tt harald.garcke@ur.de} (H. Garcke), {\tt elisabetta.rocca@unipv.it} (E. Rocca)}, P. Colli$^{\sharp,\ddag}$, H. Garcke$^{\S}$, E. Rocca$^{\sharp,\ddag}$}

\maketitle 

\begin{center}
{\small $^\sharp$  Department of Mathematics, University of Pavia, 27100 Pavia, Italy.}
\end{center}

\begin{center}
{\small $^\ddag$  Research Associate at the IMATI--C.N.R. Pavia, 27100 Pavia, Italy.}
\end{center}

\begin{center}
{\small $^\S$  Fakult{\"a}t f{\"u}r Mathematik, Universit{\"a}t Regensburg, 93040 Regensburg, Germany.}
\end{center}

\maketitle

\date{}

\begin{abstract}
We propose a new Cahn--Hilliard phase field model coupled to incompressible viscoelasticity at large strains, obtained from a diffuse interface mixture model and formulated in the Eulerian configuration. A new kind of diffusive regularization, of Allen--Cahn type, is introduced in the transport equation for the deformation gradient, together with a regularizing interface term depending on the gradient of the deformation gradient in the free energy density of the system.
The designed regularization preserves the dissipative structure of the equations. We study the global existence of a weak solution for the model. 
While standard diffusive regularizations of the transport equation for the deformation gradient presented in literature for related phase field models coupled to viscoelasticity allows the study of existence of global weak solutions  only for simplified cases, i.e. in two space dimensions and for convex elastic free energy densities of Neo--Hookean type which are independent from  the phase field variable, the present diffusive regularization allows to study more general cases. In particular, we obtain the global existence of a weak solution in three space dimensions and for generic nonlinear elastic energy densities with polynomial growth, comprising the relevant cases of polyconvex Mooney--Rivlin and Ogden elastic energies. Also, our analysis considers elastic free energy densities which depend on the phase field variable and which can possibly degenerate for some values of the phase field variable. By means of an iterative argument based on elliptic regularity bootstrap steps, we find the maximum allowed polynomial growths of the Cahn--Hilliard potential and the elastic energy density which guarantee the existence of a solution in three space dimensions. We also propose two kinds of unconditionally energy stable finite element approximations of the model, based on convex splitting ideas and on the use of a scalar auxiliary variable respectively, proving the existence and stability of discrete solutions. We finally report numerical results for different test cases with shape memory alloy type free energy, showing the interplay between phase separation and finite elasticity in determining the topology of stationary states with pure phases characterized by different elastic properties. 
\end{abstract}



\section{Introduction}
The study of Cahn--Hilliard phase field models coupled with finite viscoelasticity has gained increasing interest in the recent literature, cf., e.g., \cite{lukacova1,garcke5,agosti3,liu}. These models describe the phase separation phenomena for multiphase materials in presence of elastic interactions between the materials constituents, and may find applications e.g. in soft matter dynamics \cite{lukacova1}, tumor growth dynamics \cite{garcke5}, neurological and neuromuscular deseases (see the discussion in \cite{agosti3}).

In these works both the phase field and the viscoelasticity governing equations are formulated in the Eulerian reference configuration.
Hence, the main state variable for elasticity is the velocity field, and the deformation gradient, entering in the elasticity constitutive assumptions, is determined solving a transport equation in terms of the velocity and the velocity gradient. A similar modeling approach for evolutionary models with finite viscoelasticity was implemented also in
\cite{barrett}, which studies the Oldroyd-B model for a dilute polymeric fluid, in \cite{suli}, which studies the motion of a class of incompressible heat-conducting viscoelastic rate-type fluids with stress-diffusion, in \cite{benesova,garcke4}, which study problems in magnetoelasticity, in \cite{rubicek}, which studies an evolutionary non-isothermal viscoelastic problem, and in \cite{rubicek3}, which studies the diffusion of a solvent in a saturated hyperelastic porous solid of viscoelastic Kelvin-Voigt type at large strains.

A major challenge in studying the existence of weak solutions for phase field models coupled with finite viscoelasticity is to obtain results in three space dimensions and for generic elastic energy densities which are highly nonlinear in the deformation gradient and which may depend on the phase field variable. Indeed, in \cite{benesova,lukacova1,garcke4,garcke5}, employing a second order diffusive regularization in the transport equation for the deformation gradient (i.e. adding a term proportional to the Laplacian of the deformation gradient), the existence of global weak solutions for the model is proved in two space dimensions and for elastic energy densities of Neo--Hookean type, i.e. quadratic in the deformation gradient, and independent of the phase field variable. The diffusive regularization of the transport equation, which is needed to increase the regularity of the deformation gradient and gain compactness of its approximations, breaks the kinematic relationship between the velocity and the deformation gradient variables. In \cite{agosti3}, we designed a specific second order diffusive regularization of the transport equation for the deformation gradient, depending on both the phase field variable and on the deformation gradient, which allowed us to prove the existence of global weak solutions for elastic energy densities of Neo--Hookean type which depend on the phase field variable. The degeneracy of the elastic energy density depending on the phase field variable was not allowed in the framework of \cite{agosti3}. Also, existence results were obtained in three space dimensions by adding a further viscous regularization to the Cahn--Hilliard equation. By employing a different regularization approach of the model equations, obtained by adding a dissipative contribution to the Cauchy stress tensor which involves high order nonlinear terms in the small strain rate (i.e. the symmetric part of the velocity gradient), in \cite{rubicek} the author obtained existence and regularity results of a distributional solution for an evolutionary non-isothermal viscoelastic problem with nonlinear and compressible elastic energy. The latter regularization approach preserves the kinematic relationship between the velocity and the deformation gradient variables. Its main drawback is the introduction of high order nonlinear terms involving the velocity gradient in the definition of the Cauchy stress tensor, making numerical approximations of the model not straightforward. 

In the present paper, we introduce a new kind of diffusive regularization, of Allen–Cahn type, in the transport equation for the deformation gradient, complemented by the introduction of a regularizing interface term depending on the gradient of the deformation gradient in the free energy density of the system. The designed regularization preserves the dissipative structure of the equations. The idea for the introduction of this kind of regularization comes from an observation by Roubicek \cite[Remark 3]{rubicek2}, which suggested the possibility to consider a Cahn--Hilliard regularization of the transport equation for the deformation gradient in order to obtain stronger existence results with respect to those available in the literature and based on diffusive regularizations of the transport equation for the deformation gradient.
The introduction of an interface contribution for the deformation gradient in the free energy enhances the space and time regularity of the deformation gradient, while increasing the degree of nonlinearity of the coupled system. In particular, the Cauchy stress tensor contains second order and quadratic first order terms in the deformation gradient. In our framework, the regularity of the latter terms is obtained by adding an Allen--Cahn type second order diffusive regularization in the transport equation for the deformation gradient. Following the Allen--Cahn literature, we consider a dual mixed Allen--Cahn formulation of the transport equation for the deformation gradient, introducing a dual variable of the deformation gradient which enters also in the expression for the Cauchy stress tensor. This theoretical framework allows us to obtain the existence of global weak solutions in three space dimensions and for generic nonlinear elastic energy densities with polynomial growth, which are coupled to the phase field variable and which possibly degenerate for some values of the phase field variable. 

The resulting PDE system for the considered model is the following: 
\begin{equation}
\label{eqn:systemintro}
\begin{cases}
-\nu\Delta \mathbf{v}+\nabla s=\mu \nabla \phi+\diver\left(\mathbf{M}\mathbf{F}^T\right)+\left(\nabla \mathbf{F}\right)\odot\mathbf{M},\\ \\
\diver\mathbf{v}=0,\\ \\
\frac{\partial \mathbf{F}}{\partial t}+\left(\mathbf{v}\cdot \nabla\right)\mathbf{F}-(\nabla \mathbf{v})\mathbf{F}+\gamma \mathbf{M}=0,\\ \\
\mathbf{M}=\partial_{\mathbf{F}}w(\phi,\mathbf{F})-\lambda \Delta \mathbf{F},\\ \\
\frac{\partial \phi}{\partial t}+\mathbf{v}\cdot \nabla \phi-\diver(b(\phi)\nabla \mu)=0,\\ \\
\mu=\psi'(\phi)- \Delta \phi+\partial_{\phi}w(\phi,\mathbf{F}),
\end{cases}
\end{equation}
valid in $\Omega_T$, endowed with the boundary conditions
\begin{equation}
\label{eqn:systembcintro}
b(\phi)\nabla \mu\cdot \mathbf{n}=\nabla \phi \cdot \mathbf{n}=0, \quad \left[\nabla \mathbf{F}\right]\mathbf{n}=\mathbf{0}, \quad \mathbf{v}=\mathbf{0}, 
\end{equation}
on $\partial \Omega \times [0,T]$, and with proper initial conditions. Here, $\mathbf{v}$ is the velocity field, $s$ is the pressure, $\mathbf{F}$ is the elastic deformation gradient, $\mathbf{M}$ is its dual variable, $\phi$ is the phase field variable and $\mu$ is the chemical potential. {Moreover, $\nu$ is a physical parameter, representing the viscosity of the material, while $\gamma$ and $\lambda$ are positive regularization parameters. The function $b(\phi)$ represents a positive phase-dependent mobility, $\psi(\phi)$ represents the bulk potential of the Cahn--Hilliard equation, and $w(\phi,\mathbf{F})$ is the elastic free energy density.}

By means of an iterative argument based on elliptic regularity bootstrap steps applied to \eqref{eqn:systemintro}$_4$, we find the
maximum allowed polynomial growths of the Cahn--Hilliard potential and the elastic energy density which guarantee existence of a solutions in three space dimensions. 

We observe that the dual mixed formulation of System \eqref{eqn:systemintro} is particularly suitable to design standard and efficient numerical approximations. Hence, we also propose two kinds of unconditionally energy stable finite element approximations of the model, based on convex splitting ideas and on the use of a scalar auxiliary variable, proving the existence and stability of discrete solutions. We finally show numerical tests for different test cases with shape memory alloy type free energy, proving the interplay between phase separation and finite elasticity in determining the topology of stationary states with pure phases characterized by different elastic properties.

Hence, the novelties of the present work with respect to previous studies in literature are the following:
\begin{itemize}
\item Proof of existence of global weak solutions for a Cahn--Hilliard phase field model coupled to viscoelasticity at large strains in three space dimensions, based on a new Allen--Cahn type regularization of finite viscoelasticity, for generic finite elastic energy densities with polynomial growth, which are coupled to the phase field variable and which possibly degenerate for some values of the phase field variable; 
\item Design and implementation of standard and efficient well posed and unconditionally energy stable finite element approximations of the model.
\end{itemize}

The paper is organized as follows. In Section $2$ we introduce some notation regarding the employed tensor calculus and functional analysis. In Section $3$ we develop the model derivation. In Section $4$ we report the study of existence for a global weak solution of System \eqref{eqn:systemintro}. In Section $5$ we report the design of  finite element approximations of the model and the study of existence and stability of discrete solutions. In Section $6$ we show results of numerical simulations. Finally, Section $7$ contains concluding remarks and future perspectives. 

\section{Notation}
Let $\Omega \subset \mathbb{R}^3$ be an open bounded domain in $\mathbb{R}^3$. {\rosso Let $T>0$ denote some final time, and set $\Omega_T:=\Omega \times (0,T)$.}
We start by introducing the notation for vectorial and tensorial calculus. Given $\mathbf{a},\mathbf{b}\in \mathbb{R}^3$, we denote by $\mathbf{a}\cdot \mathbf{b} \in \mathbb{R}$ their canonical scalar product in $\mathbb{R}^3$, with associated norm $|\mathbf{a}|:=(\mathbf{a}\cdot \mathbf{a})^{\frac{1}{2}}$, and by $\mathbf{a}\otimes \mathbf{b}\in \mathbb{R}^{3\times 3}$ their tensorial product. We indicate by $\{\mathbf{e}_i\}_{i=1,2,3}$ the canonical basis of $\mathbb{R}^3$. Given two second order tensors $\mathbf{A},\mathbf{B}\in \mathbb{R}^{3\times 3}$, we denote by $\mathbf{A}\colon \mathbf{B}\in \mathbb{R}$ their Frobenius scalar product in $\mathbb{R}^{3\times 3}$, i.e. by components $\mathbf{A}\colon \mathbf{B}:=\sum_{i,j=1}^3A_{ij}B_{ij}$, with associated norm $|\mathbf{A}|:=(\mathbf{A}\colon \mathbf{A})^{\frac{1}{2}}$. We indicate by $\{\mathbf{e}_{ij}:=\mathbf{e}_i\otimes \mathbf{e}_j\}_{i,j=1,2,3}$ the canonical basis of $\mathbb{R}^{3\times 3}$. Given a second order tensor $\mathbf{A}\in \mathbb{R}^{3\times 3}$, we indicate with the notation $\mathbf{A}^{[i]}\in \mathbb{R}^3$ the $i$-th column vector of $\mathbf{A}$. Given two third order tensors $\mathbf{C},\mathbf{D}\in \mathbb{R}^{3\times 3 \times 3}$, we denote by $\mathbf{C}\mathbin{\tensorm}\mathbf{D}\in \mathbb{R}$ their scalar product in $\mathbb{R}^{3\times 3 \times 3}$, i.e. by components 
\[\mathbf{C}\mathbin{\tensorm} \mathbf{D}:=\sum_{i,j,k=1}^3C_{ijk}D_{ijk}.\]
We also introduce the operation, defined by components, 
\begin{align*}
\displaystyle & \left(\mathbf{C}\odot \mathbf{B}\right)_{k}:=\sum_{ij=1}^3\mathbf{C}_{ijk}\mathbf{B}_{ij},\\
\displaystyle &\left(\mathbf{C}\odot \mathbf{D}\right)_{kl}:=\sum_{ij=1}^3\mathbf{C}_{ijk}\mathbf{D}_{ijl},
\end{align*}
which contracts respectively a third order tensor $\mathbf{C}\in \mathbb{R}^{3\times 3 \times 3}$ and a second order tensor $\mathbf{B}\in \mathbb{R}^{3\times 3}$ to a vector $\mathbf{C}\odot \mathbf{B}\in \mathbb{R}^{3}$ and two third order tensors $\mathbf{C},\mathbf{D}\in \mathbb{R}^{3\times 3 \times 3}$ to a second order tensor $\mathbf{C}\odot \mathbf{D}\in \mathbb{R}^{3\times 3}$. 

We denote by $L^p(\Omega;K)$ and $W^{r,p}(\Omega;K)$ the standard Lebesgue and Sobolev spaces of functions defined on $\Omega$ with values in a set $K$, {\rosso where $K$ may be $\mathbb{R}$ or a multiple power of $\mathbb{R}$,} and by  $L^p(0,t;V)$ the Bochner space of functions defined on $(0,t)$ with values in a Banach space $V$, with $1\leq p \leq \infty$. If $K\equiv \mathbb{R}$, we simply write  $L^p(\Omega)$ and $W^{r,p}(\Omega)$. 
Moreover, when stating general results which are valid for both functions with scalar or vectorial or tensorial values, we write $f\in L^p$, $f\in W^{r,p}$, without specifying if $f$ is a function with scalar, vectorial or tensorial values.
For a normed space $X$, the associated norm is denoted by $||\cdot||_X$. In the case $p=2$, we use the notations $H^1:=W^{1,2}$ and $H^2:=W^{2,2}$, and we denote by $(\cdot,\cdot)$ and $||\cdot||$ the $L^2$ scalar product and induced norm between functions with scalar, vectorial or tensorial values. The dual space of a Banach space $Y$ is denoted by $Y'$. The duality pairing between $H^1(\Omega;K)$ and $\left(H^1(\Omega;K)\right)'$ is denoted by $<\cdot,\cdot>$. Moreover, we denote by $C^k(\Omega;K),C_c^k(\Omega;K)$ the spaces of continuously differentiable functions (respectively with compact support) up to order $k$ defined on $\Omega$ with values in a set $K$, and with {\rosso $C^k([0,t];V)$, $k\geq 0$}, the spaces of {\rosso continuously differentiable functions up to} order $k$ from $[0,t]$ to the space $V$. We finally introduce the spaces
\begin{align*}
& L_{\text{div}}^2(\Omega,\mathbb{R}^3):=\overline{\{\mathbf{u}\in C_c^{\infty}(\Omega,\mathbb{R}^3): \, \text{div}\mathbf{u}=0 \; \text{in} \; \Omega\}}^{||\cdot||_{L^2(\Omega;\mathbb{R}^3)}},\\
& H_{0,\text{div}}^1(\Omega,\mathbb{R}^3):=\overline{\{\mathbf{u}\in C_c^{\infty}(\Omega,\mathbb{R}^3): \, \text{div}\mathbf{u}=0 \; \text{in} \; \Omega\}}^{||\cdot||_{H_0^1(\Omega;\mathbb{R}^3)}}.
\end{align*}
\noindent
In the following, $C$ denotes a generic positive constant independent of the unknown variables, the discretization and the regularization parameters, the value of which might change from line to line; $C_1, C_2, \dots$ indicate generic positive constants whose particular value must be tracked through the calculations; $C(a,b,\dots)$ denotes a constant depending on the nonnegative parameters $a,b,\dots$.

\section{Model derivation} 
The phase field model coupled with finite viscoelasticity which we analyze in this paper is a particular case of a class of phase field models derived in our previous contribution \cite{agosti3}, obtained by considering a binary, saturated, closed and non reactive mixture of a solid elastic component and a liquid component, whose dynamics is driven by the microscopic interactions between its constituents and by their macroscopic visco-elastic behaviour. The mass and momentum balance
of the mixture were derived using a generalized form of the principles of virtual powers, giving constitutive assumptions satisfying the first and second law of thermodynamics in isothermal situations.

In particular, we consider the following form for the free energy $E$ of the system in Eulerian coordinates, associated to an arbitrary subregion of the mixture $R(t)\subset \Omega$ moving with the mixture:
\begin{align}
\label{eqn:free}
&E(\phi,\nabla \phi,\mathbf{F},\nabla \mathbf{F})=\int_{R(t)}e\left(\phi,\nabla \phi,\mathbf{F},\nabla \mathbf{F}\right)=\int_{R(t)}\left(\frac{1}{2}|\nabla \phi|^2+\psi(\phi)+w(\phi,\mathbf{F})+\frac{\lambda}{2}|\nabla \mathbf{F}|^2\right),
\end{align}
where $\phi$ is the volume concentration of the elastic phase, $\mathbf{F}$ is the deformation gradient associated to the motion of the elastic phase, $w(\phi,\mathbf{F})$ is the hyperelastic free energy density and $\psi(\phi)$ represents a bulk energy due to the mechanical interactions of the micro--components. A term proportional to
$\frac{1}{2}|\nabla \phi|^2+\psi(\phi)$ represents a surface energy contribution of the interface between the phases, expressed through a diffuse interface approach, while the term proportional to $|\nabla \mathbf{F}|^2$ is associated to elastic energy contributions from interfaces in the elastic material. Here, $\lambda>0$ is a regularization parameter, while the interface thickness parameter associated to the surface energy between the phases is taken to be equal to $1$ for ease of notation. We also consider the incompressibility constraint for the partial volume of the elastic phase.
The particular model is obtained from \cite[System (37)]{agosti3}, with the regularization parameters $\theta=\delta=0$ and employing the following simplifying assumptions:
\begin{itemize}
\item The momentum transfer in the mixture due to shear stresses in the liquid is negligible with respect to the momentum transfer between the solid and the liquid components;
\item The momentum exchange between the phases is induced by a Darcy-like flow of the liquid phase through the porous-permeable solid matrix associated to the solid phase, with constant permeability.
\end{itemize} 
With these assumptions, the PDE system takes the form
\begin{equation}
\label{eqn:system}
\begin{cases}
\displaystyle -\nu\Delta \mathbf{v}+\nabla s=\diver\left(\partial_{\mathbf{F}}w(\phi,\mathbf{F})\mathbf{F}^T\right)- \diver\left(\nabla \phi\otimes \nabla \phi+\lambda (\Delta \mathbf{F})\mathbf{F}^T+\lambda \nabla \mathbf{F} \odot \nabla \mathbf{F}\right), \\  \\ 
\displaystyle \diver\mathbf{v}=0, \\  \\ 
\displaystyle  \frac{\partial \mathbf{F}}{\partial t}+\left(\mathbf{v}\cdot \nabla\right)\mathbf{F}-(\nabla \mathbf{v})\mathbf{F}=0,\\ \\ 
\displaystyle \frac{\partial \phi}{\partial t}+\mathbf{v}\cdot \nabla \phi-\diver(b(\phi)\nabla \mu)=0,\\  \\
\displaystyle \mu=\psi'(\phi)- \Delta \phi+\partial_{\phi}w(\phi,\mathbf{F}),
\end{cases}
\end{equation}
valid in $\Omega_T$, endowed with the boundary conditions
\begin{equation}
\label{eqn:systembc}
b(\phi)\nabla \mu\cdot \mathbf{n}=\nabla \phi \cdot \mathbf{n}=0, \quad \left[\nabla \mathbf{F}\right]\mathbf{n}=\mathbf{0}, \quad \mathbf{v}=\mathbf{0}, 
\end{equation}
on $\partial \Omega \times [0,T]$, and with proper initial conditions. Here, $\nu>0$ is a viscosity parameter, $s$ is the solid pressure and $b(\phi)$ is a positive mobility function. System \eqref{eqn:system} is analogous to \cite[System (39)]{agosti3}, with \eqref{eqn:free} as the free energy of the system.

Using the relation
\begin{align}
\label{eqn:gradgradparts}
&-\diver\left(\nabla \phi \otimes \nabla \phi\right)-\lambda \diver\left(\nabla \mathbf{F} \odot \nabla \mathbf{F}\right)\\
&=-\nabla\left(\frac{1}{2}|\nabla \phi|^2+\psi(\phi)\right)+\mu\nabla \phi-\partial_{\phi}w\nabla \phi-\lambda \nabla\left(\frac{1}{2}|\nabla \mathbf{F}|^2\right)-\lambda \left(\nabla \mathbf{F}\right)\odot\Delta \mathbf{F}\\
& \notag =-\nabla\underbrace{\left(\frac{1}{2}|\nabla \phi|^2+\psi(\phi)+w(\phi,\mathbf{F})+ \frac{\lambda}{2}|\nabla \mathbf{F}|^2\right)}_{e\left(\phi,\nabla \phi,\mathbf{F},\nabla \mathbf{F}\right)}+\mu\nabla \phi+\left(\nabla \mathbf{F}\right)\odot\left(\partial_{\mathbf{F}}w-\lambda \Delta \mathbf{F}\right),
\end{align}
renaming $s\leftarrow s+e$ and adding a proper diffusive regularization, with regularization parameter $\gamma>0$, in the transport equation \eqref{eqn:system}$_3$, we get
\begin{equation}
\label{eqn:system2}
\begin{cases}
-\nu\Delta \mathbf{v}+\nabla s=\mu \nabla \phi+\diver\left(\left[\partial_{\mathbf{F}}w(\phi,\mathbf{F})-\lambda \Delta \mathbf{F}\right]\mathbf{F}^T\right)+\left(\nabla \mathbf{F}\right)\odot\left(\partial_{\mathbf{F}}w(\phi,\mathbf{F})-\lambda \Delta \mathbf{F}\right),\\ \\
\diver\mathbf{v}=0,\\ \\
\frac{\partial \mathbf{F}}{\partial t}+\left(\mathbf{v}\cdot \nabla\right)\mathbf{F}-(\nabla \mathbf{v})\mathbf{F}+\gamma \left(\partial_{\mathbf{F}}w(\phi,\mathbf{F})-\lambda \Delta \mathbf{F}\right)=0,\\ \\
\frac{\partial \phi}{\partial t}+\mathbf{v}\cdot \nabla \phi-\diver(b(\phi)\nabla \mu)=0,\\ \\
\mu=\psi'(\phi)- \Delta \phi+\partial_{\phi}w(\phi,\mathbf{F}).
\end{cases}
\end{equation}

We introduce the auxiliary variable
\[
\mathbf{M}:=\partial_{\mathbf{F}}w(\phi,\mathbf{F})-\lambda \Delta \mathbf{F},
\]
and rewrite system \eqref{eqn:system2} as (cf. also System \eqref{eqn:systemintro} in the Introduction)
\begin{equation}
\label{eqn:system3}
\begin{cases}
-\nu\Delta \mathbf{v}+\nabla s=\mu \nabla \phi+\diver\left(\mathbf{M}\mathbf{F}^T\right)+\nabla \mathbf{F}\odot \mathbf{M},\\ \\
\diver\mathbf{v}=0,\\ \\
\frac{\partial \mathbf{F}}{\partial t}+\left(\mathbf{v}\cdot \nabla\right)\mathbf{F}-(\nabla \mathbf{v})\mathbf{F}+\gamma \mathbf{M}=\mathbf{0},\\ \\
\mathbf{M}=\partial_{\mathbf{F}}w(\phi,\mathbf{F})-\lambda \Delta \mathbf{F},\\ \\
\frac{\partial \phi}{\partial t}+\mathbf{v}\cdot \nabla \phi-\diver(b(\phi)\nabla \mu)=0,\\ \\
\mu=\psi'(\phi)- \Delta \phi+\partial_{\phi}w(\phi,\mathbf{F}).
\end{cases}
\end{equation}
valid in $\Omega_T$, endowed with the boundary conditions
\begin{equation}
\label{eqn:systembc}
b(\phi)\nabla \mu\cdot \mathbf{n}=\nabla \phi \cdot \mathbf{n}=0, \quad \left[\nabla \mathbf{F}\right]\mathbf{n}=\mathbf{0}, \quad \mathbf{v}=\mathbf{0}, 
\end{equation}
on $\partial \Omega \times [0,T]$, and with initial conditions $\mathbf{F}(\mathbf{x},0)=\mathbf{F}_0(\mathbf{x})$, $\phi(\mathbf{x},0)=\phi_0(\mathbf{x})$ for $\mathbf{x}\in \Omega$.

In order to proceed, we make the following assumptions:
\begin{itemize}
\item[\bf{A0}] {$\Omega\subset \mathbb{R}^3$ is a bounded domain and the boundary $\partial \Omega$ is of class $C^{2}$};
\item[\bf{A1}] $b\in C^{0}(\mathbb{R})$ and there exist $b_0,b_1>0$ such that $b_0\leq b(r) \leq b_1,$ $\forall r \in \mathbb{R}$;
\item[\bf{A2}] $w\in C^{1}(\mathbb{R}\times \mathbb{R}^{3\times 3};\mathbb{R})$, and there exist $d_1\geq 0,d_2>0$ such that $-d_1 \leq w(r,\mathbf{T})\leq d_2\left(1+|\mathbf{T}|^p\right)$ for all $r \in \mathbb{R}, \mathbf{T}\in \mathbb{R}^{3\times 3}$,with $p\in [0,6)$. Moreover, there exist $d_3,d_4>0$ such that $\left|\partial_{\mathbf{T}}w(r,\mathbf{T})\right|\leq d_3\left(1+|\mathbf{T}|^{p-1}\right)$ and $\left|\partial_{r}w(r,\mathbf{T})\right|\leq d_4\left(1+|\mathbf{T}|^{p}\right)$ for all $r \in \mathbb{R}, \mathbf{T}\in \mathbb{R}^{3\times 3}$;
\item[\bf{A3}] $\psi \in C^{1}(\mathbb{R})$ and there exist $c_1>0,c_2\geq 0$ such that $|\psi'(r)|\leq c_1\left(|r|^l+1\right)$, $\psi(r)\geq -c_2$, $\forall r \in \mathbb{R}$, with $l\in \left[0,6+\frac{8}{\max(2,2p-6)}\right)$. Moreover, there exists a convex decomposition of $\psi=\psi_++\psi_-$, where $\psi_+$ is convex and $\psi_-$ is concave, such that $|\psi_-''(r)|\leq c_1\left(|r|^q+1\right)$, $\forall r \in \mathbb{R}$, with $q\in [0,4)$;
\item[\bf{A4}] {\rosso The initial data have the regularity $\mathbf{F}_0\in H^1(\Omega;\mathbb{R}^{3\times 3})$, $\phi_0\in H^1(\Omega)$.}
\end{itemize}
\begin{rem}
We observe that Assumption \textbf{A2} includes the situation in which $w(\phi,\mathbf{F})$ may degenerate in the variable $\phi$, i.e. $w(\bar{\phi},\mathbf{F})=0$ for specific values $\phi=\bar{\phi}$ (e.g. $\bar{\phi}=0$) and for all $ \mathbf{F}\in\mathbb{R}^{3\times 3}$. This situation could not be dealt with in the theoretical framework introduced in \cite{agosti3}. Moreover, we also observe that the growth law for $\psi$ in Assumption \textbf{A3} depends on the growth law for $w$ in Assumption \textbf{A2}.
\end{rem}
\begin{rem}
Assumption \textbf{A2} concerning the form of the elastic energy density $w(\phi,\mathbf{F})$ is quite general, and includes as particular cases many realistic elastic energy densities of nonlinear elastic materials, e.g. the generalized Ogden and the Mooney-Rivlin energy densities \cite{ogden}. These latter energy densities satisfy the polyconvexity and the coercivity properties, which are required by standard models in the theory of nonlinear elastic materials. We highlight that in our theoretical framework no polyconvexity and coercivity properties are generally required to obtain analytical results. Specifically, the Mooney--Rivlin energy density with phase-dependent elastic coefficients has the form
\begin{equation}
\label{eqn:freemr}
w(\phi,\mathbf{F})= \frac{f_1(\phi)}{2}(\mathbf{F}\colon \mathbf{F}-3)+\frac{f_2(\phi)}{2}\left(\left(\mathbf{F}\colon \mathbf{F}\right)^2-\mathbf{F}^T\mathbf{F}\colon \mathbf{F}^T\mathbf{F}-6\right),
\end{equation}
which, if $f_i\in C^{1}(\mathbb{R})$ with $0\leq f_i(r)\leq k_{1,i}$, $|f_i^{\prime}(r)|\leq k_{2,i}$, $\forall r\in \mathbb{R}$, $k_{1,i},k_{2,i}>0$ for $i=1,2$, satisfies Assumption \textbf{A2} with $p=4$. The Ogden energy density with phase-dependent elastic coefficients has the form
\begin{equation}
\label{eqn:freeog}
w(\phi,\mathbf{F})= \sum_{i=1}^{N}f_i(\phi)\left(\lambda_1^{p_i}+\lambda_2^{p_i}+\lambda_3^{p_i}-3\right),
\end{equation}
where $N\in \mathbb{N}$  and $p_1,...,p_N$ are material specific parameters and $\lambda_1,\lambda_2,\lambda_3$ are the principal stretches of deformation (i.e. $\lambda_1^2,\lambda_2^2,\lambda_3^2$ are the eigenvalues of $\mathbf{F}^T\mathbf{F}$). If $f_i\in C^{1}(\mathbb{R})$ with $0\leq f_i\leq k_{1,i}$, $|f_i^{\prime}(r)|\leq k_{2,i}$, $k_{1,i}, k_{2,i}>0$ and $0<p_i<6$ for and $i=1,\dots,N$, \eqref{eqn:freeog} satisfies Assumption \textbf{A2}. 
\end{rem}
We state here the main theorem of the present work concerning the existence of a global weak solution to \eqref{eqn:system3} in three space dimensions, which will be proved in the forthcoming sections.
\begin{thm}
\label{thm:3d}
Let the assumptions \textbf{A0}--\textbf{A4} be satisfied. Then, there exists a weak solution $(\mathbf{v}, \mathbf{F},  \mathbf{M}, \phi, \mu)$ of \eqref{eqn:system3}-\eqref{eqn:systembc}  with
\begin{align*}
&\mathbf{v}\in L^{2}\left(0,T;H_{0,\diver}^1\left(\Omega;\mathbb{R}^{3}\right)\right),\\
&\mathbf{F} \in L^{\infty}(0,T;H^1(\Omega;\mathbb{R}^{3\times 3}))\cap L^{2}(0,T;H^2(\Omega;\mathbb{R}^{3\times 3})),\\
&\mathbf{M} \in L^{2}(0,T;L^2(\Omega;\mathbb{R}^{3\times 3})),\\
&\partial_t \mathbf{F} \in L^{\frac{4}{3}-s}\left(0,T;L^2(\Omega;\mathbb{R}^{3\times 3})\right), \; \; s\in \left(0,\frac{1}{3}\right),\\
&\phi \in L^{\infty}(0,T;H^1(\Omega))\cap L^{\frac{8}{\max(2,2p-6)}}(0,T;H^2(\Omega)),\\
&\partial_t \phi \in L^{2}(0,T;\left(H^1(\Omega)\right)'),\\
&\mu \in L^{2}\left(0,T;H^1\left(\Omega\right)\right),
\end{align*}
such that
\begin{equation}
\label{eqn:continuous3d}
\begin{cases}
\displaystyle \nu \int_{\Omega}\nabla \mathbf{v}\colon \nabla \mathbf{u} =\int_{\Omega}\mu \nabla \phi \cdot \mathbf{u}-\int_{\Omega}\left(\mathbf{M}\mathbf{F}^T\right)\colon \nabla\mathbf{u}+\int_{\Omega}\left(\nabla \mathbf{F}\odot \mathbf{M}\right)\cdot \mathbf{u},\\
 \displaystyle \int_{\Omega}\partial_t \mathbf{F} \colon \boldsymbol{\Theta}+\int_{\Omega}\left({\mathbf{v}} \cdot \nabla \right)\mathbf{F}\colon \boldsymbol{\Theta}-\int_{\Omega}\left(\nabla {\mathbf{v}}\right)\mathbf{F}\colon \boldsymbol{\Theta}+\gamma \int_{\Omega}\mathbf{M}\colon \boldsymbol{\Theta}=0,\\
 \displaystyle \int_{\Omega}\mathbf{M}\colon \boldsymbol{\Pi}=\int_{\Omega}\partial_{\mathbf{F}}w(\phi,\mathbf{F})\colon \boldsymbol{\Pi}+\lambda \int_{\Omega}\nabla \mathbf{F} \mathbin{\tensorm} \nabla \boldsymbol{\Pi},\\
 \displaystyle <\partial_t\phi, q>+\int_{\Omega}\left(\mathbf{v} \cdot \nabla \phi \right) q+\int_{\Omega}b\left(\phi\right)\nabla \mu \cdot \nabla q =0,\\
\displaystyle \int_{\Omega}\mu r=\int_{\Omega}\nabla \phi \cdot \nabla r+\int_{\Omega}\psi'\left(\phi\right)r+\int_{\Omega} \partial_{\phi}w(\phi,\mathbf{F})r,
\end{cases}
\end{equation} 
for a.e. $t\in (0,T)$ and for all $\mathbf{u}\in H_{0,\diver}^1\left(\Omega;\mathbb{R}^{3}\right)$, $\boldsymbol{\Theta}\in L^2\left(\Omega;\mathbb{R}^{3\times 3}\right)$, $\boldsymbol{\Pi}\in H^1\left(\Omega;\mathbb{R}^{3\times 3}\right)$, $q, r \in H^1(\Omega)$, satisfying the initial conditions $\mathbf{F}(\cdot,0)=\mathbf{F}_0$ a.e. in $\Omega$ and $\phi(\cdot,0)=\phi_0$ a.e. in $\Omega$.
\end{thm}
In the following, we will introduce a proper truncation of the growth behavior of $w(\phi,\mathbf{F})$, depending on a truncation parameter $R>1$, and we will define a proper Faedo--Galerkin approximation of a truncated version of \eqref{eqn:system3}, proving the existence of a discrete solution and studying its convergence to a continuous weak solution, as the discretization parameter tends to zero, in three  space dimensions. We will then study the limit problem, at the continuous level, as the truncation parameter $R\to \infty$, thus removing the truncation operation and recovering an existence result for system \eqref{eqn:system3} associated to the full growth laws assumed in Assumption \textbf{A2}.
\\ \\
\noindent
Before proceeding, we state some preliminary results which will be used in the analysis.
\subsection{Preliminary lemmas}
We state here some Sobolev embedding and interpolation results which will be used in the following calculations. We start by recalling the Gagliardo-Nirenberg inequality {\rosso (see e.g. \cite{brezis})}.
\begin{lem}
\label{lem:gagliardoniremberg}
Let $\Omega \subset \mathbb{R}^3$ be a bounded domain with Lipschitz boundary and $f\in W^{m,r}\cap L^q$, $q\geq 1$, $r\leq \infty$, where $f$ can be a function with scalar, vectorial or tensorial values. For any integer $j$ with $0 \leq j < m$, suppose there is $\alpha \in \mathbb{R}$ such that
\[
j-\frac{3}{k}=\left(m-\frac{3}{r}\right)\alpha+(1-\alpha)\left(-\frac{3}{q}\right), \quad \frac{j}{m}\leq \alpha \leq 1.
\]
Then, there exists a positive constant $C$ depending on $\Omega$, m, j, q, r, and $\alpha$ such that
\begin{equation}
\label{eqn:gagliardoniremberg}
||D^jf||_{L^k}\leq C||f||_{W^{m,r}}^{\alpha}||f||_{L^q}^{1-\alpha}.
\end{equation}
\end{lem}
We also state the following interpolation results.
\begin{lem}
\label{lem:interpolation2d3d}
Let $\Omega \subset \mathbb{R}^3$ be a bounded domain with Lipschitz boundary, $p\in [1,6)$ and $f \in L^{\infty}(0,T;L^2)\cap L^2(0,T;W^{1,\frac{6}{p-1}})$, where $f(\mathbf{x},t)$, with $t\in (0,T)$, $\mathbf{x}\in \Omega$, may be a scalar, a vector or a tensor. 
Then, there exists a positive constant $C$ depending on $\Omega$ such that
{\rosso
\begin{equation}
\label{eqn:gninterpolation1}
\int_0^T||f||_{L^{2+h}}^{\frac{2(2+h)(6-p)}{3h}}\leq C\int_0^T||f||_{L^{2}}^{\frac{2[(6-p)(2+h)-3h]}{3h}}||f||_{W^{1,\frac{6}{p-1}}}^2, \; \begin{cases} h\geq 0 \; \text{if} \; 1\leq p\leq 3, \\ h\in \left[0,\frac{2(6-p)}{p-3}\right]\; \text{if} \; 3<p<6.\end{cases}
\end{equation}
Let moreover $p\in [1,6)$ and $f \in L^{\infty}(0,T;L^6)\cap L^{\frac{18-2p}{3}}(0,T;W^{1,\frac{18-2p}{3}})$, where $f(\mathbf{x},t)$, with $t\in (0,T)$, $\mathbf{x}\in \Omega$, may be a scalar, a vector or a tensor. 
Then, there exists a positive constant $C$ depending on $\Omega$ such that
\begin{equation}
\label{eqn:gninterpolation2}
\int_0^T||f||_{L^{6+h}}^{\frac{2(6+h)(6-p)}{h}}\leq C\int_0^T||f||_{L^{6}}^{\frac{2[(18-3p)(6+h)-h(9-p)]}{3h}}||f||_{W^{1,\frac{18-2p}{3}}}^{\frac{18-2p}{3}},\; \begin{cases} h\geq 0 \; \text{if} \; 1\leq p\leq \frac{9}{2}, \\ h\in \left[0,\frac{18(6-p)}{2p-9}\right]\; \text{if} \; \frac{9}{2}<p<6.\end{cases}
\end{equation}
Let finally $f \in L^{\infty}(0,T;L^6)\cap L^{s}(0,T;W^{1,6})$, with $s\geq 1$, where $f(\mathbf{x},t)$, with $t\in (0,T)$, $\mathbf{x}\in \Omega$, may be a scalar, a vector or a tensor. 
Then, there exists a positive constant $C$ depending on $\Omega$ such that
\begin{equation}
\label{eqn:gninterpolation3}
\int_0^T||f||_{L^{6+h}}^{\frac{2s(6+h)}{h}}\leq C\int_0^T||f||_{L^{6}}^{\frac{(12+h)s}{h}}||f||_{W^{1,6}}^{s},\; h\geq 0.
\end{equation} 
}%
\end{lem}
We observe that \eqref{eqn:gninterpolation1}, \eqref{eqn:gninterpolation2} and \eqref{eqn:gninterpolation3} are consequences of the Gagliardo--Nirenberg inequality \eqref{eqn:gagliardoniremberg} with $j=0$, $m=1$, $k=2+h$, $r=\frac{6}{p-1}$, $q=2$ (for \eqref{eqn:gninterpolation1}), $j=0$, $m=1$, $k=6+h$, $r=\frac{18-2p}{3}$, $q=6$ (for \eqref{eqn:gninterpolation2}) and $j=0$, $m=1$, $k=6+h$, $r=6$, $q=6$ (for \eqref{eqn:gninterpolation3}).
\noindent
We moreover recall the following interpolation inequality.
\begin{lem}
\label{lem:interpolationlp}
Let $\Omega \subset \mathbb{R}^3$ be a bounded domain and $f \in L^q$, $q\geq 1$, where $f$ can be a function with scalar, vectorial or tensorial values. Let also $s\leq r \leq q$. Then, there exists a positive constant $C$ depending on $\Omega$ such that
\begin{equation}
\label{eqn:lpinterpolation}
\int_{\Omega}|f|^r\leq \left(\int_{\Omega}|f|^s\right)^{\frac{q-r}{q-s}}\left(\int_{\Omega}|f|^q\right)^{\frac{r-s}{q-s}}.
\end{equation}
\end{lem}
We finally state the Agmon type inequality in three space dimensions {\rosso (see e.g. \cite{agmon})} which will be used in the following calculations.
\begin{lem}
\label{lem:agmon2d3d}
Let $\Omega \subset \mathbb{R}^3$, be a bounded domain with Lipschitz boundary and $f \in H^2(\Omega)$. 
Then, there exists a positive constant $C$ depending on $\Omega$ such that
\begin{equation}
\label{eqn:agmon3d}
||f||_{L^{\infty}(\Omega)}\leq C||f||_{H^1(\Omega)}^{\frac{1}{2}}||f|||_{H^2(\Omega)}^{\frac{1}{2}}.
\end{equation}
\end{lem}
\section{Existence result of a global weak solution}
We first need to obtain an existence result for a truncated system.

\noindent
Let us introduce the smooth step function $g\in C^1(\mathbb{R}^+)$ with the properties
\begin{align}
\begin{cases}
\label{eqn:propg}
 0\leq g(\cdot)\leq 1;\\
 g(r)\equiv 1 \; \text{for} \; r\leq 1; \quad  g(r)\equiv 0 \; \text{for} \; r\geq 2;\\
 |g^{\prime}(r)|\leq C_g \; \forall r\geq 0.
\end{cases}
\end{align}
Given $R>1$, we define the smooth truncation function
\begin{equation}
\label{eqn:gr}
g_R(r):=g\left(\frac{r}{R}\right)+\left(1-g\left(\frac{r}{R}\right)\right)r^{\min(0,4-p)},
\end{equation}
where $p$ is defined in Assumption \textbf{A2}.
We have that
\begin{equation}
\label{eqn:grp}
g^{\prime}_R(r)=\frac{1}{R}g^{\prime}\left(\frac{r}{R}\right)\left(1-r^{\min(0,4-p)}\right)+\min(0,4-p)\left(1-g\left(\frac{r}{R}\right)\right)r^{\min(-1,3-p)}.
\end{equation}
We observe that
\begin{equation}
\label{eqn:grbound}
0\leq g_R(r)\leq
\begin{cases}
1 \quad \text{for} \; r\leq R;\\
1+\frac{1}{R^{\max(0,p-4)}} \quad \text{for} \; R\leq r \leq 2R;\\
 \frac{1}{r^{\max(0,p-4)}}\leq \frac{1}{(2R)^{\max(0,p-4)}}  \quad \text{for}\, r \geq 2R,
\end{cases}
\end{equation}
and
\begin{equation}
\label{eqn:grpbound}
|g^{\prime}_R(r)|\leq
\begin{cases}
0 \quad \text{for} \; r\leq R,\\
\frac{2C_g}{2R}\left(1-\frac{1}{(2R)^{\max(0,p-4)}}\right)+\frac{2^{\max(1,p-3)}\min(0,4-p)}{(2R)^{\max(1,p-3)}} \quad \text{for} \; R\leq r \leq 2R,\\
0  \quad \text{for}\, r \geq 2R.
\end{cases}
\end{equation}
We introduce the truncated elastic energy density
\begin{equation}
\label{eqn:wr}
w_R(\phi,\mathbf{F}):=g_R(|\mathbf{F}|)w(\phi,\mathbf{F}).
\end{equation}
Thanks to \eqref{eqn:grbound} we have that
\begin{equation}
\label{eqn:wrbound}
-d_1\leq w_R(\phi,\mathbf{F})\leq
\begin{cases}
C(1+|2R|^p) \quad \text{for} \; |\mathbf{F}|\leq 2R;\\ \\
 C|\mathbf{F}|^{\min(0,4-p)}\left(1+|\mathbf{F}|^p\right)\leq C\frac{1}{(2R)^{\max(0,p-4)}}+C|\mathbf{F}|^{p+\min(0,4-p)}\\
 \leq C(1+|\mathbf{F}|^{p+\min(0,4-p)})  \quad \text{for}\, |\mathbf{F}| \geq 2R.
\end{cases}
\end{equation}
Given \eqref{eqn:grpbound} and the relation
\begin{equation}
\label{eqn:dfwr}
\partial_{\mathbf{F}}w_R(\phi,\mathbf{F})=g^{\prime}_R(|\mathbf{F}|)\frac{\mathbf{F}}{|\mathbf{F}|}w(\phi,\mathbf{F})+g_R(|\mathbf{F}|)\partial_{\mathbf{F}}w(\phi,\mathbf{F}),
\end{equation} 
we have that
\begin{equation}
\label{eqn:wrpbound}
|\partial_{\mathbf{F}}w_R(\phi,\mathbf{F})|\leq
\begin{cases}
C(1+|2R|^{p-1}) \quad \text{for} \; |\mathbf{F}|\leq 2R;\\ \\
 C|\mathbf{F}|^{\min(0,4-p)}\left(1+|\mathbf{F}|^{p-1}\right)\leq C\frac{1}{(2R)^{\max(0,p-4)}}+C|\mathbf{F}|^{p-1+\min(0,4-p)}\\
 \leq C(1+|\mathbf{F}|^{p+\min(-1,3-p)})  \quad \text{for}\, |\mathbf{F}| \geq 2R.
\end{cases}
\end{equation}
Thanks to \eqref{eqn:propg}, \eqref{eqn:wrbound}, \eqref{eqn:wrpbound} and the fact that $\partial_{\phi}w_R(\phi,\mathbf{F})=g_R(|\mathbf{F}|)\partial_{\phi}w(\phi,\mathbf{F})$, as a consequence of Assumption \textbf{A2} we obtain the following property for $w_R$:
\begin{itemize}
\item[\bf{A2Bis}] $w_R\in C^{1}(\mathbb{R}\times \mathbb{R}^{3\times 3};\mathbb{R})$, and there exist $d_1\geq 0,d_2>0$ such that $-d_1\leq w_R(r,\mathbf{T})\leq d_2\left(1+|\mathbf{T}|^p\right)$ for all $r \in \mathbb{R}, \mathbf{T}\in \mathbb{R}^{3\times 3}$, with $p\in [0,4]$. Moreover, there exist $d_3,d_4>0$ such that $\left|\partial_{\mathbf{T}}w(r,\mathbf{T})\right|\leq d_3\left(1+|\mathbf{T}|^{p-1}\right)$ and $\left|\partial_{r}w(r,\mathbf{T})\right|\leq d_4\left(1+|\mathbf{T}|^{p}\right)$ for all $r \in \mathbb{R}, \mathbf{T}\in \mathbb{R}^{3\times 3}$.
\end{itemize}
Finally, we define the truncated system, depending on the finite parameter $R$,
\begin{equation}
\label{eqn:system3r}
\begin{cases}
-\nu\Delta \mathbf{v}_R+\nabla q_R=\mu_R \nabla \phi_R+\diver\left(\mathbf{M}_R\mathbf{F}_R^T\right)+\nabla \mathbf{F}_R\odot \mathbf{M}_R,\\ \\
\diver\mathbf{v}_R=0,\\ \\
\frac{\partial \mathbf{F}_R}{\partial t}+\left(\mathbf{v}_R\cdot \nabla\right)\mathbf{F}_R-(\nabla \mathbf{v}_R)\mathbf{F}_R+\gamma \mathbf{M}_R=0,\\ \\
\mathbf{M}_R=\partial_{\mathbf{F}_R}w_R(\phi_R,\mathbf{F}_R)-\lambda \Delta \mathbf{F}_R,\\ \\
\frac{\partial \phi_R}{\partial t}+\mathbf{v}_R\cdot \nabla \phi_R-\diver(b(\phi_R)\nabla \mu_R)=0,\\ \\
\mu_R=\psi'(\phi_R)- \Delta \phi_R+\partial_{\phi_R}w_R(\phi_R,\mathbf{F}_R),
\end{cases}
\end{equation}
valid in $\Omega_T$, endowed with the same boundary conditions as \eqref{eqn:systembc}.
In the following, for ease of notation we will avoid to report the $R$ index for the solutions of system \eqref{eqn:system3r}.  

\subsection{Faedo--Galerkin approximation scheme}
\label{sec:approx}
We define the finite dimensional spaces which will be used to formulate the Galerkin ansatz to approximate the solutions of system \eqref{eqn:system3r}.
Let $\{\boldsymbol{\eta_i}\}_{i\in \mathbb{N}}$ be the eigenfunctions of the Stokes operator with homogeneous Dirichlet boundary conditions, i.e.
\[
P_L(-\Delta)\boldsymbol{\eta_i}=\beta_i \boldsymbol{\eta_i} \quad \text{in} \; \Omega, \quad \boldsymbol{\eta_i}=\mathbf{0} \quad \text{on} \; \partial \Omega,
\]
where $P_L:L^2(\Omega;\mathbb{R}^3)\to L_{\text{div}}^2(\Omega;\mathbb{R}^3)$ is the Leray projection operator, with $0<\beta_0\leq \beta_1 \leq \dots \leq \beta_m\to \infty$. The sequence $\{\boldsymbol{\eta_i}\}_{i\in \mathbb{N}}$ can be chosen as an orthonormal basis in $L_{\diver}^2(\Omega;\mathbb{R}^3)$ and an orthogonal basis in $H_{0,\diver}^1(\Omega;\mathbb{R}^3)$, and, thanks to Assumption \textbf{A0}, we have that $\{\boldsymbol{\eta_i}\}_{i\in \mathbb{N}}\subset H^{2}(\Omega;\mathbb{R}^3)$. We introduce the projection operator
\[
P_m^S:H_{0,\diver}^1(\Omega;\mathbb{R}^3)\to \text{span}\{\boldsymbol{\eta_0},\boldsymbol{\eta_1},\dots,\boldsymbol{\eta_m}\}.
\]
Let $\{\xi_i\}_{i\in \mathbb{N}}$ be the eigenfunctions of the Laplace operator with homogeneous Neumann boundary conditions, i.e.
\[
-\Delta \xi_i=\alpha_i \xi_i \quad \text{in} \; \Omega, \quad \nabla \xi_i \cdot \mathbf{n}=\mathbf{0} \quad \text{on} \; \partial \Omega,
\]
with $0=\alpha_0< \alpha_1 \leq \dots \leq \alpha_m\to \infty$. The sequence $\{\xi_i\}_{i\in \mathbb{N}}$ can be chosen as an orthonormal basis in $L^2(\Omega)$ and an orthogonal basis in $H^1(\Omega)$, and, thanks to Assumption \textbf{A0}, $\{\xi_i\}_{i\in \mathbb{N}}\subset H^{2}(\Omega)$. Without loss of generality, we assume $\alpha_0=0$. We introduce the projection operator
\[
P_m^L:H^1(\Omega)\to \text{span}\{\xi_0,\xi_1,\dots,\xi_m\}.
\]
Let moreover $\{\boldsymbol{\Sigma}_i\}_{i\in \mathbb{N}}$ be the tensor eigenfunctions of the Laplace operator with homogeneous Neumann boundary conditions, i.e.
\[
-\Delta \boldsymbol{\Sigma}_i=\gamma_i \boldsymbol{\Sigma}_i \quad \text{in} \; \Omega, \quad \nabla \boldsymbol{\Sigma}_i  \mathbf{n}=\mathbf{0} \quad \text{on} \; \partial \Omega,
\]
with $0=\gamma_0=\gamma_1=\dots=\gamma_8< \gamma_9 \leq \dots \leq \gamma_m\to \infty$. Note that the eigenvalues $\gamma_0=\dots=\gamma_8=0$ are associated to the eigentensors $\boldsymbol{\Sigma}_0,\dots,\boldsymbol{\Sigma}_8$ which are proportional to the tensors $\mathbf{e}_{ij}$, $i,j=1,2,3$.
The sequence $\{\boldsymbol{\Sigma}_i\}_{i\in \mathbb{N}}$ can be chosen as an orthonormal basis in $L^2(\Omega;\mathbb{R}^{3\times 3}))$ and an orthogonal basis in $H^1(\Omega;\mathbb{R}^{3\times 3}))$, and, thanks to Assumption \textbf{A0}, $\{\boldsymbol{\Sigma}_i\}_{i\in \mathbb{N}}\subset H^{2}(\Omega;\mathbb{R}^{3\times 3})$. We introduce the projection operator
\[
P_m^{L,\Sigma}:H^1(\Omega;\mathbb{R}^{3\times 3})\to \text{span}\{\boldsymbol{\Sigma}_0,\boldsymbol{\Sigma}_1,\dots,\boldsymbol{\Sigma}_m\}.
\]
We make the Galerkin ansatz $\mathbf{v}_m=\sum_i d_i^m(t)\boldsymbol{\eta_i}(\mathbf{x})$, $\mathbf{F}_m=\sum_i f_i^m(t)\boldsymbol{\Sigma_i}(\mathbf{x})$, $\mathbf{M}_m=\sum_i g_i^m(t)\boldsymbol{\Sigma_i}(\mathbf{x})$, $\phi_m=\sum_i a_i^m(t)\xi_i(\mathbf{x})$, $\mu_m=\sum_i c_i^m(t)\xi_i(\mathbf{x})$ to approximate the solutions $\mathbf{v},\mathbf{F},\mathbf{M},\phi,\mu$ of system \eqref{eqn:system3r}, and project  the equation for $\mathbf{v}_m$ onto $\text{span}\left\{\boldsymbol{\eta}_0,\boldsymbol{\eta}_1,\dots,\boldsymbol{\eta}_m\right\}$, the equations for $\mathbf{F}_m$ and $\mathbf{M}_m$ onto
$\text{span}\{\boldsymbol{\Sigma_0},\boldsymbol{\Sigma_1},\dots,\boldsymbol{\Sigma_m}\}$
and the equations for $\phi_m$ and $\mu_m$ onto $\text{span}\left\{\xi_0,\xi_1, \dots, \xi_m\right\}$, obtaining the following Galerkin approximation of system \eqref{eqn:system3r}:
\begin{equation}
\label{eqn:system3gproj}
\begin{cases}
 \displaystyle \nu \int_{\Omega}\nabla \mathbf{v}_m\colon \nabla \boldsymbol{\eta}_i =\int_{\Omega}\mu_m \nabla \phi_m \cdot \boldsymbol{\eta}_i-\int_{\Omega}\left(\mathbf{M}_m\mathbf{F}_m^T\right)\colon \nabla\boldsymbol{\eta}_i+\int_{\Omega}\left(\nabla \mathbf{F}_m\odot \mathbf{M}_m\right)\cdot \boldsymbol{\eta}_i,\\
 \displaystyle \int_{\Omega}\partial_t \mathbf{F}_m \colon \boldsymbol{\Sigma}_i+\int_{\Omega}\left({\mathbf{v}}_m \cdot \nabla \right)\mathbf{F}_m\colon \boldsymbol{\Sigma}_i-\int_{\Omega}\left(\nabla {\mathbf{v}}_m\right)\mathbf{F}_m\colon \boldsymbol{\Sigma}_i+\gamma \int_{\Omega}\mathbf{M}_m\colon \boldsymbol{\Sigma}_i=0,\\
 \displaystyle \int_{\Omega}\mathbf{M}_m\colon \boldsymbol{\Sigma}_i=\int_{\Omega}\partial_{\mathbf{F}_m}w_R(\phi_m,\mathbf{F}_m)\colon \boldsymbol{\Sigma}_i+\lambda \int_{\Omega}\nabla \mathbf{F}_m \mathbin{\tensorm} \nabla \boldsymbol{\Sigma}_i,\\
 \displaystyle \int_{\Omega}\partial_t\phi_m \xi_i+\int_{\Omega}\left(\mathbf{v}_m \cdot \nabla \phi_m\right) \xi_i+\int_{\Omega}b\left(\phi_m\right)\nabla \mu_m \cdot \nabla \xi_i =0,\\
\displaystyle \int_{\Omega}\mu_m \xi_i=\int_{\Omega}\nabla \phi_m\cdot \nabla \xi_i+\int_{\Omega}\psi'\left(\phi_m\right)\xi_i+\int_{\Omega} \partial_{\phi_m}w_R(\phi_m,\mathbf{F}_m)\xi_i,
\end{cases}
\end{equation} 
{\rosso in $[0,t]$}, with $0<t\leq T$, for $i=0, \dots,m$ and with initial conditions 
\begin{equation}
    \label{eqn:system3projic}
    \phi_m(\mathbf{x},0)=P_m^L(\phi_0), \quad \mathbf{F}_m(\mathbf{x},0)=P_m^{L,\Sigma}(\mathbf{F}_0).
\end{equation}

System \eqref{eqn:system3gproj} defines a collection of initial value problems for a system of coupled ODEs
\begin{equation}
\label{eqn:odegalerkin}
\begin{cases}
 \displaystyle \nu \beta_id_i^m&=\sum_{l,s}\left[\int_{\Omega}\xi_l \nabla \xi_s \cdot \boldsymbol{\eta}_i\right]b_l^ma_s^m-\sum_{l,s}\left[\int_{\Omega}\left(\boldsymbol{\Sigma}_l\boldsymbol{\Sigma}_s^T\right)\colon \nabla\boldsymbol{\eta}_i\right]g_l^mf_s^m\\
 &+\sum_{l,s}\left[\int_{\Omega}\left(\nabla \boldsymbol{\Sigma}_l\odot \boldsymbol{\Sigma}_s\right)\cdot \boldsymbol{\eta}_i\right]f_l^mg_s^m,\\
 \displaystyle \frac{d}{dt}f_i^m&=\sum_{l,s}\biggl(-\int_{\Omega}\left(\boldsymbol{\eta}_l \cdot \nabla \right)\boldsymbol{\Sigma}_s\colon \boldsymbol{\Sigma}_i+\int_{\Omega}\left(\nabla \boldsymbol{\eta}_l\right)\boldsymbol{\Sigma}_s\colon \boldsymbol{\Sigma}_i\biggr)d_l^mf_s^m-\gamma g_i^m,\\
 \displaystyle g_i^m&=\int_{\Omega}\partial_{\mathbf{F}_m}w_R\left(\sum_{s} a_s^m\xi_s,\sum_l f_l^m\boldsymbol{\Sigma}_l\right)\colon \boldsymbol{\Sigma}_i+\lambda \gamma_i f_i^m,\\
 \displaystyle \frac{d}{dt}a_i^m&=-\sum_{l,s}\left[\int_{\Omega}\left(\boldsymbol{\eta}_l \cdot \nabla \xi_s\right) \cdot \xi_i\right]d_l^ma_s^m-
\sum_{l}\left[\int_{\Omega}b\left(\sum_{s} a_s^m\xi_s\right)\nabla \xi_l \cdot \nabla \xi_i \right]c_l^m,\\
\displaystyle c_i^m&=\alpha_ia_i^m+\int_{\Omega}\psi'\left(\sum_{s} a_s^m\xi_s\right)\xi_i+\int_{\Omega}\partial_{\phi_m}w_R\left(\sum_{s} a_s^m\xi_s,\sum_l f_l^m\boldsymbol{\Sigma}_l\right)\xi_i,\\
\displaystyle a_i^m(0)&=\left(\phi_0,\xi_i\right), \; f_i^m(0)=\left(\mathbf{F}_0,\boldsymbol{\Sigma_i}\right), \quad i=0, \dots, m.
\end{cases}
\end{equation} 
Due to the Assumptions \textbf{A1, A2Bis, A3} on the regularity of the functions $b,\psi,w_R$ and the regularity in space of the functions $\xi_i, \boldsymbol{\eta}_i, \boldsymbol{\Sigma}_i,$, the right hand side of \eqref{eqn:odegalerkin} depends continuously on the independent variables and we can apply the Peano existence theorem to infer that there exist a sufficiently small $t_1$ with $0<t_1\leq T$ and a local solution $(d_i^m,f_i^m,g_i^m,a_i^m,c_i^m)$ of \eqref{eqn:odegalerkin}, for $i=0, \dots, m$. 
\subsection{A priori estimates} 
We now deduce a priori estimates, uniform in the discretization parameter $m$, for the solutions of system \eqref{eqn:system3gproj}, which can be rewritten, combining the equations over $i=0, \dots, m$, as
\begin{equation}
\label{eqn:system3galerkin}
\begin{cases}
 \displaystyle \nu \int_{\Omega}\nabla \mathbf{v}_m\colon \nabla \boldsymbol{\eta} =\int_{\Omega}\mu_m \nabla \phi_m \cdot \boldsymbol{\eta}-\int_{\Omega}\left(\mathbf{M}_m\mathbf{F}_m^T\right)\colon \nabla\boldsymbol{\eta}+\int_{\Omega}\left(\nabla \mathbf{F}_m\odot \mathbf{M}_m\right)\cdot \boldsymbol{\eta},\\
 \displaystyle \int_{\Omega}\partial_t \mathbf{F}_m \colon \boldsymbol{\Sigma}+\int_{\Omega}\left({\mathbf{v}}_m \cdot \nabla \right)\mathbf{F}_m\colon \boldsymbol{\Sigma}-\int_{\Omega}\left(\nabla {\mathbf{v}}_m\right)\mathbf{F}_m\colon \boldsymbol{\Sigma}+\gamma \int_{\Omega}\mathbf{M}_m\colon \boldsymbol{\Sigma}=0,\\
 \displaystyle \int_{\Omega}\mathbf{M}_m\colon \boldsymbol{\Gamma}=\int_{\Omega}\partial_{\mathbf{F}_m}w_R(\phi_m,\mathbf{F}_m)\colon \boldsymbol{\Gamma}+\lambda \int_{\Omega}\nabla \mathbf{F}_m \mathbin{\tensorm} \nabla \boldsymbol{\Gamma},\\
 \displaystyle \int_{\Omega}\partial_t\phi_m \xi+\int_{\Omega}\left(\mathbf{v}_m \cdot \nabla \phi_m\right) \xi+\int_{\Omega}b\left(\phi_m\right)\nabla \mu_m \cdot \nabla \xi =0,\\
\displaystyle \int_{\Omega}\mu_m \chi=\int_{\Omega}\nabla \phi_m\cdot \nabla \chi+\int_{\Omega}\psi'\left(\phi_m\right)\chi+\int_{\Omega} \partial_{\phi_m}w_R(\phi_m,\mathbf{F}_m)\chi,
\end{cases}
\end{equation} 
for a.e. $t \in [0,t_1]$, with $\boldsymbol{\eta} \in \text{span}\left\{\boldsymbol{\eta}_0, \boldsymbol{\eta}_1, \dots, \boldsymbol{\eta}_m\right\}$, $\boldsymbol{\Sigma}, \boldsymbol{\Gamma} \in \text{span}\{\boldsymbol{\Sigma_0},\boldsymbol{\Sigma_1},\dots,\boldsymbol{\Sigma_m}\}$ and
$\xi, \chi \in \text{span}\left\{\xi_0, \xi_1, \dots, \xi_m\right\}$, with initial conditions defined in \eqref{eqn:system3projic}. We take $\boldsymbol{\eta}=\mathbf{v}_m$, $\boldsymbol{\Sigma}=\mathbf{M}_m$, $\boldsymbol{\Gamma}=-\partial_t \mathbf{F}_m$, $\xi=\mu_m$, {\rosso $\chi=-\partial_t \phi_m$} in \eqref{eqn:system3galerkin}, sum all the equations and integrate in time between $0$ and $t\in [0,t_1]$. We get, for any $t\in [0,t_1]$,
\begin{align}
\label{eqn:apriori1}
&\notag \nu\int_0^{t}\int_{\Omega}\left|\nabla \mathbf{v}_m\right|^2+\gamma \int_0^{t}\int_{\Omega}\left|\mathbf{M}_m\right|^2+\int_0^{t}\int_{\Omega}b(\phi_m)\left|\nabla \mu_m\right|^2\\
& \notag +\int_{\Omega}\left(\frac{1}{2}|\nabla \phi_m(t)|^2+\psi(\phi_m(t))+w_R(\phi_m(t),\mathbf{F}_m(t))+\frac{\lambda}{2}|\nabla \mathbf{F}_m(t)|^2\right)\\
& = \int_{\Omega}\left(\frac{1}{2}|\nabla \phi_m(0)|^2+\psi(\phi_m(0))+w_R(\phi_m(0),\mathbf{F}_m(0))+\frac{\lambda}{2}|\nabla \mathbf{F}_m(0)|^2\right).
\end{align}
Hence, \eqref{eqn:apriori1} becomes
\begin{align}
\label{eqn:apriori}
& \notag \sup_{t\in [0,t_1]}\left(\frac{1}{2}|\nabla \phi_m|^2+\frac{\lambda}{2}|\nabla \mathbf{F}_m|^2\right)+\nu\int_0^{t_1}\int_{\Omega}\left|\nabla \mathbf{v}_m\right|^2+\gamma \int_0^{t_1}\int_{\Omega}\left|\mathbf{M}_m\right|^2\\
&+\int_0^{t_1}\int_{\Omega}b(\phi_m)\left|\nabla \mu_m\right|^2 \leq C(\mathbf{F}_0,\phi_0)+C,
\end{align}
where we used Assumptions \textbf{A2Bis}, \textbf{A3} and {\rosso \textbf{A4}}. The constants in the right hand side of \eqref{eqn:apriori} depends only on the initial data and on the domain $\Omega$ and not on the discretization parameter $m$ and on the truncation parameter $R$. Thanks to the a priori estimate \eqref{eqn:apriori}, we may extend {\rosso by continuity} the local solution of system \eqref{eqn:system3galerkin} to the interval $[0,T]$.

From the Poincar\'e inequality and from \eqref{eqn:apriori}, we have that
{\rosso
\begin{equation}
\label{eqn:vm3d}
\mathbf{v}_m\; \text{is uniformly bounded in} \; L^2(0,T;H_{0,\diver}^1(\Omega;\mathbb{R}^3))\hookrightarrow L^2(0,T;L_{\diver}^6(\Omega;\mathbb{R}^3)).
\end{equation}
Moreover, from \eqref{eqn:apriori} we have that
\begin{equation}
\label{eqn:Mm3d}
\mathbf{M}_m \; \text{is uniformly bounded in} \; L^{2}(0,T;L^2(\Omega;\mathbb{R}^{3\times 3})).
\end{equation}
We take $\boldsymbol{\Sigma}=\mathbf{e}_{l,r}$, $l,r=1,2,3$, which are proportional to the eigentensors associated to $\gamma_0,\dots,\gamma_8$. Hence, in \eqref{eqn:system3galerkin}$_2$, integrating by parts in the second and third terms, we obtain that
\begin{align}
\label{eqn:fmean}
& \displaystyle \notag \partial_t\int_{\Omega}\mathbf{F}_{m,lr}+\int_{\partial \Omega}\mathbf{F}_{m,lr}\mathbf{v}_m\cdot \mathbf{n}-\int_{\Omega}\mathbf{F}_{m,lr}\diver \mathbf{v}_m-\int_{\partial \Omega}\mathbf{v}_{m,l}\mathbf{F}_m^{[r]}\cdot \mathbf{n}\\
& \displaystyle +\int_{\Omega}\mathbf{v}_{m,l}\diver \mathbf{F}_m^{[r]}+\gamma\int_{\Omega}\mathbf{M}_{m,lr}=\mathbf{0}.
\end{align}
Integrating in time the latter equations over the interval $[0,T]$, using the facts that $\diver\mathbf{v}_m=0$, $\mathbf{v}_m=\mathbf{0}$ on $\partial \Omega$, the Cauchy-Schwarz and Young inequalities, \eqref{eqn:apriori} and Assumption {\rosso \textbf{A4}}, we obtain that
\begin{equation}
\label{eqn:massFmdisc}
\sup_{t\in [0,T]}\left|\left(\int_{\Omega}\mathbf{F}_{m,lr}\right)(t)\right|\leq C(\mathbf{F}_0,\phi_0,T), \quad \forall l,r=1,2,3.
\end{equation}
Hence, from \eqref{eqn:apriori}, \eqref{eqn:massFmdisc} and the Poincar\'e--Wirtinger inequality we deduce that
and
\begin{equation}
\label{eqn:Fm3d}
\mathbf{F}_m \; \text{is uniformly bounded in} \; L^{\infty}(0,T;H^1(\Omega;\mathbb{R}^{3\times 3}))\hookrightarrow L^{\infty}(0,T;L^6(\Omega;\mathbb{R}^{3\times 3})).
\end{equation}
Next, taking $\xi=1$ (which is a multiple of $\xi_0$) in \eqref{eqn:system3galerkin}$_4$, we obtain, using the facts that $\diver\mathbf{v}_m=0$, $\mathbf{v}_m=\mathbf{0}$ on $\partial \Omega$ and integrating by parts in the second term, that
\begin{equation}
\label{eqn:massphimdisc}
\left(\partial_t\phi_m,1\right)=0.
\end{equation}
Hence, from \eqref{eqn:apriori}, \eqref{eqn:massphimdisc} and the Poincar\'e--Wirtinger inequality we deduce that
\begin{equation}
\label{eqn:phim3d}
\phi_m\; \text{is uniformly bounded in} \; L^{\infty}(0,T;H^1(\Omega))\hookrightarrow L^{\infty}(0,T;L^6(\Omega)).
\end{equation}
\subsection{Higher order regularity for $\mathbf{F}_m$ and $\phi_m$}
We observe that, from \eqref{eqn:Fm3d} and \textbf{A2Bis}, $\mathbf{F}_m\in L^{\infty}(0,T;L^6(\Omega;\mathbb{R}^{3\times 3}))$ implies that
\begin{equation}
\label{eqn:a1}
\partial_{\mathbf{F}_m}w_R(\phi_m,\mathbf{F}_m)\in L^{\infty}(0,T;L^{2}(\Omega;\mathbb{R}^{3\times 3})).
\end{equation}
Taking $\boldsymbol{\Gamma}=-\Delta \mathbf{F}_m$ in \eqref{eqn:system3galerkin}$_3$ we get, using the Cauchy--Schwarz and Young inequalities, that
\[
\lambda ||\Delta \mathbf{F}_m||_{L^2(\Omega;\mathbb{R}^{3\times 3})}^2\leq \frac{\lambda}{2}||\Delta \mathbf{F}_m||_{L^2(\Omega;\mathbb{R}^{3\times 3})}^2+\frac{1}{\lambda} ||\mathbf{M}_m||_{L^2(\Omega;\mathbb{R}^{3\times 3})}^2+\frac{1}{\lambda} ||\partial_{\mathbf{F}_m}w_R(\phi_m,\mathbf{F}_m)||_{L^2(\Omega;\mathbb{R}^{3\times 3})}^2,
\]
which, after integration in time over the interval $[0,T]$, gives, thanks to \eqref{eqn:apriori}, \eqref{eqn:a1} and elliptic regularity theory, that
\begin{equation}
\label{eqn:a1bis}
||\mathbf{F}_m||_{L^{2}(0,T;H^2(\Omega;\mathbb{R}^{3\times 3}))}\leq C.
\end{equation}
From \eqref{eqn:Fm3d}, \eqref{eqn:a1bis} and \eqref{eqn:gninterpolation1} (applied to $\nabla \mathbf{F}_m$) with $p=4,h=\frac{4}{3}$, we get that
\begin{align}
\label{eqn:a2}
&\notag \mathbf{F}_m \; \text{is uniformly bounded in}\\
& L^{\infty}(0,T;H^{1}(\Omega;\mathbb{R}^{3\times 3}))\cap L^{2}(0,T;H^2(\Omega;\mathbb{R}^{3\times 3}))\hookrightarrow L^{\frac{10}{3}}(0,T;W^{1,\frac{10}{3}}(\Omega;\mathbb{R}^{3\times 3})).
\end{align}
From \eqref{eqn:gninterpolation2} with $p=4,h=2$ we then have that
\begin{align}
\label{eqn:a3}
&\notag \mathbf{F}_m\; \text{is uniformly bounded in}\\
& L^{\infty}(0,T;L^6(\Omega;\mathbb{R}^{3\times 3}))\cap L^{\frac{10}{3}}(0,T;W^{1,\frac{10}{3}}(\Omega;\mathbb{R}^{3\times 3}))\hookrightarrow L^{16}(0,T;L^{8}(\Omega;\mathbb{R}^{3\times 3})).
\end{align}
We also observe that, taking $p=4,h=1$ in \eqref{eqn:gninterpolation1} (applied to $\nabla \mathbf{F}_m$), we have that
\begin{align}
\label{eqn:a2bis}
&||\mathbf{F}_m||_{L^{4}(0,T;W^{1,3}(\Omega;\mathbb{R}^{3\times 3}))}\leq C,
\end{align}
uniformly in $m$, and moreover,  taking $p=4,h=1+\epsilon$, with $\epsilon \in (0,3]$ in \eqref{eqn:gninterpolation1}, we have that
\begin{align}
\label{eqn:a2tris}
&||\mathbf{F}_m||_{L^{4-\frac{8\epsilon}{3(1+\epsilon)}}(0,T;W^{1,3+\epsilon}(\Omega;\mathbb{R}^{3\times 3}))}\leq C,
\end{align}
which implies, employing the Sobolev embedding theorem and defining $k:=\frac{8\epsilon}{3(1+\epsilon)}\in (0,2]$, that
\begin{align}
\label{eqn:a2quatris}
&||\mathbf{F}_m||_{L^{4-k}(0,T;C(\bar{\Omega};\mathbb{R}^{3\times 3}))}\leq C,
\end{align}
uniformly in $m$.

From Assumption \textbf{A2Bis} and \eqref{eqn:a3}, we have that
\begin{equation}
\label{eqn:a4}
\partial_{\mathbf{\phi}_m}w_R(\phi_m,\mathbf{F}_m)\in L^{4}(0,T;L^{2}(\Omega)).
\end{equation}
Taking $\chi=-\Delta \phi_m$ in \eqref{eqn:system3galerkin}$_5$, using the Cauchy--Schwarz and Young inequalities and Assumption \textbf{A3}, we get
\begin{align}
\label{eqn:deltaphi13d}
& \displaystyle \notag ||\Delta \phi_m||^2+\int_{\Omega}\psi_+''(\phi_m)|\nabla \phi_m|^2\leq ||\nabla \mu_m||_{L^2(\Omega;\mathbb{R}^3)}||\nabla \phi_m||_{L^2(\Omega;\mathbb{R}^3)}-\int_{\Omega}\psi_-''(\phi_m)|\nabla \phi_m|^2\\
& \displaystyle \notag +\int_{\Omega}\partial_{\mathbf{\phi}_m}w(\phi_m,\mathbf{F}_m)\Delta \phi_m\leq ||\nabla \mu_m||_{L^2(\Omega;\mathbb{R}^3)}||\nabla \phi_m||_{L^2(\Omega;\mathbb{R}^3)}+\frac{1}{2}||\nabla \phi_m||_{L^2(\Omega;\mathbb{R}^3)}^2\\
& \displaystyle +C\int_{\Omega}\left(1+|\phi_m|^q\right)|\nabla \phi_m|^2+\frac{1}{4}||\Delta \phi_m||^2+||\partial_{\mathbf{\phi}_m}w(\phi_m,\mathbf{F}_m)||^2.
\end{align} 
We use \eqref{eqn:agmon3d}, elliptic regularity theory and the Young inequality, observing that $\frac{4}{q}>1$ when $q<4$, to write
\begin{align*}
&\int_{\Omega}|\phi_m|^q|\nabla \phi_m|^2\leq ||\phi_m||_{L^{\infty}(\Omega)}^q||\nabla \phi_m||_{L^2(\Omega;\mathbb{R}^3)}^2\leq C||\phi_m||_{H^1(\Omega;\mathbb{R}^3)}^{\frac{4+q}{2}}\left(||\phi_m||^{\frac{q}{2}}+||\Delta \phi_m||^{\frac{q}{2}}\right) \\
&\leq C||\phi_m||_{H^1(\Omega;\mathbb{R}^3)}^{2+q}+||\phi_m||_{H^1(\Omega;\mathbb{R}^3)}^{\frac{2(q+4)}{4-q}}+\frac{1}{4}||\Delta \phi_m||^2.
\end{align*}
Using this inequality in \eqref{eqn:deltaphi13d}, together with the Young inequality, we obtain that
\begin{align}
\label{eqn:deltaphi23d}
& \displaystyle \notag ||\Delta \phi_m||^2\leq ||\nabla \mu_m||_{L^2(\Omega;\mathbb{R}^3)}||\nabla \phi_m||_{L^2(\Omega;\mathbb{R}^3)}+C||\phi_m||_{H^1(\Omega;\mathbb{R}^3)}^2+C||\phi_m||_{H^1(\Omega;\mathbb{R}^3)}^{2+q}\\
& \displaystyle +||\phi_m||_{H^1(\Omega;\mathbb{R}^3)}^{\frac{2(q+4)}{4-q}}
+||\partial_{\mathbf{\phi}_m}w(\phi_m,\mathbf{F}_m)||^2+\frac{1}{2}||\Delta \phi_m||^2.
\end{align}
Taking the square of \eqref{eqn:deltaphi23d} and integrating in time over the interval $(0,T)$, using \eqref{eqn:apriori}, Assumption \textbf{A1} and \eqref{eqn:a4}, we infer that
\begin{align}
\label{eqn:a7}
&\phi_m\; \text{is uniformly bounded in}\; L^{\infty}(0,T;H^1(\Omega))\cap L^{4}(0,T;H^2(\Omega)).
\end{align}
Since $L^{\infty}(0,T;H^1(\Omega))\cap L^{4}(0,T;H^2(\Omega))\hookrightarrow L^{\infty}(0,T;L^6(\Omega))\cap L^{4}(0,T;W^{1,6}(\Omega))$, using \eqref{eqn:gninterpolation3} with $s=h=4$ we get that 
\begin{equation}
\label{eqn:phimintermediate}
||\phi_m||_{L^{20}(0,T;L^{10}(\Omega))}\leq C.
\end{equation}
Taking moreover $\chi= 1$ in \eqref{eqn:system3galerkin}$_5$, we obtain, using Assumptions \textbf{A2Bis, A3} and \eqref{eqn:Fm3d}, that
\begin{align}
\label{eqn:mumdiscrete}
& \displaystyle \notag\left|(\mu_m,1)\right|=\left|\int_{\Omega}\psi'\left(\phi_m\right)+\int_{\Omega}\partial_{\phi_m}w_R(\phi_m,\mathbf{F}_m)\right|\\
& \displaystyle \leq C||\phi_m||_{L^l(\Omega)}^l+C||\mathbf{F}_m||_{L^p(\Omega;\mathbb{R}^{3\times 3})}^p+C\leq C||\phi_m||_{L^l(\Omega)}^l+C(\mathbf{F}_0,\phi_0)+C,
\end{align}
where $p\in [0,4]$, $l \in [0,10)$. Taking the square of \eqref{eqn:mumdiscrete} and integrating in time over the interval $[0,T]$, using also \eqref{eqn:phimintermediate}, we obtain that $|(\mu_m,1)|$ is bounded in $L^{2}(0,T)$ and consequently, by the Poincar\'e--Wirtinger inequality and \eqref{eqn:apriori},
\begin{equation}
\label{eqn:mum3d}
\mu_m\; \text{is uniformly bounded in} \; L^{2}(0,T;H^1(\Omega))\hookrightarrow L^{2}(0,T;L^6(\Omega)).
\end{equation}
\subsection{Dual estimates}
We now deduce a priori estimates, uniform in $m$, for the time derivatives of $\mathbf{F}_m$ and $\phi_m$ in \eqref{eqn:system3galerkin}. Multiplying \eqref{eqn:system3galerkin}$_2$ by a time function $\zeta \in L^{\frac{2(4-k)}{2-k}}(0,T)$, with $k \in (0,2)$, choosing $\boldsymbol{\Sigma}=P_m^{L,\Sigma}\boldsymbol{\Pi}$, with a generic $\boldsymbol{\Pi}\in L^2(\Omega;\mathbb{R}^{3\times 3})$, and integrating in time over the interval $(0,T)$, we get
\begin{align}
\label{eqn:fmdual0}
&\notag \int_0^T\int_{\Omega}\partial_t \mathbf{F}_m \colon \boldsymbol{\Pi}\zeta \leq \int_0^T ||\mathbf{v}_m||_{L^6(\Omega;\mathbb{R}^3)}||\nabla \mathbf{F}_m||_{L^3(\Omega;\mathbb{R}^{3\times 3\times 3})}||P_m^{L,\Sigma}\boldsymbol{\Pi}||_{L^{2}(\Omega;\mathbb{R}^{3\times 3})}|\zeta|\\
& \notag +\int_0^T ||\nabla \mathbf{v}_m||_{L^{2}(\Omega;\mathbb{R}^{3\times 3})}||\mathbf{F}_m||_{L^{\infty}(\Omega;\mathbb{R}^{3\times 3})}||P_m^{L,\Sigma}\boldsymbol{\Pi}||_{L^{2}(\Omega;\mathbb{R}^{3\times 3})}|\zeta|\\
& \notag +\gamma \int_0^T||\mathbf{M}_m||_{L^{2}(\Omega;\mathbb{R}^{3\times 3})}||\boldsymbol{\Pi}||_{L^{2}(\Omega;\mathbb{R}^{3\times 3})}|\zeta| \\
& \notag \leq C 
\biggl(||\mathbf{v}_m||_{L^2(0,T;L^6(\Omega;\mathbb{R}^3))}||\nabla \mathbf{F}_m||_{L^{4}(0,T;L^3(\Omega;\mathbb{R}^{3\times 3\times 3}))} \\
& \notag+||\nabla \mathbf{v}_m||_{L^2(0,T;L^{2}(\Omega;\mathbb{R}^{3\times 3}))}||\mathbf{F}_m||_{L^{4-k}(0,T;C(\bar{\Omega};\mathbb{R}^{3\times 3}))}+
||\mathbf{M}_m||_{L^2(0,T;L^{2}(\Omega;\mathbb{R}^{3\times 3}))}\biggr)\\
&  \times ||\boldsymbol{\Pi}||_{L^{2}(\Omega;\mathbb{R}^{3\times 3})}||\zeta||_{L^{\frac{2(4-k)}{2-k}}(0,T)}\leq C||\boldsymbol{\Pi}||_{L^{2}(\Omega;\mathbb{R}^{3\times 3})}||\zeta||_{L^{\frac{2(4-k)}{2-k}}(0,T)},
\end{align}
where, in the last step, we used \eqref{eqn:apriori}, \eqref{eqn:a2bis} and \eqref{eqn:a2quatris}.
Hence, we deduce that
{\rosso
\begin{equation}
\label{eqn:fmdual}
||\partial_t \mathbf{F}_m||_{L^{\frac{4}{3}-s}\left(0,T;L^{2}(\Omega;\mathbb{R}^{3\times 3}) \right)}\leq C,
\end{equation}
where we set $s:=\frac{2k}{3(6-k)}\in \left(0,\frac{1}{3}\right)$. Indeed, it's easy to check that
\[
\frac{2-k}{2(4-k)}+\frac{3}{4-3s}=1
\]
is verified if $s:=\frac{2k}{3(6-k)}$.

Moreover, multiplying \eqref{eqn:system3galerkin}$_4$ by a time function $\zeta \in L^{2}(0,T)$, choosing $\xi=P_m^L(\pi)$, with a generic $\pi \in H^1(\Omega)$, and integrating in time over the interval $(0,T)$, we obtain
\begin{align}
\label{eqn:phimdual0}
& \notag \int_0^T\int_{\Omega}\partial_t\phi_m \pi \zeta \leq \int_0^T ||\mathbf{v}_m||_{L^6(\Omega;\mathbb{R}^3)} ||\nabla \phi_m||_{L^2(\Omega;\mathbb{R}^{3\times})}||P_m^L(\pi)||_{L^{3}(\Omega)}|\zeta|\\
& \notag +\int_0^T||b\left(\phi_m\right)\nabla \mu_m||_{L^{2}(\Omega;\mathbb{R}^2)} ||\nabla P_m^L(\pi)||_{L^{2}(\Omega;\mathbb{R}^2)}|\zeta|\\
& \notag \leq C\left(||\mathbf{v}_m||_{L^2(0,T;L^6(\Omega;\mathbb{R}^3))}||\nabla \phi_m||_{L^{\infty}(0,T;L^2(\Omega))}+||b\left(\phi_m\right)\nabla \mu_m||_{L^2(0,T;L^{2}(\Omega;\mathbb{R}^2))}\right)\\
& \times ||\pi||_{H^1(\Omega)}||\zeta||_{L^2(0,T)}.
\end{align}
Hence, using \eqref{eqn:apriori} and Assumption \textbf{A1} we have that
\begin{equation}
\label{eqn:phimdual}
||\partial_t \phi_m||_{L^{2}\left(0,T;\left(H^1(\Omega)\right)'\right)}\leq C.
\end{equation}
\subsection{Passage to the limit as $m\to \infty$}
Collecting the results \eqref{eqn:vm3d}, \eqref{eqn:Mm3d}, \eqref{eqn:Fm3d}, \eqref{eqn:phim3d}, \eqref{eqn:a2}, \eqref{eqn:a2tris}, \eqref{eqn:a7}, \eqref{eqn:mum3d}, \eqref{eqn:fmdual} and \eqref{eqn:phimdual}, which are uniform in $m$, from the Banach--Alaoglu and the Aubin--Lions lemma, we finally obtain the convergence properties, up to subsequences of the solutions, which we still label by the index $m$, as follows:
\begin{align}
\label{eqn:conv3d1} & \mathbf{v}_m \rightharpoonup \mathbf{v} \quad \text{in} \quad L^{2}\left(0,T;H_{0,\diver}^1\left(\Omega;\mathbb{R}^{3}\right)\right),\\
\label{eqn:conv3d2} & \mathbf{M}_m {\rightharpoonup} \mathbf{M} \quad \text{in} \quad L^{2}(0,T;L^2(\Omega;\mathbb{R}^{3\times 3}))\\
\label{eqn:conv3d3} & \mathbf{F}_m \overset{\ast}{\rightharpoonup} \mathbf{F} \quad \text{in} \quad L^{\infty}(0,T;H^1(\Omega;\mathbb{R}^{3\times 3})),\\
\label{eqn:conv3d3bis} & \mathbf{F}_m {\rightharpoonup} \mathbf{F} \quad \text{in} \quad L^{4-k}(0,T;W^{1,3+\epsilon}(\Omega;\mathbb{R}^{3\times 3}))\cap L^{2}(0,T;H^2(\Omega;\mathbb{R}^{3\times 3})), \; \; \epsilon\in \left(0,3\right],\\
\label{eqn:conv3d4} & \partial_t \mathbf{F}_m \rightharpoonup \partial_t \mathbf{F} \quad \text{in} \quad L^{\frac{4}{3}-s}\left(0,T;L^2(\Omega;\mathbb{R}^{3\times 3})\right), \; \; s\in \left(0,\frac{1}{3}\right),\\
\label{eqn:conv3d5} & \phi_m \overset{\ast}{\rightharpoonup} \phi \quad \text{in} \quad L^{\infty}(0,T;H^1(\Omega))\cap L^{4}(0,T;H^2(\Omega)),\\
\label{eqn:conv3d6} & \mu_m \rightharpoonup \mu \quad \text{in} \quad L^{2}\left(0,T;H^1\left(\Omega\right)\right),\\
\label{eqn:conv3d7} & \partial_t \phi_m \rightharpoonup \partial_t \phi \quad \text{in} \quad L^{2}(0,T;\left(H^1(\Omega)\right)'),\\
\label{eqn:conv3d8} & \mathbf{F}_m \to \mathbf{F} \quad \text{in} \quad C^{0}(0,T;L^p(\Omega;\mathbb{R}^{3\times 3})) \cap L^{2}(0,T;W^{1,p}(\Omega;\mathbb{R}^{3\times 3})), \;\; p\in [1,6), \;\; \text{and} \;\; \text{a.e. in} \; \; \Omega_T,\\
\label{eqn:conv3d8bis} & \mathbf{F}_m \to \mathbf{F} \quad \text{in} \quad L^{4-k}(0,T;C(\bar{\Omega};\mathbb{R}^{3\times 3})),\\
\label{eqn:conv3d9} & \phi_m \to \phi \quad \text{in} \quad C^{0}(0,T;L^p(\Omega)) \cap L^{4}(0,T;W^{1,p}(\Omega)), \;\; p\in [1,6), \;\; \text{and} \;\; \text{a.e. in} \; \; \Omega_T,
\end{align}
with $k\in (0,2]$, as $m\to \infty$. We note that \eqref{eqn:conv3d8bis} follows from the compact embedding $W^{1,3+\epsilon}\subset C^0$ and \eqref{eqn:conv3d3bis}, \eqref{eqn:conv3d4}.
With the convergence results \eqref{eqn:conv3d1}--\eqref{eqn:conv3d9}, we can pass to the limit in the system \eqref{eqn:system3galerkin} as $m\to \infty$.
Let's take $\boldsymbol{\eta}=\boldsymbol{\eta}_m=P_m^{S}(\mathbf{u})$, for arbitrary $\mathbf{u}\in H_{0,\diver}^1(\Omega;\mathbb{R}^3)$, 
$\boldsymbol{\Sigma}=\boldsymbol{\Sigma}_m=P_m^{L,\Sigma}(\boldsymbol{\Theta})$, for arbitrary $\boldsymbol{\Theta}\in L^2(\Omega;\mathbb{R}^{3\times 3})$, $\boldsymbol{\Gamma}=\boldsymbol{\Gamma}_m=P_m^{L,\Sigma}(\boldsymbol{\Pi})$, for arbitrary $\boldsymbol{\Pi}\in H^1(\Omega;\mathbb{R}^{3\times 3})$,
$\xi=\xi_m=P_m^L(q)$, for arbitrary $q\in H^1(\Omega)$, $\chi=\chi_m=P_m^L(r)$, for arbitrary $r\in H^1(\Omega)$, multiply the equations by a function $\omega \in C_0^{\infty}([0,T])$ and integrate over the time interval $[0,T]$. This gives
\begin{equation}
\label{eqn:limit3d1}
\begin{cases}
\displaystyle \nu \int_0^T\omega \int_{\Omega}\nabla \mathbf{v}_m\colon \nabla \boldsymbol{\eta}_m= \int_0^T\omega\int_{\Omega}\mu_m \nabla \phi_m \cdot \boldsymbol{\eta}_m-\int_0^T\omega\int_{\Omega}\left(\mathbf{M}_m\mathbf{F}_m^T\right)\colon \nabla\boldsymbol{\eta}_m\\ 
\displaystyle +\int_0^T\omega\int_{\Omega}\left(\nabla \mathbf{F}_m\odot \mathbf{M}_m\right)\cdot \boldsymbol{\eta}_m,\\ \\
\displaystyle  \int_0^T\omega \int_{\Omega}\partial_t \mathbf{F}_m\colon \boldsymbol{\Sigma}_m +\int_0^T\omega\int_{\Omega}\left({\mathbf{v}}_m \cdot \nabla \right)\mathbf{F}_m\colon \boldsymbol{\Sigma}_m-\int_0^T\omega\int_{\Omega}\left(\nabla {\mathbf{v}}_m\right)\mathbf{F}_m\colon \boldsymbol{\Sigma}_m\\
\displaystyle +\gamma \int_0^T\omega\int_{\Omega}\mathbf{M}_m\colon \boldsymbol{\Sigma}_m=0,\\ \\
\displaystyle \int_0^T\omega \int_{\Omega}\mathbf{M}_m\colon \boldsymbol{\Gamma}_m=\int_0^T\omega \int_{\Omega}\partial_{\mathbf{F}_m}w_R(\phi_m,\mathbf{F}_m)\colon \boldsymbol{\Gamma}_m+\lambda \int_0^T\omega \int_{\Omega}\nabla \mathbf{F}_m \mathbin{\tensorm} \nabla \boldsymbol{\Gamma}_m,\\ \\
 \displaystyle \int_0^T\omega\int_{\Omega}\partial_t\phi_m \xi_m+\int_0^T\omega\int_{\Omega}\left(\mathbf{v}_m \cdot \nabla \phi_m\right) \xi_m+\int_0^T\omega\int_{\Omega}b\left(\phi_m\right)\nabla \mu_m \cdot \nabla \xi_m=0,\\ \\
\displaystyle \int_0^T\omega\int_{\Omega}\mu_m \chi_m=\int_0^T\omega\int_{\Omega}\nabla \phi_m\cdot \nabla \chi_m+\int_0^T\omega\int_{\Omega}\psi'\left(\phi_m\right)\chi_m\\
\displaystyle +\int_0^T\omega\int_{\Omega}\partial_{\phi_m}w_R(\phi_m,\mathbf{F}_m)\chi_m.
\end{cases}
\end{equation} 
We observe that 
\begin{align}
\label{eqn:projstrong}
\begin{cases}
\boldsymbol{\eta}_m=P_m^{S}(\mathbf{u})\rightarrow \mathbf{u}\; \; \text{in} \; \; H^1(\Omega;\mathbb{R}^3),\\
\boldsymbol{\Sigma}_m=P_m^{L,\Sigma}(\boldsymbol{\Theta})\rightarrow \boldsymbol{\Theta}\; \; \text{in} \; \; L^2(\Omega;\mathbb{R}^{3\times 3}),\\
\boldsymbol{\Gamma}_m=P_m^{L,\Sigma}(\boldsymbol{\Pi})\rightarrow \boldsymbol{\Pi}\; \; \text{in} \; \; H^1(\Omega;\mathbb{R}^{3\times 3}),\\
\xi_m=P_m^L(q)\rightarrow q, \chi_m=P_m^L(r)\rightarrow r \; \; \text{in} \; \; H^1(\Omega).
\end{cases}
\end{align}
Thanks to \eqref{eqn:conv3d1} and \eqref{eqn:projstrong}$_1$, we have
\begin{equation}
\label{eqn:lim1}
\nu \int_0^T\omega \int_{\Omega}\nabla \mathbf{v}_m\colon \nabla \boldsymbol{\eta}_m\to \nu \int_0^T\omega \int_{\Omega}\nabla \mathbf{v}\colon \nabla \mathbf{u},
\end{equation}
as $m\to \infty$. Owing to \eqref{eqn:conv3d9}, \eqref{eqn:projstrong}$_1$ and the fact that $\boldsymbol{\eta}_m$ converges in $H^1(\Omega;\mathbb{R}^3)\hookrightarrow L^{6}(\Omega;\mathbb{R}^3)$, we have that $\omega \nabla \phi_m \cdot \boldsymbol{\eta}_m\to \omega \nabla \phi \cdot \mathbf{u}$ strongly in $L^2(0,T;L^{\frac{3}{2}}(\Omega))$. Indeed, it turns out that
\begin{align}
\label{eqn:lim1bis}
& \notag \int_0^T|\omega|^2\left(\int_{\Omega} |\nabla \phi_m\cdot \boldsymbol{\eta}_m -\nabla \phi \cdot \mathbf{u}|^{\frac{3}{2}}\right)^{\frac{4}{3}}=
\int_0^T|\omega|^2\left(\int_{\Omega} |\nabla (\phi_m-\phi)\cdot \boldsymbol{\eta}_m +\nabla \phi \cdot (\boldsymbol{\eta}_m-\mathbf{u})|^{\frac{3}{2}}\right)^{\frac{4}{3}}\\
& \notag \leq C \int_0^T|\omega|^2||\nabla (\phi_m-\phi)||_{L^2(\Omega;\mathbf{R}^3)}^2||\boldsymbol{\eta}_m||_{L^6(\Omega;\mathbf{R}^3)}^2 +C\int_0^T|\omega|^2||\nabla \phi||_{L^2(\Omega;\mathbf{R}^3)}^2||\boldsymbol{\eta}_m-\mathbf{u}||_{L^6(\Omega;\mathbf{R}^3)}^2\\
& \notag \leq C||\omega||_{L^{\infty}(0,T)}^2||\nabla \phi_m-\nabla \phi||_{L^2(0,T;L^2(\Omega;\mathbb{R}^{3}))}^2||\boldsymbol{\eta}_m||_{L^6(\Omega;\mathbb{R}^{3})}^2\\
& +C||\omega||_{L^{\infty}(0,T)}^2||\nabla \phi||_{L^2(0,T;L^2(\Omega;\mathbb{R}^{3}))}^2||\boldsymbol{\eta}_m-\mathbf{u}||_{H^1(\Omega;\mathbb{R}^{3})}^2\to 0,
\end{align}
as $m\to \infty$.
Hence, using \eqref{eqn:conv3d6}, by the product of weak--strong convergence we have that 
\begin{equation}
\label{eqn:lim2}
\int_0^T\omega\int_{\Omega}\mu_m \nabla \phi_m \cdot \boldsymbol{\eta}_m\to \int_0^T\omega\int_{\Omega}\mu \nabla \phi \cdot \mathbf{u},
\end{equation}
as $m\to \infty$. For what concerns the second term on the right hand side of \eqref{eqn:limit3d1}$_1$, thanks to \eqref{eqn:conv3d8bis} and \eqref{eqn:projstrong}$_1$ we have that $\omega \nabla \boldsymbol{\eta}_m \mathbf{F}_m\to \omega \nabla \mathbf{u} \mathbf{F}$ strongly in $L^2(0,T;L^2(\Omega;\mathbb{R}^{3\times 3}))$. Indeed, it holds that
\begin{align*}
&\int_0^T|\omega|^2\int_{\Omega}|\nabla\boldsymbol{\eta}_m\mathbf{F}_m-\nabla\mathbf{u}\mathbf{F}|^2=\int_0^T|\omega|^2\int_{\Omega}|\nabla\boldsymbol{\eta}_m(\mathbf{F}_m-\mathbf{F})+\nabla(\boldsymbol{\eta}_m-\mathbf{u})\mathbf{F}|^2\\
&\leq C \int_0^T|\omega|^2||\mathbf{F}_m-\mathbf{F}||_{L^{\infty}(\Omega;\mathbf{R}^{3\times 3})}^2\int_{\Omega}|\nabla\boldsymbol{\eta}_m|^2+C\int_0^T|\omega|^2||\mathbf{F}||_{L^{\infty}(\Omega;\mathbf{R}^{3\times 3})}^2\int_{\Omega}|\nabla(\boldsymbol{\eta}_m-\mathbf{u})|^2\\
& \leq C||\omega||_{L^{\infty}(0,T)}^2||\mathbf{F}_m-\mathbf{F}||_{L^2(0,T;L^{\infty}(\Omega;\mathbb{R}^{3\times 3}))}^2||\nabla\boldsymbol{\eta}_m||_{L^2(\Omega;\mathbb{R}^{3\times 3})}^2\\
&+C||\omega||_{L^{\infty}(0,T)}^2||\mathbf{F}||_{L^2(0,T;L^{\infty}(\Omega;\mathbb{R}^{3\times 3}))}^2||\boldsymbol{\eta}_m-\mathbf{u}||_{H^1(\Omega;\mathbb{R}^{3\times 3})}^2\to 0,
\end{align*}
as $m\to \infty$. Hence, using \eqref{eqn:conv3d2}, by the product of weak--strong convergence we have that 
\begin{equation}
\label{eqn:lim3}
\int_0^T\omega\int_{\Omega}\left(\mathbf{M}_m\mathbf{F}_m^T\right)\colon \nabla\boldsymbol{\eta}_m\to \int_0^T\omega\int_{\Omega}\left(\mathbf{M}\mathbf{F}^T\right)\colon \nabla\mathbf{u},
\end{equation}
as $m\to \infty$. 
For what concerns the third term on the right hand side of \eqref{eqn:limit3d1}$_1$, thanks to \eqref{eqn:conv3d8}, \eqref{eqn:projstrong}$_1$ and the fact that $\boldsymbol{\eta}_m\in H^1(\Omega;\mathbb{R}^3)\hookrightarrow L^{6}(\Omega;\mathbb{R}^3)$, we have that $\omega \nabla \mathbf{F}_m \boldsymbol{\eta}_m \to \omega \nabla \mathbf{F} \mathbf{u}$ strongly in $L^2(0,T;L^2(\Omega;\mathbb{R}^{3\times 3}))$. Indeed,
\begin{align*}
&\int_0^T|\omega|^2\int_{\Omega} |\nabla \mathbf{F}_m\boldsymbol{\eta}_m-\nabla \mathbf{F}\mathbf{u}|^2=\int_0^T|\omega|^2\int_{\Omega} |\nabla (\mathbf{F}_m-\mathbf{F})\boldsymbol{\eta}_m+\nabla \mathbf{F}(\boldsymbol{\eta}_m-\mathbf{u})|^2\\
& \leq C\int_0^T|\omega|^2||\nabla (\mathbf{F}_m-\mathbf{F})||_{L^3(\Omega;\mathbb{R}^{3\times 3\times 3})}^2||\boldsymbol{\eta}_m||_{L^6(\Omega;\mathbb{R}^{3})}^2\\
&+C\int_0^T|\omega|^2||\nabla \mathbf{F}||_{L^3(\Omega;\mathbb{R}^{3\times 3\times 3})}^2||\boldsymbol{\eta}_m-\mathbf{u}||_{L^6(\Omega;\mathbb{R}^{3})}^2\\
& \leq C||\omega||_{L^{\infty}(0,T)}^2||\nabla (\mathbf{F}_m-\mathbf{F})||_{L^2(0,T;L^3(\Omega;\mathbb{R}^{3\times 3\times 3}))}^2||\boldsymbol{\eta}_m||_{L^6(\Omega;\mathbb{R}^{3})}^2\\
&+C||\omega||_{L^{\infty}(0,T)}^2||\nabla \mathbf{F}||_{L^2(0,T;L^3(\Omega;\mathbb{R}^{3\times 3\times 3}))}^2||\boldsymbol{\eta}_m-\mathbf{u}||_{H^1(\Omega;\mathbb{R}^{3})}^2\to 0,
\end{align*}
as $m\to \infty$. Hence, using \eqref{eqn:conv3d2}, by the product of weak--strong convergence we have that 
\begin{equation}
\label{eqn:lim3}
\int_0^T\omega\int_{\Omega}\left(\nabla \mathbf{F}_m\odot \mathbf{M}_m\right)\cdot \boldsymbol{\eta}_m\to \int_0^T\omega\int_{\Omega}\left(\nabla \mathbf{F}\odot \mathbf{M}\right)\cdot \mathbf{u},
\end{equation}
as $m\to \infty$. 

Now we consider the second equation in \eqref{eqn:limit3d1}. Since
\[\omega \boldsymbol{\Sigma}_m,\omega \boldsymbol{\Theta}\in C^0(0,T;L^2(\Omega;\mathbb{R}^{3\times 3}))\hookrightarrow L^{4+\frac{9s}{1-3s}}\left(0,T;L^2(\Omega;\mathbb{R}^{3\times 3})\right),\] 
with $s\in \left(0,\frac{1}{3}\right)$, we get from \eqref{eqn:conv3d4} and \eqref{eqn:projstrong}$_2$ that
\[
\int_0^T\omega\int_{\Omega}\partial_t \mathbf{F}_m\colon \boldsymbol{\Sigma}_m \to \int_0^T\omega\int_{\Omega}\partial_t \mathbf{F}\colon \boldsymbol{\Theta},
\]
as $m\to \infty$. For what concerns the second term in \eqref{eqn:limit3d1}$_2$, we observe, thanks to \eqref{eqn:conv3d8} and \eqref{eqn:projstrong}$_2$, that $\omega \nabla \mathbf{F}_m\odot \boldsymbol{\Sigma}_m \to \omega \nabla \mathbf{F}\odot \boldsymbol{\Theta}$ strongly in $L^2(0,T;L^{\frac{6}{5}}(\Omega;\mathbb{R}^{3\times 3}))$. Indeed, we infer that
\begin{align*}
&\int_0^T|\omega|^2\left(\int_{\Omega} \left|\nabla\mathbf{F}_m\odot \boldsymbol{\Sigma}_m-\nabla\mathbf{F}\odot \boldsymbol{\Theta}\right|^{\frac{6}{5}}\right)^{\frac{5}{3}}\\
&=\int_0^T|\omega|^2\left(\int_{\Omega} \left|\nabla(\mathbf{F}_m-\mathbf{F})\odot \boldsymbol{\Sigma}_m+\nabla\mathbf{F}\odot (\boldsymbol{\Sigma}_m-\boldsymbol{\Theta})\right|^{\frac{6}{5}}\right)^{\frac{5}{3}}\\
&\leq \int_0^T|\omega|^2||\nabla(\mathbf{F}_m-\mathbf{F})||_{L^3(\Omega;\mathbb{R}^{3\times 3\times 3})}^2||\boldsymbol{\Sigma}_m||_{L^2(\Omega;\mathbb{R}^{3\times 3})}^2\\
& +\int_0^T|\omega|^2||\nabla \mathbf{F}||_{L^3(\Omega;\mathbb{R}^{3\times 3\times 3})}^2||\boldsymbol{\Sigma}_m-\boldsymbol{\Theta}||_{L^2(\Omega;\mathbb{R}^{3\times 3})}^2\\
& \leq C||\omega||_{L^{\infty}(0,T)}^2||\nabla (\mathbf{F}_m-\mathbf{F})||_{L^2(0,T;L^3(\Omega;\mathbb{R}^{3\times 3\times 3}))}^2||\boldsymbol{\Sigma}_m||_{L^2(\Omega;\mathbb{R}^{3\times 3})}^2\\
& + C||\omega||_{L^{\infty}(0,T)}^2||\nabla \mathbf{F}||_{L^2(0,T;L^3(\Omega;\mathbb{R}^{3\times 3\times 3}))}^2||\boldsymbol{\Sigma}_m-\boldsymbol{\Theta}||_{L^2(\Omega;\mathbb{R}^{3\times 3})}^2\to 0,
\end{align*}
as $m\to \infty$. Hence, using \eqref{eqn:conv3d1}, by the product of weak--strong convergence we have that 
\begin{equation}
\label{eqn:lim4}
\int_0^T\omega\int_{\Omega}\left({\mathbf{v}}_m \cdot \nabla \right)\mathbf{F}_m\colon \boldsymbol{\Sigma}_m\to \int_0^T\omega\int_{\Omega}\left({\mathbf{v}} \cdot \nabla \right)\mathbf{F}\colon \boldsymbol{\Theta},
\end{equation}
as $m\to \infty$. 
For what concerns the third term in \eqref{eqn:limit3d1}$_2$, thanks to \eqref{eqn:conv3d8bis} and \eqref{eqn:projstrong}$_2$ we have that $\omega \boldsymbol{\Sigma}_m \mathbf{F}_m^T\to \omega \boldsymbol{\Theta} \mathbf{F}^T$ strongly in $L^2(0,T;L^2(\Omega;\mathbb{R}^{3\times 3}))$. Indeed,
\begin{align*}
&\int_0^T|\omega|^2\int_{\Omega}|\boldsymbol{\Sigma}_m\mathbf{F}_m^T-\boldsymbol{\Theta}\mathbf{F}^T|^2=\int_0^T|\omega|^2\int_{\Omega}|\boldsymbol{\Sigma}_m(\mathbf{F}_m-\mathbf{F})^T+(\boldsymbol{\Sigma}_m-\boldsymbol{\Theta})\mathbf{F}^T|^2\\
&\leq C \int_0^T|\omega|^2||\mathbf{F}_m-\mathbf{F}||_{L^{\infty}(\Omega;\mathbf{R}^{3\times 3})}^2\int_{\Omega}|\boldsymbol{\Sigma}_m|^2+C\int_0^T|\omega|^2||\mathbf{F}||_{L^{\infty}(\Omega;\mathbf{R}^{3\times 3})}^2\int_{\Omega}|\boldsymbol{\Sigma}_m-\boldsymbol{\Theta}|^2\\
& \leq C||\omega||_{L^{\infty}(0,T)}^2||\mathbf{F}_m-\mathbf{F}||_{L^2(0,T;L^{\infty}(\Omega;\mathbb{R}^{3\times 3}))}^2||\boldsymbol{\Sigma}_m||_{L^2(\Omega;\mathbb{R}^{3\times 3})}^2\\
&+C||\omega||_{L^{\infty}(0,T)}^2||\mathbf{F}||_{L^2(0,T;L^{\infty}(\Omega;\mathbb{R}^{3\times 3}))}^2||\boldsymbol{\Sigma}_m-\boldsymbol{\Theta}||_{L^2(\Omega;\mathbb{R}^{3\times 3})}^2\to 0,
\end{align*}
as $m\to \infty$. Hence, using \eqref{eqn:conv3d1}, by the product of weak--strong convergence we have that 
\begin{equation}
\label{eqn:lim5}
\int_0^T\omega\int_{\Omega}\left(\nabla {\mathbf{v}}_m\right)\mathbf{F}_m\colon \boldsymbol{\Sigma}_m\to \int_0^T\omega\int_{\Omega}\left(\nabla {\mathbf{v}}\right)\mathbf{F}\colon \boldsymbol{\Theta},
\end{equation}
as $m\to \infty$. 
Finally, thanks to \eqref{eqn:conv3d2} and to \eqref{eqn:projstrong}$_2$, we have the limit
\begin{equation}
\label{eqn:lim6}
 \int_0^T\omega \int_{\Omega}\mathbf{M}_m\colon \boldsymbol{\Sigma}_m\to \int_0^T\omega \int_{\Omega}\mathbf{M}\colon \boldsymbol{\Theta},
\end{equation}
as $m\to \infty$.

We consider the third equation in \eqref{eqn:limit3d1}. The limit of the first term in \eqref{eqn:limit3d1}$_3$ can be studied similarly to the limit of the last term in \eqref{eqn:limit3d1}$_2$. 
Concerning the second term in \eqref{eqn:limit3d1}$_3$, we employ a similar argument as in \cite[Remark $1$]{garckelam}. Using \eqref{eqn:conv3d8}, \eqref{eqn:conv3d9} and the fact that $\partial_{\mathbf{F}_m}w_R\in C(\mathbb{R}\times \mathbb{R}^{3\times 3};\mathbb{R}^{3\times 3})$, we have that $\partial_{\mathbf{F}_m}w_R(\phi_m,\mathbf{F}_m)\to \partial_{\mathbf{F}}w_R(\phi,\mathbf{F})$ a.e. in $\Omega_T$. We observe that, given \eqref{eqn:conv3d8} and \eqref{eqn:gninterpolation3} with $s=2,h=\frac{24}{5}$, we get that
\begin{equation}
\label{eqn:conv3d8tris} \mathbf{F}_m \to \mathbf{F} \quad \text{in} \quad L^{q}(0,T;L^{\frac{6}{5}q}(\Omega;\mathbb{R}^{3\times 3})), \;\; q\in [1,9).
\end{equation}
Hence, we can prove that $|\mathbf{F}_m|^q\omega \boldsymbol{\Gamma}_m\to |\mathbf{F}|^q\omega \boldsymbol{\Pi}$ strongly in $L^1(0,T;L^1(\Omega;\mathbb{R}^{3\times 3}))$ for $q\in [1,9)$. Indeed
\begin{align*}
&\int_0^T|\omega|\int_{\Omega} \left||\mathbf{F}_m|^q\boldsymbol{\Gamma}_m-|\mathbf{F}|^q\boldsymbol{\Pi}\right|=\int_0^T|\omega|\int_{\Omega} \left|\left(|\mathbf{F}_m|^q-|\mathbf{F}|^q\right)\boldsymbol{\Gamma}_m+|\mathbf{F}|^q\left(\boldsymbol{\Gamma}_m-\boldsymbol{\Pi}\right)\right|\\
&=\int_0^T|\omega|\int_{\Omega} \left|\left(|\mathbf{F}+(\mathbf{F}_m-\mathbf{F})|^q-|\mathbf{F}|^q\right)\boldsymbol{\Gamma}_m+|\mathbf{F}|^q\left(\boldsymbol{\Gamma}_m-\boldsymbol{\Pi}\right)\right|\\
&\leq C\int_0^T|\omega|\int_{\Omega} |\mathbf{F}_m-\mathbf{F}|^q|\boldsymbol{\Gamma}_m|+\int_0^T|\omega|\int_{\Omega}|\mathbf{F}|^q|\boldsymbol{\Gamma}_m-\boldsymbol{\Pi}|\\
& \leq C\int_0^T|\omega|\,||\mathbf{F}_m-\mathbf{F}||_{L^{\frac{6}{5}q}(\Omega;\mathbb{R}^{3\times 3})}^q||\boldsymbol{\Gamma}_m||_{L^{6}(\Omega;\mathbb{R}^{3\times 3})}+\int_0^T|\omega|\,||\mathbf{F}||_{L^{\frac{6}{5}q}(\Omega;\mathbb{R}^{3\times 3})}^q||\boldsymbol{\Gamma}_m-\boldsymbol{\Pi}||_{L^{6}(\Omega;\mathbb{R}^{3\times 3})}\\
& \leq C||\omega||_{L^{\infty}(0,T)}||\mathbf{F}_m-\mathbf{F}||_{L^{q}(0,T;L^{\frac{6}{5}q}(\Omega;\mathbb{R}^{3\times 3}))}||\boldsymbol{\Gamma}_m||_{L^6(\Omega;\mathbb{R}^{3\times 3})}\\
& + C||\omega||_{L^{\infty}(0,T)}||\mathbf{F}||_{L^{q}(0,T;L^{\frac{6}{5}q}(\Omega;\mathbb{R}^{3\times 3}))}||\boldsymbol{\Gamma}_m-\boldsymbol{\Pi}||_{H^1(\Omega;\mathbb{R}^{3\times 3})}\to 0,
\end{align*}
as $m\to \infty$, thanks to \eqref{eqn:conv3d8tris} and \eqref{eqn:projstrong}$_3$.
Then, thanks to the growth behavior in \textbf{A2Bis}, applying a generalized form of the Lebesgue convergence theorem and using also  \eqref{eqn:projstrong}$_3$ we have that
\begin{equation}
\label{eqn:lim7}
 \int_0^T\omega \int_{\Omega}\partial_{\mathbf{F}_m}w_R(\phi_m,\mathbf{F}_m)\colon \boldsymbol{\Gamma}_m\to \int_0^T\omega \int_{\Omega}\partial_{\mathbf{F}}w_R(\phi,\mathbf{F})\colon \boldsymbol{\Pi},
\end{equation}
as $m\to \infty$.
With similar arguments, thanks to the growth behavior in \textbf{A2Bis}, \eqref{eqn:conv3d8tris} and  \eqref{eqn:projstrong}$_4$, we can also conclude that
\begin{equation}
\label{eqn:lim8}
 \int_0^T\omega \int_{\Omega}\partial_{\phi_m}w_R(\phi_m,\mathbf{F}_m)\chi_m\to \int_0^T\omega \int_{\Omega}\partial_{\phi}w_R(\phi,\mathbf{F})r.
\end{equation}
We also observe that, given \eqref{eqn:conv3d9} and \eqref{eqn:gninterpolation3} with $s=4,h=\frac{48}{5}$, we get that
\begin{equation}
\label{eqn:conv3d9bis} \phi_m \to \phi \quad \text{in} \quad L^{q}(0,T;L^{\frac{6}{5}q}(\Omega)), \;\; q\in [1,13).
\end{equation}
Hence, given the growth law and regularity in Assumption \textbf{A3}, together with  \eqref{eqn:projstrong}$_4$, applying similarly a generalized form of the Lebesgue convergence theorem we get that
\begin{equation}
\label{eqn:lim9}
 \int_0^T\omega \int_{\Omega}\psi'\left(\phi_m\right)\chi_m\to \int_0^T\omega \int_{\Omega}\psi'\left(\phi \right)r,
\end{equation}
as $m\to \infty$. We also deduce, thanks to \eqref{eqn:conv3d3} and  \eqref{eqn:projstrong}$_3$, that
\begin{equation}
\label{eqn:lim10}
 \int_0^T\omega \int_{\Omega}\nabla \mathbf{F}_m \mathbin{\tensorm} \nabla \boldsymbol{\Gamma}_m\to \int_0^T\omega \int_{\Omega}\nabla \mathbf{F}\mathbin{\tensorm} \nabla \boldsymbol{\Pi},
\end{equation}
as $m\to \infty$. 

Considering the fourth equation in \eqref{eqn:limit3d1}, since
\[\omega \xi_m,\omega q\in C^0(0,T;H^1(\Omega;\mathbb{R}^{3\times 3}))\hookrightarrow L^{2}\left(0,T;H^1(\Omega)\right),\] 
we get from \eqref{eqn:conv3d7} and \eqref{eqn:projstrong}$_4$ that
\[
\int_0^T\omega\int_{\Omega}\partial_t\phi_m \xi_m \to \int_0^T\omega<\partial_t\phi, q>,
\]
as $m\to \infty$. Concerning the second term in \eqref{eqn:limit3d1}$_4$, we obtain, with similar calculations as in \eqref{eqn:lim1bis}, that $\omega \nabla \phi_m \xi_m\to \omega \nabla \phi q$ strongly in $L^2(0,T;L^{\frac{3}{2}}(\Omega))$. Hence, thanks to \eqref{eqn:conv3d1},
\begin{equation}
\label{eqn:lim11}
 \int_0^T\omega \int_{\Omega}\left(\mathbf{v}_m \cdot \nabla \phi_m\right) \xi_m\to \int_0^T\omega \int_{\Omega}\left(\mathbf{v} \cdot \nabla \phi\right) q,
\end{equation}
as $m\to \infty$. 

As for the last term in \eqref{eqn:limit3d1}$_4$, considering \eqref{eqn:conv3d9} and Assumption \textbf{A1}, $b(\phi_m)\to b(\phi)$ a.e. in $\Omega_T$ and is uniformly bounded, hence by applying the Lebesgue convergence theorem and \eqref{eqn:projstrong}$_4$ we obtain that $b(\phi_m)\nabla \xi_m\to b(\phi)\nabla q$ strongly in $L^{2}\left(0,T;L^2(\Omega;\mathbb{R}^3)\right)$. Then, by \eqref{eqn:conv3d9} we have that 
\begin{equation}
\label{eqn:lim12}
 \int_0^T\omega \int_{\Omega}b\left(\phi_m\right)\nabla \mu_m \cdot \nabla \xi_m \to \int_0^T\omega \int_{\Omega}b\left(\phi\right)\nabla \mu \cdot \nabla q,
\end{equation}
as $m\to \infty$.

Lastly, considering \eqref{eqn:conv3d5}, \eqref{eqn:conv3d6} and \eqref{eqn:projstrong}$_4$, we obtain that
\begin{equation}
\label{eqn:lim13}
 \int_0^T\omega \int_{\Omega}\mu_m \chi_m \to \int_0^T\omega \int_{\Omega}\mu r,
\end{equation}
and
\begin{equation}
\label{eqn:lim14}
 \int_0^T\omega \int_{\Omega}\nabla \phi_m\cdot \nabla \chi_m \to \int_0^T\omega \int_{\Omega}\nabla \phi\cdot \nabla r,
\end{equation}
as $m\to \infty$. We also recall \eqref{eqn:lim8} and \eqref{eqn:lim9} to pass to the limit in \eqref{eqn:limit3d1}$_5$.

We need finally to prove that the initial conditions hold. Due to \eqref{eqn:conv3d8} and \eqref{eqn:conv3d9}, in particular to the facts that $\mathbf{F}_m\to \mathbf{F}$ strongly in $C^0([0,T];L^2(\Omega;\mathbb{R}^{3\times 3}))$ and $\phi_m\to \phi$ strongly in $C^0([0,T];L^2(\Omega))$, and the facts that $\mathbf{F}_m(0)=P_m^{L,\Sigma}(\mathbf{F}_0)\to \mathbf{F}_0$ strongly in $L^2(\Omega,\mathbb{R}^{3\times 3})$ and $\phi_m(0)=P_m^L(\phi_0)\to \phi_0$ strongly in $L^2(\Omega)$, we have that $\mathbf{F}(0)=\mathbf{F}_0$ and $\phi(0)=\phi_0$ a.e. in $\Omega$. We have thus proved Theorem \ref{thm:3d}.

\subsection{Passage to the limit as $R\to \infty$}
\label{sec:rtoinfinity}
We recall that the limit point $(\mathbf{v}_R,\mathbf{F}_R,\mathbf{M}_R,\phi_R,\mu_R)$ of system \eqref{eqn:limit3d1}, as $m\to \infty$, satisfies the system
\begin{equation}
\label{eqn:continuous3dR}
\begin{cases}
\displaystyle \nu \int_{\Omega}\nabla \mathbf{v}_R\colon \nabla \mathbf{u} =\int_{\Omega}\mu_R \nabla \phi_R \cdot \mathbf{u}-\int_{\Omega}\left(\mathbf{M}_R\mathbf{F}_R^T\right)\colon \nabla\mathbf{u}+\int_{\Omega}\left(\nabla \mathbf{F}_R\odot \mathbf{M}_R\right)\cdot \mathbf{u},\\
 \displaystyle \int_{\Omega}\partial_t \mathbf{F}_R \colon \boldsymbol{\Theta}+\int_{\Omega}\left({\mathbf{v}_R} \cdot \nabla \right)\mathbf{F}_R\colon \boldsymbol{\Theta}-\int_{\Omega}\left(\nabla {\mathbf{v}_R}\right)\mathbf{F}_R\colon \boldsymbol{\Theta}+\gamma \int_{\Omega}\mathbf{M}_R\colon \boldsymbol{\Theta}=0,\\
 \displaystyle \int_{\Omega}\mathbf{M}_R\colon \boldsymbol{\Pi}=\int_{\Omega}\partial_{\mathbf{F}_R}w_R(\phi_R,\mathbf{F}_R)\colon \boldsymbol{\Pi}+\lambda \int_{\Omega}\nabla \mathbf{F}_R \mathbin{\tensorm} \nabla \boldsymbol{\Pi},\\
 \displaystyle <\partial_t\phi_R, q>+\int_{\Omega}\left(\mathbf{v}_R \cdot \nabla \phi_R \right) q+\int_{\Omega}b\left(\phi_R\right)\nabla \mu_R \cdot \nabla q =0,\\
\displaystyle \int_{\Omega}\mu_R r=\int_{\Omega}\nabla \phi_R \cdot \nabla r+\int_{\Omega}\psi'\left(\phi_R\right)r+\int_{\Omega} \partial_{\phi_R}w_R(\phi_R,\mathbf{F}_R)r,
\end{cases}
\end{equation} 
a.e. in $(0,T)$ and for all $\mathbf{u}\in H_{0,\diver}^1\left(\Omega;\mathbb{R}^{3}\right)$, $\boldsymbol{\Theta}\in L^2\left(\Omega;\mathbb{R}^{3\times 3}\right)$, $\boldsymbol{\Pi}\in H^1\left(\Omega;\mathbb{R}^{3\times 3}\right)$, $q, r \in H^1(\Omega)$, as well as the initial conditions $\mathbf{F}_R(\cdot,0)=\mathbf{F}_0$ a.e. in $\Omega$ and $\phi_R(\cdot,0)=\phi_0$ a.e. in $\Omega$. We note that in \eqref{eqn:continuous3dR} we have restored the index $R$ to indicate the presence of the truncation. The aim of this section is to study the limit problem of system \eqref{eqn:continuous3dR} as $R\to \infty$.

We start by searching for growth laws for $w_R(\phi_R,\mathbf{F}_R)$ which are uniform in $R$. From \eqref{eqn:grbound} and \eqref{eqn:grpbound} we deduce that
\begin{equation}
\label{eqn:grboundu}
0\leq g_R(r)\leq C \quad \forall r\in \mathbb{R},
\end{equation}
uniformly in $R$, and also that
\begin{equation}
\label{eqn:grpboundu}
|g^{\prime}_R(r)|\leq
\begin{cases}
0 \quad \text{for} \; r\leq R;\\
\frac{C}{r}+\frac{C}{r^{\max(1,p-3)}}\leq \frac{C}{r} \quad \text{for} \; R\leq r \leq 2R;\\
0  \quad \text{for}\, r \geq 2R.
\end{cases}
\end{equation}
Then, given the definition \eqref{eqn:wr} and formula \eqref{eqn:dfwr}, we obtain that
\begin{equation}
\label{eqn:wrboundu}
-d_1\leq w_R(\phi,\mathbf{F})\leq C(1+|\mathbf{F}|^{p}),
\end{equation}
and
\begin{equation}
\label{eqn:wrpboundu}
|\partial_{\mathbf{F}}w_R(\phi,\mathbf{F})|\leq C(1+|\mathbf{F}|^{p-1}), \;\;|\partial_{\phi}w_R(\phi,\mathbf{F})|\leq C(1+|\mathbf{F}|^{p}),
\end{equation}
uniformly in $R$. Hence, Assumption \textbf{A2} is valid for the truncated elastic energy density $w_R$ uniformly in $R$.

The a priori estimate \eqref{eqn:apriori} is uniform in the truncation parameter $R$, as the weak lower semicontinuity of the involved norms translates to the limit. Hence, also arguing as in \eqref{eqn:fmean}-\eqref{eqn:phim3d}, we infer that
\begin{equation}
\label{eqn:vm3du}
\mathbf{v}_R\; \text{is uniformly bounded in} \; L^2(0,T;H_{0,\diver}^1(\Omega;\mathbb{R}^3))\hookrightarrow L^2(0,T;L_{\diver}^6(\Omega;\mathbb{R}^3)),
\end{equation}
\begin{equation}
\label{eqn:Mm3du}
\mathbf{M}_R \; \text{is uniformly bounded in} \; L^{2}(0,T;L^2(\Omega;\mathbb{R}^{3\times 3})),
\end{equation}
\begin{equation}
\label{eqn:Fm3du}
\mathbf{F}_R \; \text{is uniformly bounded in} \; L^{\infty}(0,T;H^1(\Omega;\mathbb{R}^{3\times 3}))\hookrightarrow L^{\infty}(0,T;L^6(\Omega;\mathbb{R}^{3\times 3})),
\end{equation}
\begin{equation}
\label{eqn:phim3du}
\phi_R\; \text{is uniformly bounded in} \; L^{\infty}(0,T;H^1(\Omega))\hookrightarrow L^{\infty}(0,T;L^6(\Omega)),
\end{equation}
\begin{equation}
\label{eqn:mum3dnablau}
\nabla \mu_R\; \text{is uniformly bounded in} \; L^{2}(0,T;L^2(\Omega;\mathbb{R}^3)).
\end{equation}

We now derive higher order estimates for $\mathbf{F}_R$ and $\phi_R$ by employing an iterative bootstrap argument based on elliptic regularity theory. We observe that, if $p\leq 4$, the truncation operation is not active, i.e. $g_R(|\mathbf{F}_R|)\equiv 1$. We thus specialize to the case $4<p<6$. 
We observe, from \eqref{eqn:Fm3du} and Assumption \textbf{A2}, that 
\begin{equation}
\label{eqn:a1u}
\partial_{\mathbf{F}_R}w_R(\phi_R,\mathbf{F}_R)\; \text{is uniformly bounded in} \;  L^{\infty}(0,T;L^{\frac{6}{p-1}}(\Omega;\mathbb{R}^{3\times 3})).
\end{equation}
Then, from equation \eqref{eqn:continuous3dR}$_3$, \eqref{eqn:Mm3du}, elliptic regularity theory and \eqref{eqn:gninterpolation1} with $h=\frac{2(6-p)}{3}$ we deduce that
\begin{align}
\label{eqn:a2u}
& \notag \mathbf{F}_R\; \text{is uniformly bounded in} \;  \\
& L^{\infty}(0,T;H^{1}(\Omega;\mathbb{R}^{3\times 3}))\cap L^{2}(0,T;W^{2,\frac{6}{p-1}}(\Omega;\mathbb{R}^{3\times 3}))\hookrightarrow L^{\frac{18-2p}{3}}(0,T;W^{1,\frac{18-2p}{3}}(\Omega;\mathbb{R}^{3\times 3})).
\end{align}
From \eqref{eqn:gninterpolation2} with $h=2(6-p)$ and \eqref{eqn:Fm3du} we also have that
\begin{align}
\label{eqn:a3u}
& \notag \mathbf{F}_R\; \text{is uniformly bounded in} \; \\
& L^{\infty}(0,T;L^6(\Omega;\mathbb{R}^{3\times 3}))\cap L^{\frac{18-2p}{3}}(0,T;W^{1,\frac{18-2p}{3}}(\Omega;\mathbb{R}^{3\times 3}))\hookrightarrow L^{18-2p}(0,T;L^{18-2p}(\Omega;\mathbb{R}^{3\times 3})).
\end{align}
We split the argument for different sub-intervals of the interval $(4,6)$.
\subsubsection{$4< p \leq 5$: one bootstrap step} 
From \eqref{eqn:a3u} and Assumption \textbf{A2} we get that
\begin{align}
\label{eqn:a4u}
& \notag \partial_{\mathbf{F}_R}w_R(\phi_R,\mathbf{F}_R)\; \text{is uniformly bounded in} \; \\
& L^{\frac{18-2p}{p-1}}(0,T;L^{\frac{18-2p}{p-1}}(\Omega;\mathbb{R}^{3\times 3}))\hookrightarrow L^{2}(0,T;L^{2}(\Omega;\mathbb{R}^{3\times 3})).
\end{align}
Hence, from equation \eqref{eqn:continuous3dR}$_3$, \eqref{eqn:Mm3du}, and applying again elliptic regularity theory we obtain that
\begin{align}
\label{eqn:a5u}
& \notag \mathbf{F}_R\; \text{is uniformly bounded in} \;  \\
& L^{\infty}(0,T;H^{1}(\Omega;\mathbb{R}^{3\times 3}))\cap L^{2}(0,T;H^2(\Omega;\mathbb{R}^{3\times 3}))\hookrightarrow L^{10}(0,T;L^{10}(\Omega;\mathbb{R}^{3\times 3})), 
\end{align}
where we used \eqref{eqn:gninterpolation3} with $s=2,h=4$, which implies that
\begin{align}
\label{eqn:a6u}
& \partial_{\mathbf{\phi}_R}w_R(\phi_R,\mathbf{F}_R)\; \text{is uniformly bounded in} \; L^{\frac{10}{p}}(0,T;L^{\frac{10}{p}}(\Omega))\hookrightarrow L^{2}(0,T;L^{2}(\Omega)),
\end{align}
and, proceeding as in \eqref{eqn:deltaphi23d},
\begin{align}
\label{eqn:a7u}
\phi_R\; \text{is uniformly bounded in} \;  L^{\infty}(0,T;H^1(\Omega))\cap L^{2}(0,T;H^2(\Omega)).
\end{align}

\subsubsection{$5< p \leq \frac{16}{3}$: two bootstrap steps} 
We observe that in this interval \eqref{eqn:a3u} does not imply that $\partial_{\mathbf{F}_R}w_R(\phi_R,\mathbf{F}_R)$ is uniformly bounded in $L^{2}(0,T;L^{2}(\Omega;\mathbb{R}^{3\times 3}))$. Hence, we modify \eqref{eqn:a3u} using \eqref{eqn:gninterpolation2} with $h=\frac{6(6-p)}{2p-7}$, obtaining that
\begin{align}
\label{eqn:a8u}
& \notag \mathbf{F}_R\; \text{is uniformly bounded in} \;  \\
& L^{\infty}(0,T;L^6(\Omega;\mathbb{R}^{3\times 3}))\cap L^{\frac{18-2p}{3}}(0,T;W^{1,\frac{18-2p}{3}}(\Omega;\mathbb{R}^{3\times 3}))\hookrightarrow L^{2(p-1)}(0,T;L^{\frac{6(p-1)}{2p-7}}(\Omega;\mathbb{R}^{3\times 3})).
\end{align} 
Hence, it holds that
\begin{equation}
\label{eqn:a9u}
\partial_{\mathbf{F}_R}w_R(\phi_R,\mathbf{F}_R)\; \text{is uniformly bounded in} \;  L^{2}(0,T;L^{\frac{6}{2p-7}}(\Omega;\mathbb{R}^{3\times 3})),
\end{equation}
and we perform a further bootstrap argument employing elliptic regularity  theory in equation \eqref{eqn:continuous3dR}$_3$.
Introducing the new exponent $p':=2p-6$, we have that
\begin{equation}
\label{eqn:a10u}
\partial_{\mathbf{F}_R}w_R(\phi_R,\mathbf{F}_R)\; \text{is uniformly bounded in} \;  L^{2}(0,T;L^{\frac{6}{p'-1}}(\Omega;\mathbb{R}^{3\times 3})),
\end{equation}
and, proceeding as in \eqref{eqn:a2u} and \eqref{eqn:a3u},
\begin{align}
\label{eqn:a10ubis}
& \notag \mathbf{F}_R\; \text{is uniformly bounded in} \;  \\
& L^{18-2p'}(0,T;L^{18-2p'}(\Omega;\mathbb{R}^{3\times 3}))\equiv L^{30-4p}(0,T;L^{30-4p}(\Omega;\mathbb{R}^{3\times 3})).
\end{align}
Then finally, it turns out that
\begin{align}
\label{eqn:a11u}
& \notag \partial_{\mathbf{F}_R}w_R(\phi_R,\mathbf{F}_R) \; \text{is uniformly bounded in} \; \\
& L^{\frac{30-4p}{p-1}}(0,T;L^{\frac{30-4p}{p-1}}(\Omega;\mathbb{R}^{3\times 3}))\hookrightarrow L^{2}(0,T;L^{2}(\Omega;\mathbb{R}^{3\times 3})).
\end{align}
Hence, from equation \eqref{eqn:continuous3dR}$_3$, \eqref{eqn:Mm3du}, and applying again elliptic regularity  theory we obtain that
\begin{align}
\label{eqn:a12u}
& \notag \mathbf{F}_R \; \text{is uniformly bounded in} \; \\
& L^{\infty}(0,T;H^{1}(\Omega;\mathbb{R}^{3\times 3}))\cap L^{2}(0,T;H^2(\Omega;\mathbb{R}^{3\times 3}))\hookrightarrow L^{\frac{8p}{2p-6}}(0,T;L^{2p}(\Omega;\mathbb{R}^{3\times 3})), 
\end{align}
where we have applied \eqref{eqn:gninterpolation3}, with $s=2,h=2p-6$. Using Assumption \textbf{A2}, \eqref{eqn:a12u} implies that
\begin{align}
\label{eqn:a13u}
\partial_{\mathbf{\phi}_R}w_R(\phi_R,\mathbf{F}_R)\; \text{is uniformly bounded in} \;  L^{\frac{8}{2p-6}}(0,T;L^{2}(\Omega)).
\end{align}
In conclusion, taking \eqref{eqn:deltaphi23d} to the power of $\frac{4}{2p-6}$ and integrating in time over the interval $(0, T)$, we obtain that
\begin{align}
\label{eqn:a14u}
& \notag \phi_R\; \text{is uniformly bounded in} \;  \\
& L^{\infty}(0,T;H^1(\Omega))\cap L^{\frac{8}{2p-6}}(0,T;H^2(\Omega))\hookrightarrow L^{\frac{12p-20}{3(2p-6)}}(0,T;W^{1,\frac{12p-20}{3(2p-6)}}(\Omega)),
\end{align}
where we have used \eqref{eqn:lpinterpolation}, with $s=2,q=6$.
We observe that $\frac{12p-20}{3(2p-6)}>2$ for any $p>3$.

\subsubsection{$\frac{16}{3}< p \leq \frac{11}{2}$: three bootstrap steps}
We observe that in this interval \eqref{eqn:a10ubis} does not imply that \[\partial_{\mathbf{F}_R}w_R(\phi_R,\mathbf{F}_R)\; \text{is uniformly bounded in} \;  L^{2}(0,T;L^{2}(\Omega;\mathbb{R}^{3\times 3})).\]
Hence, we modify \eqref{eqn:a10ubis} using \eqref{eqn:gninterpolation2} with $h=\frac{6(6-p')}{p+p'-7}$, obtaining that
\begin{align}
\label{eqn:a15u}
& \notag \mathbf{F}_R \; \text{is uniformly bounded in} \; \\
&L^{\infty}(0,T;L^6(\Omega;\mathbb{R}^{3\times 3}))\cap L^{\frac{18-2p'}{3}}(0,T;W^{1,\frac{18-2p'}{3}}(\Omega;\mathbb{R}^{3\times 3}))\hookrightarrow L^{2(p-1)}(0,T;L^{\frac{6(p-1)}{p+p'-7}}(\Omega;\mathbb{R}^{3\times 3})).
\end{align} 
Hence, we infer that
\begin{equation}
\label{eqn:a16u}
\partial_{\mathbf{F}_R}w_R(\phi_R,\mathbf{F}_R)\; \text{is uniformly bounded in} \; L^{2}(0,T;L^{\frac{6}{p+p'-7}}(\Omega;\mathbb{R}^{3\times 3})).
\end{equation}
Introducing the new exponent $p'':=p+p'-6$, we have that
\begin{equation}
\label{eqn:a17u}
\partial_{\mathbf{F}_R}w_R(\phi_R,\mathbf{F}_R)\; \text{is uniformly bounded in} \;  L^{2}(0,T;L^{\frac{6}{p''-1}}(\Omega;\mathbb{R}^{3\times 3})),
\end{equation}
and, proceeding as in \eqref{eqn:a2u} and \eqref{eqn:a3u},
\begin{align*}
& \mathbf{F}_R\; \text{is uniformly bounded in} \; \\
&L^{18-2p''}(0,T;L^{18-2p''}(\Omega;\mathbb{R}^{3\times 3}))\equiv L^{42-6p}(0,T;L^{42-6p}(\Omega;\mathbb{R}^{3\times 3})).
\end{align*}
Then, we have that
\begin{align}
\label{eqn:a18u}
& \notag \partial_{\mathbf{F}_R}w_R(\phi_R,\mathbf{F}_R)\; \text{is uniformly bounded in} \;  \\
&L^{\frac{42-6p}{p-1}}(0,T;L^{\frac{42-6p}{p-1}}(\Omega;\mathbb{R}^{3\times 3}))\hookrightarrow L^{2}(0,T;L^{2}(\Omega;\mathbb{R}^{3\times 3})).
\end{align}
Hence, from equation \eqref{eqn:continuous3dR}$_3$, \eqref{eqn:Mm3du}, and applying again elliptic regularity theory we obtain, analogously to \eqref{eqn:a12u}, that
\begin{align}
\label{eqn:a19u}
& \notag \mathbf{F}_R\; \text{is uniformly bounded in} \;  \\
& L^{\infty}(0,T;H^{1}(\Omega;\mathbb{R}^{3\times 3}))\cap L^{2}(0,T;H^2(\Omega;\mathbb{R}^{3\times 3}))\hookrightarrow L^{\frac{8p}{2p-6}}(0,T;L^{2p}(\Omega;\mathbb{R}^{3\times 3})), 
\end{align}
which implies, analogously to \eqref{eqn:a14u}, that
\begin{align}
\label{eqn:a20u}
& \notag \phi_R\; \text{is uniformly bounded in} \;  \\
&L^{\infty}(0,T;H^1(\Omega))\cap L^{\frac{8}{2p-6}}(0,T;H^2(\Omega))\hookrightarrow L^{\frac{12p-20}{3(2p-6)}}(0,T;L^{\frac{12p-20}{3(2p-6)}}(\Omega)).
\end{align}

\subsubsection{Iterative argument up to $p<6$}
Based on the previous observations, we search for an iteration argument which let us apply $n$ bootstrap steps to obtain higher order regularity of $\mathbf{F}_R$ and $\phi_R$ in proper contracting intervals for $p$, up to $p<6$.
Let $n\in \mathbb{N}_+$, and 
\begin{equation}
\label{eqn:a21u}
\frac{10+6(n-1)}{(n-1)+2}<p\leq \frac{10+6n}{n+2},
\end{equation}
where $n$ is the number of bootstrap steps which must be applied in the interval defined in \eqref{eqn:a21u} to obtain higher order regularity of $\mathbf{F}_R$ and $\phi_R$.
We iteratively define the sequence $\{p^n\}_{n\in \mathbb{N}}$, where $n$ is an index (not an exponent), with
\begin{align}
\label{eqn:iteration}
&p^0=p;\\
& \notag p^n=p+p^{n-1}-6=n(p-6)+p,
\end{align}
 and $p$ as in \eqref{eqn:a21u}. Iterating $n$ bootstrap steps, we have that
\[
\mathbf{F}_R \; \text{is uniformly bounded in} \;  L^{18-2p^n}(0,T;L^{18-2p^n}(\Omega;\mathbb{R}^{3\times 3})),
\]
and hence
\begin{align}
\label{eqn:a22u}
& \notag \partial_{\mathbf{F}_R}w_R(\phi_R,\mathbf{F}_R)\; \text{is uniformly bounded in} \;  \\
& L^{\frac{18-2p^n}{p-1}}(0,T;L^{\frac{18-2p^n}{p-1}}(\Omega;\mathbb{R}^{3\times 3}))\hookrightarrow L^{2}(0,T;L^{2}(\Omega;\mathbb{R}^{3\times 3})).
\end{align}
Note that 
\[
\frac{18-2p^n}{p-1}=\frac{18-2[n(p-6)+p]}{p-1}=2 \;\; \text{if} \;\; p=\frac{10+6n}{n+2}.
\]
Hence, from \eqref{eqn:gninterpolation2}, analogously to \eqref{eqn:a12u}, it follows that
\begin{align}
\label{eqn:a23u}
& \notag \mathbf{F}_R\; \text{is uniformly bounded in} \;  \\
&L^{\infty}(0,T;H^{1}(\Omega;\mathbb{R}^{3\times 3}))\cap L^{2}(0,T;H^2(\Omega;\mathbb{R}^{3\times 3}))\hookrightarrow L^{\frac{8p}{2p-6}}(0,T;L^{2p}(\Omega;\mathbb{R}^{3\times 3})), 
\end{align}
which implies, analogously to \eqref{eqn:a14u}, that
\begin{align}
\label{eqn:a24u}
& \notag \phi_R\; \text{is uniformly bounded in} \;  \\
& L^{\infty}(0,T;H^1(\Omega))\cap L^{\frac{8}{2p-6}}(0,T;H^2(\Omega))\;\hookrightarrow L^{\frac{12p-20}{3(2p-6)}}(0,T;L^{\frac{12p-20}{3(2p-6)}}(\Omega)).
\end{align}
We note that, thanks to \eqref{eqn:a23u}, the bound \eqref{eqn:a2tris} is still valid uniformly in $R$.\\
\noindent
Taking $n\to \infty$, we obtain this result for $p\to 6$. We observe that for $p=6$, $p^n=p=6$, $\partial_{\mathbf{F}_R}w_R(\phi_R,\mathbf{F}_R)\in L^{\frac{6}{5}}(0,T;L^{\frac{6}{5}}(\Omega;\mathbb{R}^{3\times 3}))$ and the argument breaks down.

We highlight that \eqref{eqn:a23u} and \eqref{eqn:a24u} are valid for any $p\in (4,6)$.
Since $L^{\infty}(0,T;H^1(\Omega))\cap L^{\frac{8}{2p-6}}(0,T;H^2(\Omega))\hookrightarrow L^{\infty}(0,T;L^6(\Omega))\cap L^{\frac{8}{2p-6}}(0,T;W^{1,6}(\Omega))$, using \eqref{eqn:gninterpolation3} with $s=h=\frac{8}{2p-6}$ we get that 
\begin{equation}
\label{eqn:phimintermediate2}
||\phi_R||_{L^{2q}(0,T;L^{q}(\Omega))}\leq C,\quad q=6+\frac{8}{2p-6}.
\end{equation}
Then, with similar calculations as in \eqref{eqn:mumdiscrete} and using Assumptions \textbf{A2, A3}, we obtain that  $|(\mu_R,1)|$ is bounded in $L^{2}(0,T)$ and consequently, by the Poincar\'e--Wirtinger inequality and \eqref{eqn:apriori},
\begin{equation}
\label{eqn:mum3du}
\mu_R\; \text{is uniformly bounded in} \; L^{2}(0,T;H^1(\Omega)) \; \hookrightarrow L^{2}(0,T;L^6(\Omega)).
\end{equation}
We finally observe that, thanks to \eqref{eqn:vm3du}--\eqref{eqn:mum3dnablau}, \eqref{eqn:a23u}, \eqref{eqn:a24u} and \eqref{eqn:mum3du} the estimates \eqref{eqn:fmdual} and \eqref{eqn:phimdual} translate to the limit.

Collecting the obtained estimates, which are uniform in $R$, 
from the Banach--Alaoglu and the Aubin--Lions lemma, we finally obtain the convergence properties, up to subsequences of the solutions, which we still label by the index $R$,
\begin{align}
\label{eqn:conv3d1u} & \mathbf{v}_R \rightharpoonup \mathbf{v} \quad \text{in} \quad L^{2}\left(0,T;H_{0,\diver}^1\left(\Omega;\mathbb{R}^{3}\right)\right),\\
\label{eqn:conv3d2u} & \mathbf{M}_R {\rightharpoonup} \mathbf{M} \quad \text{in} \quad L^{2}(0,T;L^2(\Omega;\mathbb{R}^{3\times 3})),\\
\label{eqn:conv3d3u} & \mathbf{F}_R \overset{\ast}{\rightharpoonup} \mathbf{F} \quad \text{in} \quad L^{\infty}(0,T;H^1(\Omega;\mathbb{R}^{3\times 3})),\\
\label{eqn:conv3d3ubis} & \mathbf{F}_R {\rightharpoonup} \mathbf{F} \quad \text{in} \quad L^{4-k}(0,T;W^{1,3+\epsilon}(\Omega;\mathbb{R}^{3\times 3}))\cap L^{2}(0,T;H^2(\Omega;\mathbb{R}^{3\times 3})), \; \; \epsilon\in \left(0,3\right],\\
\label{eqn:conv3d4u} & \partial_t \mathbf{F}_R \rightharpoonup \partial_t \mathbf{F} \quad \text{in} \quad L^{\frac{4}{3}-s}\left(0,T;L^2(\Omega;\mathbb{R}^{3\times 3})\right), \; \; s\in \left(0,\frac{1}{3}\right),\\
\label{eqn:conv3d5u} & \phi_R \overset{\ast}{\rightharpoonup} \phi \quad \text{in} \quad L^{\infty}(0,T;H^1(\Omega))\cap L^{\frac{8}{2p-6}}(0,T;H^2(\Omega)),\\
\label{eqn:conv3d6u} & \mu_R \rightharpoonup \mu \quad \text{in} \quad L^{2}\left(0,T;H^1\left(\Omega\right)\right),\\
\label{eqn:conv3d7u} & \partial_t \phi_R \rightharpoonup \partial_t \phi \quad \text{in} \quad L^{2}(0,T;\left(H^1(\Omega)\right)'),\\
\label{eqn:conv3d8u} & \mathbf{F}_R \to \mathbf{F} \quad \text{in} \quad C^{0}(0,T;L^p(\Omega;\mathbb{R}^{3\times 3})) \cap L^{2}(0,T;W^{1,p}(\Omega;\mathbb{R}^{3\times 3})), \;\; p\in [1,6), \;\; \text{and} \;\; \text{a.e. in} \; \; \Omega_T,\\
\label{eqn:conv3d8ubis} & \mathbf{F}_R \to \mathbf{F} \quad \text{in} \quad L^{4-k}(0,T;C(\bar{\Omega};\mathbb{R}^{3\times 3})), \\
\label{eqn:conv3d9u} & \phi_R \to \phi \quad \text{in} \quad C^{0}(0,T;L^p(\Omega)) \cap L^{\frac{8}{2p-6}}(0,T;W^{1,p}(\Omega)), \;\; p\in [1,6), \;\; \text{and} \;\; \text{a.e. in} \; \; \Omega_T,
\end{align}
with $k\in (0,2]$, as $R\to \infty$.
With the convergence results \eqref{eqn:conv3d1u}--\eqref{eqn:conv3d9u}, we can pass to the limit in the system \eqref{eqn:continuous3dR} as $R\to \infty$, with essentially the same calculations as the ones employed to pass to the limit in \eqref{eqn:system3galerkin} as $m\to \infty$, with fixed test functions $\mathbf{u}\in H_{0,\diver}^1\left(\Omega;\mathbb{R}^{3}\right)$, $\boldsymbol{\Theta}\in L^2\left(\Omega;\mathbb{R}^{3\times 3}\right)$, $\boldsymbol{\Pi}\in H^1\left(\Omega;\mathbb{R}^{3\times 3}\right)$, $q, r \in H^1(\Omega)$.
We only point out the fact that, thanks to \eqref{eqn:conv3d9u} and to \eqref{eqn:lpinterpolation}, with $s=2,q=6$, we have that
\begin{equation}
\label{eqn:conv3d9ubis}
\nabla \phi_R \to \nabla \phi \quad \text{in} \quad L^{q}(0,T;L^{q}(\Omega;\mathbb{R}^3)), \quad q\in \left[1,\frac{12p-20}{3(2p-6)}\right),
\end{equation}
and, since, as already observed, $\frac{12p-20}{3(2p-6)}>2$ for any $p>3$, 
\begin{equation}
\label{eqn:conv3d9utris}
\nabla \phi_R \to \nabla \phi \quad \text{in} \quad L^{2}(0,T;L^{2}(\Omega;\mathbb{R}^3)).
\end{equation}
Hence, thanks to \eqref{eqn:conv3d9utris}, we can pass to the limit as $R\to \infty$ e.g. in the first term on the right hand side of \eqref{eqn:continuous3dR}$_1$ with similar calculations to the ones employed in \eqref{eqn:lim2} (with fixed test functions).
\noindent
Also, given \eqref{eqn:conv3d9u} and \eqref{eqn:gninterpolation3} with $s=\frac{8}{2p-6},h=\frac{48}{5(p-3)}$, we get that
\begin{equation}
\label{eqn:conv3d9uquatris} \phi_m \to \phi \quad \text{in} \quad L^{q}(0,T;L^{\frac{6}{5}q}(\Omega)), \;\; q\in \left[1,\frac{5p-7}{p-3}\right).
\end{equation}
Hence, given the growth law and regularity in Assumption \textbf{A3}, and observing that
\[
\frac{5p-7}{p-3}>6+\frac{8}{2p-6}, \quad \forall p<7,
\]
applying a generalized form of the Lebesgue convergence theorem as in \eqref{eqn:lim9} we can pass to the limit as $R\to \infty$ in the second term on the right hand side of \eqref{eqn:continuous3dR}$_5$. We have thus proved Theorem \ref{thm:3d} also for $p\in (4,6)$.

\section{Finite element approximations of the model}
In this section we propose two unconditionally gradient stable fully discrete finite element approximations of system \eqref{eqn:system3}-\eqref{eqn:systembc}. The gradient stability of the approximation schemes guarantees that the discrete system maintains the dissipative nature of the continuous system, with the possibility of defining a discrete energy which is a decreasing Lyapunov functional in time for the discrete solutions.
\begin{rem}
\label{rem:finiteexistence}
We observe that an existence result for System \eqref{eqn:system3r}, i.e. the system with the truncated elastic energy density with polynomial growth up to $p\leq 4$, could be in principle obtained by studying the convergence of one of the finite element approximations introduced in this section, instead of using the Faedo--Galerkin approximation introduced in Section \ref{sec:approx}. Since the existence result obtained in Section \ref{sec:rtoinfinity} with polynomial growth of the elastic energy density up to $p<6$ is obtained using elliptic regularity theory techniques, which are not available in the discrete systems associated to standard finite element approximations, the general existence result expressed in Theorem \ref{thm:3d} cannot be obtained by studying the convergence of the forthcoming finite element approximations.   
\end{rem}

Let $ h>0 $ be a discretization parameter and let $ \mathcal{T}_h$ be a quasi-uniform conforming decomposition of the domain $ \Omega\subset \mathbb{R}^3 $ into $ 3 $-simplices $ K $, with $ h_K=\mathrm{diam}(K) $ and $ h = \max_{K\in\mathcal{T}_h}h_K $. We introduce the following finite-element spaces for scalar, vector and matrix valued functions: 
\begin{align*}
	S_h & := \{s_h\in C^0(\bar{\Omega}): s_h|_K\in\mathbb{P}_1(K), \forall K\in \mathcal{T}_h\} \subset H^1(\Omega),\\
	V_h & := \{\mathbf{v}_h\in C^0(\bar{\Omega};\mathbb{R}^3)\cap H_0^1(\Omega;\mathbb{R}^3): \mathbf{v}_h|_K\in[\mathbb{P}_2(K)]^3, \forall K\in \mathcal{T}_h\},\\
	X_h & := \{\mathbf{A}_h\in C^0(\bar{\Omega};\mathbb{R}^{3\times 3}): \mathbf{A}_h|_K\in[\mathbb{P}_1(K)]^{3\times 3}, \forall K\in \mathcal{T}_h\}\subset H^1(\Omega;\mathbb{R}^{3\times 3}),\\
\end{align*}
where $ \mathbb{P}_l(K) $, $l\in \mathbb{N}$, stands for the space of polynomials of total order $l$ in $ K $. 
The spaces $V_h$ and $R_h:=S_h\cap L_0^2(\Omega)$, with $L_0^2(\Omega):=\{f\in L^2(\Omega)|\int_{\Omega}f=0\}$, constitute the lowest order Taylor-Hood elements \cite{taylorhood}, which are stable elements (i.e. satisfy the discrete Ladyzhenskaya-Babu\v{s}ka-Brezzi stability condition) to approximate the velocity and pressure variables respectively in the Stokes system.
We also introduce the $ L^2 $-projection operators $ P_h^S:L^2(\Omega) \rightarrow S_h $, $ P_h^V:L^2(\Omega;\mathbb{R}^3) \rightarrow V_h $ and $ P_h^X:L^2(\Omega;\mathbb{R}^{3\times 3}) \rightarrow X_h $.
We set $\Delta t = T /N$ for a $N\in \mathbb{N}$, and $ t_n=n\Delta t, n=0,\dots,N $. A finite element approximation of a continuous field $f(\mathbf{x},t_n)$ at time $t_n$ will be indicated by $f_h^n$. Given the initial data $\phi_0\in H^1(\Omega)$, $\mathbf{F}_0\in H^1(\Omega;\mathbb{R}^{3\times 3})$, we set $\phi_h^0=P_h^S(\phi_0)$ and $\mathbf{F}_h^0=P_h^X(\mathbf{F}_0)$.

The first fully discrete approximation scheme is a generalization to system \eqref{eqn:system3} of the convex splitting methods which are widely used to derive unconditionally gradient stable schemes for the Cahn--Hilliard equations \cite{ellstu,eyre}. In order to derive such a scheme, we need to assume a particular form for the elastic energy density $w(\phi,\mathbf{F})$, i.e. we make the following assumption:
\begin{itemize}
\item[$\bf{A2_{h}}$] There exist functions $f\in C^{1}(\mathbb{R})$, with $-k_{f1}\leq f(r) \leq k_{f2}$, $|f^{\prime}(r)|\leq k_{f3}$, $k_{f1},k_{f2},k_{f3}\geq 0$, for all $r\in \mathbb{R}$, functions $g,h\in C^1(\mathbb{R}^{3\times 3})$, with $0\leq g(\mathbf{T})\leq k_g(1+|\mathbf{T}|^p)$, $k_g>0$, $-k_{h1}\leq h(\mathbf{T})\leq k_{h2}(1+|\mathbf{T}|^p)$ and $h(\mathbf{T})\sim k_{h2}|\mathbf{T}|^p$ for $|\mathbf{T}|\to +\infty$, $k_{h1}\geq 0,k_{h2}>0$ and $k_{h2}\geq k_{f1}k_g$, $p\in [0,6)$, for all $\mathbf{T}\in \mathbb{R}^{3\times 3}$, and function $m\in C^1(\mathbb{R})$, with $-k_m\leq m(s)$, $k_m\geq 0$ for all $s\in \mathbb{R}$, such that
\begin{equation}
\label{eqn:wassum}
w(\phi,\mathbf{F})=f(\phi)g(\mathbf{F})+h(\mathbf{F})+m(\phi),
\end{equation}
\end{itemize}
where $m^{\prime}(\cdot)$ satisfies the same assumptions as $\psi^{\prime}(\cdot)$ in \textbf{A3}.
We note that \eqref{eqn:wassum}, with the assumed properties of $f$, $g$, $h$ and $j$, satisfies Assumption \textbf{A2} with $d_1=k_{h1}+k_m$. We moreover observe that the term $m(\phi)$ in \eqref{eqn:wassum} can be incorporated in the term $\psi(\phi)$, so that in the following analysis we will consider as new Cahn--Hilliard potential $\psi(\phi)\leftarrow \psi(\phi)+m(\phi)$.  
\begin{rem}
\label{rem:discdeg}
Under the assumptions $f\in C^{1}(\mathbb{R})$, with $0\leq f(r) \leq k_1$, $|f^{\prime}(r)|\leq k_2$, $k_1,k_2\geq 0$, for all $r\in \mathbb{R}$, and $g\in C^1(\mathbb{R}^{3\times 3})$, with $0\leq g(\mathbf{T})\leq k_g(1+|\mathbf{T}|^p)$, $k_g>0$, for all $\mathbf{T}\in \mathbb{R}^{3\times 3}$, Hypothesis \textbf{A2} could be satisfied also in the case $w(\phi,\mathbf{F})=f(\phi)g(\mathbf{F})$. This means that, when the function $f$ in \eqref{eqn:wassum} is positive, we could consider in the present framework also an elastic energy density which degenerates with the variable $\phi$.
\end{rem}
We introduce the following convex splittings
\begin{align}
& \psi(\phi)=\psi_+(\phi)+\psi_-(\phi),\\
& \notag f(\phi)=f_+(\phi)+f_-(\phi),\\
& \notag g(\mathbf{F})=g_+(\mathbf{F})+g_-(\mathbf{F}),\\
& \notag h(\mathbf{F})=h_+(\mathbf{F})+h_-(\mathbf{F}),
\end{align}
into convex (indicated by the $+$ index) and concave (indicated by the $-$ index) parts. 
\begin{rem}
 \label{rem:convsplit}
 The existence of a convex splitting for a generic function $g\in C^1(\mathbb{R}^{3\times 3})$ is guaranteed if e.g. we make the non restrictive assumption that, for any $\mathbf{F}_1,\mathbf{F}_2\in \mathbb{R}^{3\times 3}$, there exists a $c_0>0$ such that
 \begin{equation}
  \label{ass:conv}
  \left(g^{\prime}(\mathbf{F}_2)-g^{\prime}(\mathbf{F}_1)\right)\colon (\mathbf{F}_2-\mathbf{F}_1)\geq -c_0|\mathbf{F}_2-\mathbf{F}_1|^2.
 \end{equation}
 Indeed, let us choose
 \[
  g_+(\mathbf{F}):=g(\mathbf{F})+\frac{c_0}{2}|\mathbf{F}|^2.
 \]
 Hence,
 \[
  \left(g^{\prime}_+(\mathbf{F}_2)-g^{\prime}_+(\mathbf{F}_1)\right)\colon (\mathbf{F}_2-\mathbf{F}_1)=\left(g^{\prime}(\mathbf{F}_2)-g^{\prime}(\mathbf{F}_1)\right)\colon (\mathbf{F}_2-\mathbf{F}_1)+c_0|\mathbf{F}_2-\mathbf{F}_1|^2\geq 0,
 \]
 which implies that $g^{\prime}_+$ is a monotone function, and therefore $g_+$ is convex. We then set
 \[
  g_-(\mathbf{F}):=-\frac{c_0}{2}|\mathbf{F}|^2,
 \]
whith $g_-$ a concave function.
 \end{rem}
 We also introduce the notation $[f(\cdot)]_+:=\max(0,f(\cdot))$ for the positive part of a given function $f\in C^0(\mathbb{R})$.
We then consider the following fully discretized approximation of system \eqref{eqn:system3}:
\\ \\
\noindent
\textbf{Problem} $\mathbf{P_h^I}$: for $n=0,\dots,N-1$, given $\phi_h^n\in S_h$, $\mathbf{F}_h^n\in X_h$,\\
 find $(\mathbf{v}_h^{n+1},s_h^{n+1},\mathbf{F}_h^{n+1},\mathbf{M}_h^{n+1},\phi_h^{n+1},\mu_h^{n+1})\in V_h\times R_h \times X_h\times X_h\times S_h\times S_h$ such that, \\
for all $(\mathbf{u}_h,p_h,\boldsymbol{\Theta}_h,\boldsymbol{\Pi}_h,\xi_h,\chi_h)\in V_h\times R_h \times X_h\times X_h\times S_h\times S_h$,
\begin{equation}
\label{eqn:galerkinfe1}
\begin{cases}
\displaystyle \nu \int_{\Omega}\nabla \mathbf{v}_h^{n+1}\colon \nabla \mathbf{u}_h-\int_{\Omega}s_h^{n+1}\diver \mathbf{u}_h =\int_{\Omega}\mu_h^{n+1} \nabla \phi_h^n \cdot \mathbf{u}_h-\int_{\Omega}\left(\mathbf{M}_h^{n+1}(\mathbf{F}_h^n)^T\right)\colon \nabla\mathbf{u}_h\\
\displaystyle+\int_{\Omega}\left(\nabla \mathbf{F}_h^n\odot \mathbf{M}_h^{n+1}\right)\cdot \mathbf{u}_h,\\ \\
\displaystyle \int_{\Omega}p_h\diver \mathbf{v}_h^{n+1}=0,\\ \\
 \displaystyle \int_{\Omega} \left(\mathbf{F}_h^{n+1}-\mathbf{F}_h^{n}\right) \colon \boldsymbol{\Theta}_h+\Delta t \int_{\Omega}\left({\mathbf{v}_h^{n+1}} \cdot \nabla \right)\mathbf{F}_h^n\colon \boldsymbol{\Theta}_h-\Delta t\int_{\Omega}\left(\nabla {\mathbf{v}_h^{n+1}}\right)\mathbf{F}_h^n\colon \boldsymbol{\Theta}_h\\
 \displaystyle+\Delta t\gamma \int_{\Omega}\mathbf{M}_h^{n+1}\colon \boldsymbol{\Theta}_h=0,\\ \\
 \displaystyle \int_{\Omega}\mathbf{M}_h^{n+1}\colon \boldsymbol{\Pi}_h=\int_{\Omega}\left(h_+^{\prime}(\mathbf{F}_h^{n+1})+h_-^{\prime}(\mathbf{F}_h^{n})\right)\colon \boldsymbol{\Pi}_h+\int_{\Omega}[f(\phi_h^n)]_+\left(g_+^{\prime}(\mathbf{F}_h^{n+1})+g_-^{\prime}(\mathbf{F}_h^{n})\right)\colon \boldsymbol{\Pi}_h\\
 \displaystyle +\int_{\Omega}(f(\phi_h^n)-[f(\phi_h^n)]_+)\left(g_+^{\prime}(\mathbf{F}_h^{n})+g_-^{\prime}(\mathbf{F}_h^{n+1})\right)\colon \boldsymbol{\Pi}_h +\lambda \int_{\Omega}\nabla \mathbf{F}_h^{n+1} \mathbin{\tensorm} \nabla \boldsymbol{\Pi}_h,\\ \\
 \displaystyle \int_{\Omega}\left(\phi_h^{n+1}-\phi_h^n\right)\xi_h+\Delta t\int_{\Omega}\left(\mathbf{v}_h^{n+1} \cdot \nabla \phi_h^n \right) \xi_h+\Delta t\int_{\Omega}b\left(\phi_h^n\right)\nabla \mu_h^{n+1} \cdot \nabla \xi_h =0,\\ \\
\displaystyle \int_{\Omega}\mu_h^{n+1} \chi_h=\int_{\Omega}\nabla \phi_h^{n+1} \cdot \nabla \chi_h+\int_{\Omega}\left(\psi_+^{\prime}\left(\phi_h^{n+1}\right)+\psi_-^{\prime}\left(\phi_h^{n}\right)\right)\chi_h\\
 \displaystyle+\int_{\Omega}\left(f_+^{\prime}(\phi_h^{n+1})g(\mathbf{F}_h^{n+1})+f_-^{\prime}(\phi_h^n)g(\mathbf{F}_h^{n+1})\right)\chi_h.
\end{cases}
\end{equation}
\begin{rem}
 \label{rem:a2hgen}
 The approximation scheme \eqref{eqn:galerkinfe1} could be straightforwardly generalized to the case in which \eqref{eqn:wassum} has the form 
 \begin{equation}
\label{eqn:wassum2}
w(\phi,\mathbf{F})=\sum_{i=1}^mf_i(\phi)g_i(\mathbf{F})+h(\mathbf{F})+m(\phi),
\end{equation}
where $m\in \mathbb{N}$, with the functions in \eqref{eqn:wassum2} satisfying the Assumption $\bf{A2_{h}}$. Assumption $\bf{A2_{h}}$ could also be relaxed by assuming the functions $g_i$ to satisfy the property that $|g_i(\mathbf{T})|\leq k_{h2,i}(1+|\mathbf{T}|^p)$, $p\in [0,6)$, for all $\mathbf{T}\in \mathbb{R}^{3\times 3}$, without constraints on their positivity, by employing a separation of the positive and negative parts of $g_i$ in \eqref{eqn:galerkinfe1}$_6$ as done for the function $f$ in \eqref{eqn:galerkinfe1}$_4$. 
\end{rem}
The existence of a solution to \eqref{eqn:galerkinfe1} and the gradient stability of the approximation scheme are given in the following lemma.
\begin{lem}
\label{lem:gradstab1}
For all $n=0,\dots,N-1$, given $\phi_h^n\in S_h$, $\mathbf{F}_h^n\in X_h$, with $\phi_h^0=P_h^S(\phi_0)$ and $\mathbf{F}_h^0=P_h^X(\mathbf{F}_0)$,
 there exists a solution $(\mathbf{v}_h^{n+1},s_h^{n+1},\mathbf{F}_h^{n+1},\mathbf{M}_h^{n+1},\phi_h^{n+1},\mu_h^{n+1})\in V_h\times R_h \times X_h\times X_h\times S_h\times S_h$ to system \eqref{eqn:galerkinfe1}, which satisfies the following stability bound:
\begin{align}
\label{eqn:stabbound1}
&\max_{n=0\to N-1}\left(\frac{1}{2}||\nabla \phi_h^{n+1}||_{L^2(\Omega;\mathbb{R}^3)}^2+\frac{\lambda}{2}||\nabla \mathbf{F}_h^{n+1}||_{L^2(\Omega;\mathbb{R}^{3\times 3 \times 3})}^2\right)+\frac{(\Delta t)^2}{2}\sum_{n=0}^{N-1}\left|\left|\frac{\nabla (\phi_h^{n+1}-\phi_h^n)}{(\Delta t)^2}\right|\right|_{L^2(\Omega;\mathbb{R}^3)}^2\\
& \notag+\frac{\lambda(\Delta t)^2}{2}\sum_{n=0}^{N-1}\left|\left|\frac{\nabla (\mathbf{F}_h^{n+1}-\mathbf{F}_h^{n})}{(\Delta t)^2}\right|\right|_{L^2(\Omega;\mathbb{R}^{3\times 3 \times 3})}^2+ \Delta t\nu \sum_{n=0}^{N-1}||\nabla \mathbf{v}_h^{n+1}||_{L^2(\Omega;\mathbb{R}^{3\times 3})}^2\\
& \notag +\Delta t\gamma \sum_{n=0}^{N-1}|| \mathbf{M}_h^{n+1}||_{L^2(\Omega;\mathbb{R}^{3\times 3})}^2+\Delta t \sum_{n=0}^{N-1}\int_{\Omega}b\left(\phi_h^n\right)\nabla \mu_h^{n+1} \cdot \nabla \mu_h^{n+1}\leq C(\phi_h^0,\mathbf{F}_h^0)+C.
\end{align}
Moreover, the function 
\begin{equation}
\label{eqn:lyapunov1}
\int_{\Omega}\left(\psi(\phi_h^{n+1})+h(\mathbf{F}_h^{n+1})+f(\phi_h^{n+1})g(\mathbf{F}_h^{n+1})\right)+\frac{1}{2}||\nabla \phi_h^{n+1}||_{L^2(\Omega;\mathbb{R}^3)}^2+\frac{\lambda}{2}||\nabla \mathbf{F}_h^{n+1}||_{L^2(\Omega;\mathbb{R}^{3\times 3 \times 3})}^2
\end{equation}
is a decreasing (Lyapunov) function for the discrete solutions.
\end{lem}
\begin{pf}
Following \cite{ghk}, we are going to prove the existence of a solution to System \eqref{eqn:galerkinfe1} by applying \cite[Lemma II.1.4]{temam77}, which relies on proper stability bounds for the discrete system. In order to proceed, we first consider a reduced expression of System \eqref{eqn:galerkinfe1}, where the incompressibility constraint \eqref{eqn:system3}$_2$ is enforced in the definition of the finite element space for the variable $\mathbf{v}_h^{n+1}$, i.e. we consider $\mathbf{v}_h^{n+1},\mathbf{u}_h\in V_{h,\text{div}}$, where
\[
V_{h,\text{div}}:=\{\mathbf{v}_h\in V_h:\int_{\Omega}\text{div}\mathbf{v}_hp_h=0 \; \forall p_h\in S_h\}.
\]
We also introduce the variables
\[
\mathbf{G}_h^{n+1}:=\mathbf{F}_h^{n+1}-\mathbf{F}_h^{n},
\]
\[
z_h^{n+1}:=\phi_h^{n+1}-\phi_h^{n},
\]
\[
y_h^{n+1}:=\mu_h^{n+1}-\frac{1}{|\Omega|}\int_{\Omega}\mu_h^{n+1}.
\]
We observe that, taking $\xi_h\equiv 1$ in \eqref{eqn:galerkinfe1}$_5$, after integration by parts in the second term on the left hand side and using $\mathbf{v}_h^{n+1}\in V_{h,\text{div}}$, we have that
\[
\int_{\Omega}\phi_h^{n+1}=\int_{\Omega}\phi_h^{n}.
\]
Hence, the variables $z_h^{n+1}$ and $y_h^{n+1}$ belong to the space
\[
S_{h,0}:=\{q_h\in S_h:\int_{\Omega}q_h=0\}.
\]
We then require equations \eqref{eqn:galerkinfe1}$_5$ and \eqref{eqn:galerkinfe1}$_6$ to be valid only for test functions $\xi_h,\chi_h\in S_{h,0}$, substituting $\phi_h^{n+1}=z_h^{n+1}+\phi_h^{n}$ and $\mu_h^{n+1}=y_h^{n+1}+\frac{1}{|\Omega|}\int_{\Omega}\mu_h^{n+1}$ in their expressions. We further substitute $\mathbf{F}_h^{n+1}=\mathbf{G}_h^{n+1}+\mathbf{F}_h^{n}$ in \eqref{eqn:galerkinfe1}$_4$ and \eqref{eqn:galerkinfe1}$_5$ and we add to \eqref{eqn:galerkinfe1}$_4$ a regularizing term $\alpha \int_{\Omega}\mathbf{G}_h^{n+1}\colon \boldsymbol{\Pi}_h$, with $\alpha>0$. The latter term is introduced to be able to control, in the forthcoming stability estimates, the mass of $\mathbf{G}_h^{n+1}$, which is not equal to zero. In a second step we will study the limit problem $\alpha\to 0$. The reduced system which we are going  to study for the time being is thus the following: 
given $\phi_h^n\in S_h$, $\mathbf{F}_h^n\in X_h$, find $(\mathbf{v}_{h,\alpha}^{n+1},\mathbf{G}_{h,\alpha}^{n+1},\mathbf{M}_{h,\alpha}^{n+1},z_{h,\alpha}^{n+1},y_{h,\alpha}^{n+1})\in V_{h,\text{div}}\times X_h\times X_h\times S_{h,0}\times S_{h,0}$ such that,
for all $(\mathbf{u}_h,\boldsymbol{\Theta}_h,\boldsymbol{\Pi}_h,\xi_h,\chi_h)\in V_{h,\text{div}}\times X_h\times X_h\times S_{h,0}\times S_{h,0}$,
\begin{equation}
\label{eqn:galerkinfe1red}
\begin{cases}
\displaystyle \nu \int_{\Omega}\nabla \mathbf{v}_{h,\alpha}^{n+1}\colon \nabla \mathbf{u}_h=\int_{\Omega}y_{h,\alpha}^{n+1} \nabla \phi_h^n \cdot \mathbf{u}_h-\int_{\Omega}\left(\mathbf{M}_{h,\alpha}^{n+1}(\mathbf{F}_h^n)^T\right)\colon \nabla\mathbf{u}_h\\
\displaystyle+\int_{\Omega}\left(\nabla \mathbf{F}_h^n\odot \mathbf{M}_{h,\alpha}^{n+1}\right)\cdot \mathbf{u}_h,\\ \\
 \displaystyle \int_{\Omega} \mathbf{G}_{h,\alpha}^{n+1} \colon \boldsymbol{\Theta}_h+\Delta t \int_{\Omega}\left({\mathbf{v}_{h,\alpha}^{n+1}} \cdot \nabla \right)\mathbf{F}_h^n\colon \boldsymbol{\Theta}_h-\Delta t\int_{\Omega}\left(\nabla {\mathbf{v}_{h,\alpha}^{n+1}}\right)\mathbf{F}_h^n\colon \boldsymbol{\Theta}_h\\
 \displaystyle+\Delta t\gamma \int_{\Omega}\mathbf{M}_{h,\alpha}^{n+1}\colon \boldsymbol{\Theta}_h=0,\\ \\
 \displaystyle \int_{\Omega}\mathbf{M}_{h,\alpha}^{n+1}\colon \boldsymbol{\Pi}_h=\alpha \int_{\Omega}\mathbf{G}_{h,\alpha}^{n+1}\colon\boldsymbol{\Pi}_h\\
 +\int_{\Omega}\left(h_+^{\prime}(\mathbf{G}_{h,\alpha}^{n+1}+\mathbf{F}_{h}^{n})+h_-^{\prime}(\mathbf{F}_h^{n})\right)\colon \boldsymbol{\Pi}_h+\int_{\Omega}[f(\phi_h^n)]_+\left(g_+^{\prime}(\mathbf{G}_{h,\alpha}^{n+1}+\mathbf{F}_{h}^{n})+g_-^{\prime}(\mathbf{F}_h^{n})\right)\colon \boldsymbol{\Pi}_h\\
 \displaystyle +\int_{\Omega}(f(\phi_h^n)-[f(\phi_h^n)]_+)\left(g_+^{\prime}(\mathbf{F}_h^{n})+g_-^{\prime}(\mathbf{G}_{h,\alpha}^{n+1}+\mathbf{F}_{h}^{n})\right)\colon \boldsymbol{\Pi}_h +\lambda \int_{\Omega}\nabla \left(\mathbf{G}_{h,\alpha}^{n+1}+\mathbf{F}_{h}^{n}\right) \mathbin{\tensorm} \nabla \boldsymbol{\Pi}_h,\\ \\
 \displaystyle \int_{\Omega}z_{h,\alpha}^{n+1}\xi_h+\Delta t\int_{\Omega}\left(\mathbf{v}_{h,\alpha}^{n+1} \cdot \nabla \phi_h^n \right) \xi_h+\Delta t\int_{\Omega}b\left(\phi_h^n\right)\nabla y_{h,\alpha}^{n+1} \cdot \nabla \xi_h =0,\\ \\
\displaystyle \int_{\Omega}y_{h,\alpha}^{n+1} \chi_h=\int_{\Omega}\nabla \left(z_{h,\alpha}^{n+1}+\phi_{h}^{n}\right) \cdot \nabla \chi_h+\int_{\Omega}\left(\psi_+^{\prime}\left(z_{h,\alpha}^{n+1}+\phi_{h}^{n}\right)+\psi_-^{\prime}\left(\phi_h^{n}\right)\right)\chi_h\\
 \displaystyle+\int_{\Omega}\left(f_+^{\prime}(z_{h,\alpha}^{n+1}+\phi_{h}^{n})g(\mathbf{G}_{h,\alpha}^{n+1}+\mathbf{F}_{h}^{n})+f_-^{\prime}(\phi_h^n)g(\mathbf{G}_{h,\alpha}^{n+1}+\mathbf{F}_{h}^{n})\right)\chi_h,
\end{cases}
\end{equation}
where the subscript $\alpha$ indicates the dependence of the solutions of \eqref{eqn:galerkinfe1red} on the parameter $\alpha$.
We now introduce the finite dimensional Hilbert space
\[
X:=V_{h,\text{div}}\times X_h\times X_h\times S_{h,0}\times S_{h,0},
\]
endowed with the inner product
\begin{align*}
&((\mathbf{v}_1,\mathbf{G}_1,\mathbf{M}_1,z_1,y_1),(\mathbf{v}_2,\mathbf{G}_2,\mathbf{M}_2,z_2,y_2))_X\\
&:=\int_{\Omega}\nabla \mathbf{v}_1 \colon \nabla \mathbf{v}_2 +\int_{\Omega}\mathbf{G}_1 \colon \mathbf{G}_2+\int_{\Omega}\mathbf{M}_1 \colon \mathbf{M}_2+\int_{\Omega}\nabla z_1 \cdot \nabla z_2+\int_{\Omega}\nabla y_1 \cdot \nabla y_2,
\end{align*}
and the corresponding induced norm $||\cdot||_X^2:=(\cdot,\cdot)_X$. \\
\noindent
We define the continuous map $\mathcal{P}:X\rightarrow X$ which associates to $(\mathbf{v},\mathbf{G},\mathbf{M},z,y)\in X$ an element $\mathcal{P}(X)\in X$ such that, for all $(\bar{\mathbf{v}},\bar{\mathbf{G}},\bar{\mathbf{M}},\bar{z},\bar{y})\in X$,
\begin{align}
\label{eqn:mapx}
&\displaystyle (\mathcal{P}(\mathbf{v},\mathbf{G},\mathbf{M},z,y),(\bar{\mathbf{v}},\bar{\mathbf{G}},\bar{\mathbf{M}},\bar{z},\bar{y}))_X=
\Delta t \nu \int_{\Omega}\nabla \mathbf{v}\colon \nabla \bar{\mathbf{v}}-\Delta t
\int_{\Omega}y \nabla \phi_h^n \cdot \bar{\mathbf{v}}+\Delta t\int_{\Omega}\left(\mathbf{M}(\mathbf{F}_h^n)^T\right)\colon \nabla\bar{\mathbf{v}}\\
&\displaystyle \notag -\Delta t\int_{\Omega}\left(\nabla \mathbf{F}_h^n\right)^T\mathbf{M}\cdot \bar{\mathbf{v}}+\int_{\Omega} \mathbf{G} \colon \bar{\mathbf{M}}+\Delta t \int_{\Omega}\left({\mathbf{v}} \cdot \nabla \right)\mathbf{F}_h^n\colon \bar{\mathbf{M}}-\Delta t\int_{\Omega}\left(\nabla \mathbf{v}\right)\mathbf{F}_h^n\colon \bar{\mathbf{M}}\\
&\displaystyle \notag
+\Delta t\gamma \int_{\Omega}\mathbf{M}\colon \bar{\mathbf{M}}
 -\int_{\Omega}\mathbf{M}\colon \bar{\mathbf{G}}+\alpha \int_{\Omega}\mathbf{G}\colon\bar{\mathbf{G}}+\int_{\Omega}\left(h_+^{\prime}(\mathbf{G}+\mathbf{F}_h^n)+h_-^{\prime}(\mathbf{F}_h^{n})\right)\colon \bar{\mathbf{G}}\\
  &\displaystyle \notag +\int_{\Omega}[f(\phi_h^n)]_+\left(g_+^{\prime}(\mathbf{G}+\mathbf{F}_h^n)+g_-^{\prime}(\mathbf{F}_h^{n})\right)\colon \bar{\mathbf{G}}+\int_{\Omega}(f(\phi_h^n)-[f(\phi_h^n)]_+)\left(g_+^{\prime}(\mathbf{F}_h^{n})+g_-^{\prime}(\mathbf{G}+\mathbf{F}_h^n)\right)\colon \bar{\mathbf{G}} \\
  &\displaystyle \notag +\lambda \int_{\Omega}\nabla \left(\mathbf{G}+\mathbf{F}_h^n\right) \mathbin{\tensorm} \nabla \bar{\mathbf{G}}+\int_{\Omega}z\,\bar{y}+\Delta t\int_{\Omega}\left(\mathbf{v} \cdot \nabla \phi_h^n \right) \bar{y}+\Delta t\int_{\Omega}b\left(\phi_h^n\right)\nabla y \cdot \nabla \bar{y}\\
&\displaystyle \notag -\int_{\Omega}y \bar{z}+\int_{\Omega}\nabla \left(z+\phi_h^n\right) \cdot \nabla \bar{z}+\int_{\Omega}\left(\psi_+^{\prime}\left(z+\phi_h^n\right)+\psi_-^{\prime}\left(\phi_h^{n}\right)\right)\bar{z}\\
&\displaystyle \notag +\int_{\Omega}\left(f_+^{\prime}(z+\phi_h^n)g(\mathbf{G}+\mathbf{F}_h^n)+f_-^{\prime}(\phi_h^n)g(\mathbf{G}+\mathbf{F}_h^n)\right)\bar{z}.
\end{align}
A zero of the map $\mathcal{P}$, if it exists, is a solution to System \eqref{eqn:galerkinfe1red}. In the following we are going to prove that, if $||(\mathbf{v},\mathbf{G},\mathbf{M},z,y)||_X=R>0$ is large enough, then $(\mathcal{P}(\mathbf{v},\mathbf{G},\mathbf{M},z,y),(\mathbf{v},\mathbf{G},\mathbf{M},z,y))_X>0$. Then,
\cite[Lemma II.1.4]{temam77} implies that $\mathcal{P}$ admits a root $(\mathbf{v}^*,\mathbf{G}^*,\mathbf{M}^*,z^*,y^*)\in X$, with $||(\mathbf{v}^*,\mathbf{G}^*,\mathbf{M}^*,z^*,y^*)||_X\leq R$.
Indeed, we have that
\begin{align*}
&\displaystyle (\mathcal{P}(\mathbf{v},\mathbf{G},\mathbf{M},z,y),(\mathbf{v},\mathbf{G},\mathbf{M},z,y))_X=
\Delta t \nu \int_{\Omega}\nabla \mathbf{v}\colon \nabla {\mathbf{v}}
+\Delta t\gamma \int_{\Omega}\mathbf{M}\colon {\mathbf{M}}
 +\alpha \int_{\Omega}\mathbf{G}\colon{\mathbf{G}}\\
&\displaystyle+\int_{\Omega}\left(h_+^{\prime}(\mathbf{G}+\mathbf{F}_h^n)+h_-^{\prime}(\mathbf{F}_h^{n})\right)\colon \left({\mathbf{G}}+\mathbf{F}_h^{n}-\mathbf{F}_h^{n}\right)\\
&\displaystyle +\int_{\Omega}[f(\phi_h^n)]_+\left(g_+^{\prime}(\mathbf{G}+\mathbf{F}_h^n)+g_-^{\prime}(\mathbf{F}_h^{n})\right)\colon \left({\mathbf{G}}+\mathbf{F}_h^{n}-\mathbf{F}_h^{n}\right)\\
  &\displaystyle +\int_{\Omega}(f(\phi_h^n)-[f(\phi_h^n)]_+)\left(g_+^{\prime}(\mathbf{F}_h^{n})+g_-^{\prime}(\mathbf{G}+\mathbf{F}_h^n)\right)\colon \left({\mathbf{G}}+\mathbf{F}_h^{n}-\mathbf{F}_h^{n}\right) \\
  & \displaystyle +\lambda \int_{\Omega}\nabla \left(\mathbf{G}+\mathbf{F}_h^n\right) \mathbin{\tensorm} \nabla {\mathbf{G}}+\Delta t\int_{\Omega}b\left(\phi_h^n\right)\nabla y \cdot \nabla {y}+\int_{\Omega}\nabla \left(z+\phi_h^n\right) \cdot \nabla {z}\\
  &\displaystyle +\int_{\Omega}\left(\psi_+^{\prime}\left(z+\phi_h^n\right)+\psi_-^{\prime}\left(\phi_h^{n}\right)\right)(z+\phi_h^n-\phi_h^n)\\
  &\displaystyle +\int_{\Omega}\left(f_+^{\prime}(z+\phi_h^n)g(\mathbf{G}+\mathbf{F}_h^{n})+f_-^{\prime}(\phi_h^n)g(\mathbf{G}+\mathbf{F}_h^{n})\right)\left(z+\phi_h^n-\phi_h^n\right).
\end{align*}
Given the facts that $f_+,g_+,\psi_+$ are convex functions and $f_-,g_-,\psi_-$ are concave functions of their arguments, given the assumption of the positivity of $g$, using Assumption \textbf{A1} and the relation
$a\star (a+b)=\frac{1}{2}|a|^2+\frac{1}{2}|a+b|^2-\frac{1}{2}|b|^2$ (with $a,b$ any scalar, vector or matrix elements with the corresponding scalar product indicated by the symbol $\star$), we get that
\begin{align}
\label{eqn:mapxcalc1}
&\displaystyle \Delta t\nu ||\nabla \mathbf{v}||_{L^2(\Omega;\mathbb{R}^{3\times 3})}^2+\Delta t\gamma || \mathbf{M}||_{L^2(\Omega;\mathbb{R}^{3\times 3})}^2+\alpha ||\mathbf{G}||_{L^2(\Omega;\mathbb{R}^{3\times 3})}^2\\
&\displaystyle \notag +\int_{\Omega}\left(h(\mathbf{G}+\mathbf{F}_h^n)+f(\phi_h^n)g(\mathbf{G}+\mathbf{F}_h^n)+f(z+\phi_h^n)g(\mathbf{G}+\mathbf{F}_h^n)+\psi\left(z+\phi_h^n\right)\right)\\
&\displaystyle \notag +\frac{\lambda}{2} ||\nabla \mathbf{G}||_{L^2(\Omega;\mathbb{R}^{3\times 3 \times 3})}^2+\Delta t b_0||\nabla y||_{L^2(\Omega;\mathbb{R}^{3})}^2+\frac{1}{2}||\nabla z||_{L^2(\Omega;\mathbb{R}^{3})}^2\leq (\mathcal{P}(\mathbf{v},\mathbf{G},\mathbf{M},z,y),(\mathbf{v},\mathbf{G},\mathbf{M},z,y))_X\\
&\displaystyle \notag  +\int_{\Omega}\left(h(\mathbf{F}_h^n)+f(\phi_h^n)g(\mathbf{F}_h^n)+f(\phi_h^n)g(\mathbf{G}+\mathbf{F}_h^n)+\psi\left(\phi_h^n\right)\right)+\frac{\lambda}{2} ||\nabla \mathbf{F}_h^n||_{L^2(\Omega;\mathbb{R}^{3\times 3 \times 3})}^2+\frac{1}{2}||\nabla \phi_h^n||_{L^2(\Omega;\mathbb{R}^{3})}^2.
\end{align}
Now, using Assumptions $\bf{A2_h}$ and \textbf{A3}, we get that
\begin{align}
\label{eqn:mapxcalc2}
&\displaystyle \Delta t\nu ||\nabla \mathbf{v}||_{L^2(\Omega;\mathbb{R}^{3\times 3})}^2+\Delta t\gamma || \mathbf{M}||_{L^2(\Omega;\mathbb{R}^{3\times 3})}^2+\alpha ||\mathbf{G}||_{L^2(\Omega;\mathbb{R}^{3\times 3})}^2+\frac{\lambda}{2} ||\nabla \mathbf{G}||_{L^2(\Omega;\mathbb{R}^{3\times 3 \times 3})}^2\\
&\displaystyle \notag +\Delta t b_0||\nabla y||_{L^2(\Omega;\mathbb{R}^{3})}^2+\frac{1}{2}||\nabla z||_{L^2(\Omega;\mathbb{R}^{3})}^2-C\left(\mathbf{F}_h^n,\phi_h^n\right)\leq (\mathcal{P}(\mathbf{v},\mathbf{G},\mathbf{M},z,y),(\mathbf{v},\mathbf{G},\mathbf{M},z,y))_X,
\end{align}
where $C\left(\mathbf{F}_h^n,\phi_h^n\right)$ is a positive constant which depends on $\mathbf{F}_h^n$ and $\phi_h^n$. Hence,
\[
(\mathcal{P}(\mathbf{v},\mathbf{G},\mathbf{M},z,y),(\mathbf{v},\mathbf{G},\mathbf{M},z,y))_X\geq C||(\mathbf{v},\mathbf{G},\mathbf{M},z,y)||_X-C\left(\mathbf{F}_h^n,\phi_h^n\right)>0
\]
if $||(\mathbf{v},\mathbf{G},\mathbf{M},z,y)||_X\geq R$ and $R$ is large enough, and \cite[Lemma II.1.4]{temam77} implies that $\mathcal{P}$ admits a root $(\mathbf{v}^*,\mathbf{G}^*,\mathbf{M}^*,z^*,y^*)\in X$, which is a solution of \eqref{eqn:galerkinfe1red}. 

With the aim of recovering a solution for the original System \eqref{eqn:galerkinfe1}, we take the limit in \eqref{eqn:galerkinfe1red} as $\alpha\to 0$. We thus need to obtain a-priori estimates for system \eqref{eqn:galerkinfe1red} which are uniform in the parameter $\alpha$. 
Let us take $\mathbf{u}_h=\Delta t \mathbf{v}_{h,\alpha}^{n+1}$, $\boldsymbol{\Theta}_h=\mathbf{M}_{h,\alpha}^{n+1}$, $\boldsymbol{\Pi}_h=-\mathbf{G}_{h,\alpha}^{n+1}$, $\xi_h=y_{h,\alpha}^{n+1}$ and $\chi_h=-z_{h,\alpha}^{n+1}$ in \eqref{eqn:galerkinfe1red}. Summing all the equations, using again the identity $a\star (a+b)=\frac{1}{2}|a|^2+\frac{1}{2}|a+b|^2-\frac{1}{2}|b|^2$ (with $a,b$ any scalar, vector or matrix elements with the corresponding scalar product indicated by the symbol $\star$), the convexity of the functions $f_+,g_+,\psi_+$ and the concavity of the functions $f_-,g_-,\psi_-$, given also the assumption of the positivity of $g$, we get, similarly to \eqref{eqn:mapxcalc1}, that
\begin{align}
\label{eqn:discalpha2}
&\displaystyle\Delta t\nu ||\nabla \mathbf{v}_{h,\alpha}^{n+1}||_{L^2(\Omega;\mathbb{R}^{3\times 3})}^2+\Delta t\gamma || \mathbf{M}_{h,\alpha}^{n+1}||_{L^2(\Omega;\mathbb{R}^{3\times 3})}^2+\alpha ||\mathbf{G}_{h,\alpha}^{n+1}||_{L^2(\Omega;\mathbb{R}^{3\times 3})}^2+\frac{\lambda}{2}||\nabla \mathbf{G}_{h,\alpha}^{n+1}||_{L^2(\Omega;\mathbb{R}^{3\times 3 \times 3})}^2\\
& \displaystyle \notag +\Delta t b_0||\nabla y_{h,\alpha}^{n+1}||_{L^2(\Omega;\mathbb{R}^{3})}^2 +\frac{1}{2}||\nabla z_{h,\alpha}^{n+1}||_{L^2(\Omega;\mathbb{R}^{3})}^2\leq C(\phi_h^{n},\mathbf{F}_h^{n}),
\end{align}
uniformly in $\alpha$.
Moreover, similarly to \eqref{eqn:fmean}, taking $\boldsymbol{\Theta}_h=\mathbf{e}_{l,r}$, $l,r=1,2,3$, in \eqref{eqn:galerkinfe1red}$_2$ and using the fact that $\mathbf{v}_{h,\alpha}^{n+1}\in V_{h,\text{div}}$ and \eqref{eqn:discalpha2}, we obtain that
\begin{equation}
\label{eqn:discalpha3}
\left|\int_{\Omega}\mathbf{G}_{h,\alpha lr}^{n+1}\right|=\left|\Delta t \gamma \int_{\Omega}\mathbf{M}_{h,\alpha lr}^{n+1}\right|\leq C|| \mathbf{M}_{h,\alpha}^{n+1}||_{L^2(\Omega;\mathbb{R}^{3\times 3})}\leq  C(\phi_h^{n},\mathbf{F}_h^{n}),
\end{equation}
uniformly in $\alpha$. From \eqref{eqn:discalpha2} and \eqref{eqn:discalpha3} we conclude, thanks to the Poincar\'e--Wirtinger inequality, that
\begin{equation}
\label{eqn:discalpha4}
||\mathbf{G}_{h,\alpha}^{n+1}||_{H^1(\Omega;\mathbb{R}^{3\times 3})}^2\leq C(\phi_h^{n},\mathbf{F}_h^{n}),
\end{equation}
uniformly in $\alpha$. Then, thanks to the Bolzano--Weierstrass theorem and \eqref{eqn:discalpha2}-\eqref{eqn:discalpha4}, given any sequence $\alpha \to 0$, we can identify a subsequence $\alpha' \to 0$ and a limit point \[(\mathbf{v}_{h,\alpha'}^{n+1},\mathbf{G}_{h,\alpha'}^{n+1},\mathbf{M}_{h,\alpha'}^{n+1},z_{h,\alpha'}^{n+1},y_{h,\alpha'}^{n+1})\to (\mathbf{v}_{h}^{n+1},\mathbf{G}_{h}^{n+1},\mathbf{M}_{h}^{n+1},z_{h}^{n+1},y_{h}^{n+1})\]
in $X$ which satisfies \eqref{eqn:galerkinfe1red} as $\alpha' \to 0$.

Finally, we define 
\[(\mathbf{v}_{h}^{n+1},\mathbf{F}_{h}^{n+1},\mathbf{M}_{h}^{n+1},\phi_{h}^{n+1},\mu_{h}^{n+1}):=(\mathbf{v}_{h}^{n+1},\mathbf{G}_{h}^{n+1}+\mathbf{F}_{h}^{n},\mathbf{M}_{h}^{n+1},z_{h}^{n+1}+\phi_h^n,y_{h}^{n+1}+\bar{\rho}),\]
with 
\[\bar{\rho}:=\frac{1}{|\Omega|}\left(\int_{\Omega}\left(\psi_+^{\prime}\left(\phi_{h}^{n+1}\right)+\psi_-^{\prime}\left(\phi_h^{n}\right)\right)+\int_{\Omega}\left(f_+^{\prime}(\phi_{h}^{n+1})g(\mathbf{F}_{h}^{n+1})+f_-^{\prime}(\phi_h^n)g(\mathbf{F}_{h}^{n+1})\right)\right).\]
Observing moreover that \eqref{eqn:galerkinfe1red}$_1$ defines a linear functional from $V_h$ to $\mathbb{R}$ which vanishes on $V_{h,\text{div}}$, we conclude by standard arguments (see e.g. \cite[Section 1.2, Chapter III]{gr}) that there exists a unique $s_h^{n+1}\in R_h$ such that 
 \[
 (\mathbf{v}_{h}^{n+1},s_h^{n+1},\mathbf{F}_{h}^{n+1},\mathbf{M}_{h}^{n+1},\phi_{h}^{n+1},\mu_{h}^{n+1})\in V_h\times R_h\times X_h \times X_h \times S_h \times X_h
 \]
 is a solution of System \eqref{eqn:galerkinfe1}.

Taking $\mathbf{u}_h=\Delta t \mathbf{v}_h^{n+1}$, $p_h=\Delta t s_h^{n+1}$, $\boldsymbol{\Theta}_h=\mathbf{M}_h^{n+1}$, $\boldsymbol{\Pi}_h=-(\mathbf{F}_h^{n+1}-\mathbf{F}_h^{n})$, $\xi_h=\mu_h^{n+1}$ and $\chi_h=-(\phi_h^{n+1}-\phi_h^n)$ in \eqref{eqn:galerkinfe1}, with similar calculations as in \eqref{eqn:discalpha2} we obtain that
\begin{align}
\label{eqn:discreteenergy1}
& \displaystyle \int_{\Omega}\left(f(\phi_h^{n+1})g(\mathbf{F}_h^{n+1})+\psi\left(\phi_h^{n+1}\right)+h\left(\mathbf{F}_h^{n+1}\right)\right)+\frac{1}{2}||\nabla \phi_h^{n+1}||_{L^2(\Omega;\mathbb{R}^{3})}^2+\frac{\lambda}{2}||\nabla \mathbf{F}_h^{n+1}||_{L^2(\Omega;\mathbb{R}^{3\times 3 \times 3})}^2\\
& \notag \displaystyle +\frac{1}{2}\left|\left|\nabla \left(\phi_h^{n+1}-\phi_h^{n}\right)\right|\right|_{L^2(\Omega;\mathbb{R}^{3})}^2+\frac{\lambda}{2}\left|\left|\nabla \left(\mathbf{F}_h^{n+1}-\mathbf{F}_h^{n}\right)\right|\right|_{L^2(\Omega;\mathbb{R}^{3\times 3 \times 3})}^2\\
& \notag \displaystyle +\Delta t\nu ||\nabla \mathbf{v}_h^{n+1}||_{L^2(\Omega;\mathbb{R}^{3\times 3})}^2+\Delta t\gamma || \mathbf{M}_h^{n+1}||_{L^2(\Omega;\mathbb{R}^{3\times 3})}^2+\Delta t \int_{\Omega}b\left(\phi_h^n\right)\nabla \mu_h^{n+1}\cdot \nabla \mu_h^{n+1}\\
& \notag \displaystyle  \leq \int_{\Omega}\left(f(\phi_h^{n})g(\mathbf{F}_h^{n})+\psi\left(\phi_h^{n}\right)+h\left(\mathbf{F}_h^{n}\right)\right)+\frac{1}{2}||\nabla \phi_h^{n}||_{L^2(\Omega;\mathbb{R}^{3})}^2+\frac{\lambda}{2}||\nabla \mathbf{F}_h^{n}||_{L^2(\Omega;\mathbb{R}^{3\times 3 \times 3})}^2,
\end{align}
which gives that \eqref{eqn:lyapunov1} is a decreasing (Lyapunov) function for the discrete solutions. Summing \eqref{eqn:discreteenergy1} from $n=0\to m$, for $m=0\to N-1$, and using Assumptions \textbf{A2,A3}, we finally get \eqref{eqn:stabbound1}.
\end{pf}
\begin{rem}
\label{rem:galerkin1dt}
We observe that System \eqref{eqn:galerkinfe1} admits a solution and is unconditionally gradient stable for any value of $\Delta t>0$, with no requirements on the smallness of $\Delta t$ (with respect to $h$ and the model parameters). 
\end{rem}
\begin{rem}
\label{rem:convdecoupled}
System \eqref{eqn:galerkinfe1} is fully coupled and nonlinear, and can be solved e.g. by means of a Newton method, which at each iteration requires the solution of the full tangent algebraic system. In order to decrease the computational demand for the solution of System \eqref{eqn:galerkinfe1}, we practically solve a fixed-point iteration scheme which decouples \eqref{eqn:galerkinfe1}$_1$-\eqref{eqn:galerkinfe1}$_2$, \eqref{eqn:galerkinfe1}$_3$-\eqref{eqn:galerkinfe1}$_4$ and \eqref{eqn:galerkinfe1}$_5$-\eqref{eqn:galerkinfe1}$_6$, i.e. we solve the following problem:
\\ \\
\noindent
\textbf{Problem} $\mathbf{P_h^{ID}}$: for $n=0,\dots,N-1$, given $\phi_h^n\in S_h$, $\mathbf{F}_h^n\in X_h$, and for $k=0,\dots,K-1$, given $\phi_h^{n+1,0}=\phi_h^{n}$, $\mathbf{F}_h^{n+1,0}=\mathbf{F}_h^{n}$, $\mu_h^{n+1,0}=\mu_h^{n}$, $\mathbf{M}_h^{n+1,0}=\mathbf{M}_h^{n}$,
 find 
\[(\mathbf{v}_h^{n+1,k+1},s_h^{n+1,k+1},\mathbf{F}_h^{n+1,k+1},\mathbf{M}_h^{n+1,k+1},\phi_h^{n+1,k+1},\mu_h^{n+1,k+1})\in V_h\times R_h \times X_h\times X_h\times S_h\times S_h\] 
such that,
for all $(\mathbf{u}_h,p_h,\boldsymbol{\Theta}_h,\boldsymbol{\Pi}_h,\xi_h,\chi_h)\in V_h\times R_h \times X_h\times X_h\times S_h\times S_h$,
\begin{equation}
\label{eqn:galerkinfe1dec}
\begin{cases}
\displaystyle \nu \int_{\Omega}\nabla \mathbf{v}_h^{n+1,k+1}\colon \nabla \mathbf{u}_h-\int_{\Omega}s_h^{n+1,k+1}\diver \mathbf{u}_h =\int_{\Omega}\mu_h^{n+1,k} \nabla \phi_h^n \cdot \mathbf{u}_h-\int_{\Omega}\left(\mathbf{M}_h^{n+1,k}(\mathbf{F}_h^n)^T\right)\colon \nabla\mathbf{u}_h\\
\displaystyle+\int_{\Omega}\left(\nabla \mathbf{F}_h^n\odot \mathbf{M}_h^{n+1,k}\right)\cdot \mathbf{u}_h,\\ \\
\displaystyle \int_{\Omega}p_h\diver \mathbf{v}_h^{n+1,k+1}=0,\\ \\
 \displaystyle \int_{\Omega} \left(\mathbf{F}_h^{n+1,k+1}-\mathbf{F}_h^{n}\right) \colon \boldsymbol{\Theta}_h+\Delta t \int_{\Omega}\left({\mathbf{v}_h^{n+1,k+1}} \cdot \nabla \right)\mathbf{F}_h^n\colon \boldsymbol{\Theta}_h-\Delta t\int_{\Omega}\left(\nabla {\mathbf{v}_h^{n+1,k+1}}\right)\mathbf{F}_h^n\colon \boldsymbol{\Theta}_h\\
 \displaystyle+\Delta t\gamma \int_{\Omega}\mathbf{M}_h^{n+1,k+1}\colon \boldsymbol{\Theta}_h=0,\\ \\
 \displaystyle \int_{\Omega}\mathbf{M}_h^{n+1,k+1}\colon \boldsymbol{\Pi}_h=\int_{\Omega}\left(h_+^{\prime}(\mathbf{F}_h^{n+1,k+1})+h_-^{\prime}(\mathbf{F}_h^{n})\right)\colon \boldsymbol{\Pi}_h\\
 \displaystyle +\int_{\Omega}[f(\phi_h^n)]_+\left(g_+^{\prime}(\mathbf{F}_h^{n+1,k+1})+g_-^{\prime}(\mathbf{F}_h^{n})\right)\colon \boldsymbol{\Pi}_h\\
 \displaystyle +\int_{\Omega}(f(\phi_h^n)-[f(\phi_h^n)]_+)\left(g_+^{\prime}(\mathbf{F}_h^{n})+g_-^{\prime}(\mathbf{F}_h^{n+1,k+1})\right)\colon \boldsymbol{\Pi}_h +\lambda \int_{\Omega}\nabla \mathbf{F}_h^{n+1,k+1} \mathbin{\tensorm} \nabla \boldsymbol{\Pi}_h,\\ \\
 \displaystyle \int_{\Omega}\left(\phi_h^{n+1,k+1}-\phi_h^n\right)\xi_h+\Delta t\int_{\Omega}\left(\mathbf{v}_h^{n+1,k+1} \cdot \nabla \phi_h^n \right) \xi_h+\Delta t\int_{\Omega}b\left(\phi_h^n\right)\nabla \mu_h^{n+1,k+1} \cdot \nabla \xi_h =0,\\ \\
\displaystyle \int_{\Omega}\mu_h^{n+1,k+1} \chi_h=\int_{\Omega}\nabla \phi_h^{n+1,k+1} \cdot \nabla \chi_h+\int_{\Omega}\left(\psi_+^{\prime}\left(\phi_h^{n+1,k+1}\right)+\psi_-^{\prime}\left(\phi_h^{n}\right)\right)\chi_h\\
 \displaystyle+\int_{\Omega}\left(f_+^{\prime}(\phi_h^{n+1,k+1})g(\mathbf{F}_h^{n+1,k+1})+f_-^{\prime}(\phi_h^n)g(\mathbf{F}_h^{n+1,k+1})\right)\chi_h,
\end{cases}
\end{equation}
and set
\begin{align}
\label{eqn:galerkinfe1dec1}
& \displaystyle \mathbf{v}_h^{n+1}=\mathbf{v}_h^{n+1,K+1}, \, s_h^{n+1}=s_h^{n+1,K+1}, \, \mathbf{F}_h^{n+1}=\mathbf{F}_h^{n+1,K+1}, \, \mathbf{M}_h^{n+1}=\mathbf{M}_h^{n+1,K+1},\\ & \displaystyle \notag \phi_h^{n+1}=\phi_h^{n+1,K+1}, \, \mu_h^{n+1}=\mu_h^{n+1,K+1}.
\end{align}
The value of $K$ is chosen as the first value at which
\begin{align*}
& ||\mathbf{v}_h^{n+1,K+1}-\mathbf{v}_h^{n+1,K}||_{L^{\infty}(\Omega,\mathbb{R}^3)}+||s_h^{n+1,K+1}-s_h^{n+1,K}||_{L^{\infty}(\Omega)}\\
&+||\mathbf{F}_h^{n+1,K+1}-\mathbf{F}_h^{n+1,K}||_{L^{\infty}(\Omega,\mathbb{R}^{3\times 3})}+||\mathbf{F}_h^{n+1,K+1}-\mathbf{F}_h^{n+1,K}||_{L^{\infty}(\Omega,\mathbb{R}^{3\times 3})}\\
&+||\phi_h^{n+1,K+1}-\phi_h^{n+1,K}||_{L^{\infty}(\Omega)}+||\mu_h^{n+1,K+1}-\mu_h^{n+1,K}||_{L^{\infty}(\Omega)}<\text{tol},
\end{align*}
where $\text{tol}>0$ is a small tolerance for the convergence of the algorithm. System \eqref{eqn:galerkinfe1dec} is decoupled and can be solved, at each iteration $k$, in the following order:
\newpage
\begin{itemize}
\item[- \textbf{Step 1}] Given $\mu_h^{n+1,k}$ and $\mathbf{M}_h^{n+1,k}$, solve \eqref{eqn:galerkinfe1dec}$_1$-\eqref{eqn:galerkinfe1dec}$_2$, which is a linear saddle-point problem;
\item[- \textbf{Step 2}] Given $\mathbf{v}_h^{n+1,k+1}$ from Step 1, solve \eqref{eqn:galerkinfe1dec}$_3$-\eqref{eqn:galerkinfe1dec}$_4$ by means of a Newton method;
\item[- \textbf{Step 3}] Given $\mathbf{F}_h^{n+1,k+1}$ from Step 2, solve \eqref{eqn:galerkinfe1dec}$_5$-\eqref{eqn:galerkinfe1dec}$_6$ by means of a Newton method.
\end{itemize}
System \eqref{eqn:galerkinfe1dec} defines a map $\mathcal{M}:S_h\times X_h\rightarrow S_h\times X_h$, $\mathcal{M}(\mu_h^{n+1,k},\mathbf{M}_h^{n+1,k})=(\mu_h^{n+1,k+1},\mathbf{M}_h^{n+1,k+1})$. A fixed point of $\mathcal{M}$ is a solution of System \eqref{eqn:galerkinfe1}, which exists thanks to Lemma \ref{lem:gradstab1}. The convergence of the iterations \eqref{eqn:galerkinfe1dec}, for $k\geq 0$, together with the uniqueness of the fixed point, could be proved under proper smallness conditions on $\Delta t$, with respect to $h$ and to the model parameters, which guarantee that $\mathcal{M}$ is a contraction.
\end{rem}
\textbf{Problems} $\mathbf{P_h^I}$ and $\mathbf{P_h^{ID}}$ are gradient stable only if Assumption \eqref{eqn:wassum} is satisfied.
In order to derive a more general approximation scheme, which is also decoupled as $\mathbf{P_h^{ID}}$, we design a second approximation scheme, based on the fully decoupled scalar auxiliary variable scheme recently introduced in \cite{lishen} for the Cahn--Hilliard Navier--Stokes system and generalized here to our model. The latter scheme will be unconditionally gradient stable given a general form of $w(\phi,\mathbf{F})$.
\noindent
We thus introduce the auxiliary variable
\[
\beta:=\sqrt{\int_{\Omega}(j(\phi,\mathbf{F}))+k}, \quad j(\phi,\mathbf{F}):=\psi(\phi)+w(\phi,\mathbf{F})
\]
with $k>c_2+d_1$, so that $\beta>0$. We rewrite system \eqref{eqn:system3} as
\begin{equation}
\label{eqn:system3beta}
\begin{cases}
\displaystyle-\nu\Delta \mathbf{v}+\nabla s=\frac{\beta}{\sqrt{\int_{\Omega}(j(\phi,\mathbf{F}))+k}}\mu \nabla \phi +\frac{\beta}{\sqrt{\int_{\Omega}(j(\phi,\mathbf{F}))+k}}\diver\left(\mathbf{M}\mathbf{F}^T\right)\\
\displaystyle+\frac{\beta}{\sqrt{\int_{\Omega}(j(\phi,\mathbf{F}))+k}}\nabla \mathbf{F}\odot \mathbf{M},\\ \\
\diver\mathbf{v}=0,\\ \\
\displaystyle\frac{\partial \mathbf{F}}{\partial t}+\frac{\beta}{\sqrt{\int_{\Omega}(j(\phi,\mathbf{F}))+k}}\left(\mathbf{v}\cdot \nabla\right)\mathbf{F}-\frac{\beta}{\sqrt{\int_{\Omega}(j(\phi,\mathbf{F}))+k}}(\nabla \mathbf{v})\mathbf{F}+\gamma \mathbf{M}=0,\\ \\
\displaystyle\mathbf{M}=\frac{\beta}{\sqrt{\int_{\Omega}(j(\phi,\mathbf{F}))+k}}\partial_{\mathbf{F}}w(\phi,\mathbf{F})-\lambda \Delta \mathbf{F},\\ \\
\displaystyle \frac{\partial \phi}{\partial t}+\frac{\beta}{\sqrt{\int_{\Omega}(j(\phi,\mathbf{F}))+k}}\mathbf{v}\cdot \nabla \phi-\diver(b(\phi)\nabla \mu)=0,\\ \\
\displaystyle \mu=\frac{\beta}{\sqrt{\int_{\Omega}(j(\phi,\mathbf{F}))+k}}\psi'(\phi)- \Delta \phi+\frac{\beta}{\sqrt{\int_{\Omega}(j(\phi,\mathbf{F}))+k}}\partial_{\phi}w(\phi,\mathbf{F}),\\
\displaystyle \frac{\partial \beta}{\partial t}=\frac{1}{2\sqrt{\int_{\Omega}(j(\phi,\mathbf{F}))+k}}\left(\int_{\Omega}\partial_{\phi}j(\phi,\mathbf{F})\frac{\partial \phi}{\partial t}+\int_{\Omega}\partial_{\mathbf{F}}w(\phi,\mathbf{F})\colon \frac{\partial \mathbf{F}}{\partial t}\right).
\end{cases}
\end{equation}
Setting $\beta^0=\sqrt{\int_{\Omega}(j(\phi_0,\mathbf{F}_0))+k}$, we then consider the following fully discretized approximation of system \eqref{eqn:system3beta}:
\\ \\
\noindent
\textbf{Problem} $\mathbf{P_h^{II}}$: for $n=0,\dots,N-1$, given $\mathbf{v}_h^{n}\in V_h$, $\phi_h^n,\mu_h^n\in S_h, \beta_h^n \in \mathbb{R}, \mathbf{F}_h^n,\mathbf{M}_h^n\in X_h$, \\
find $(\mathbf{v}_h^{n+1},s_h^{n+1},\mathbf{F}_h^{n+1},\mathbf{M}_h^{n+1},\phi_h^{n+1},\mu_h^{n+1})\in V_h\times R_h \times X_h\times X_h\times S_h\times S_h$ and $\beta_h^{n+1}\in \mathbb{R}$ such that,
for all $(\mathbf{u}_h,p_h,\boldsymbol{\Theta}_h,\boldsymbol{\Pi}_h,\xi_h,\chi_h)\in V_h\times R_h \times X_h\times X_h\times S_h\times S_h$,
\begin{equation}
\label{eqn:galerkinfe2}
\begin{cases}
\displaystyle \nu \int_{\Omega}\nabla \mathbf{v}_h^{n+1}\colon \nabla \mathbf{u}_h-\int_{\Omega}s_h^{n+1}\diver \mathbf{u}_h =\frac{\beta_h^{n+1}}{\sqrt{\int_{\Omega}(j(\phi_h^n,\mathbf{F}_h^n))+k}}\int_{\Omega}\mu_h^{n} \nabla \phi_h^n \cdot \mathbf{u}_h\\
\displaystyle -\frac{\beta_h^{n+1}}{\sqrt{\int_{\Omega}(j(\phi_h^n,\mathbf{F}_h^n))+k}}\int_{\Omega}\left(\mathbf{M}_h^{n}(\mathbf{F}_h^n)^T\right)\colon \nabla\mathbf{u}_h+\frac{\beta_h^{n+1}}{\sqrt{\int_{\Omega}(j(\phi_h^n,\mathbf{F}_h^n))+k}}\int_{\Omega}\left(\nabla \mathbf{F}_h^n\odot \mathbf{M}_h^{n}\right)\cdot \mathbf{u}_h,\\ \\
\displaystyle \int_{\Omega}p_h\diver \mathbf{v}_h^{n+1}=0,\\ \\
 \displaystyle \int_{\Omega} \left(\mathbf{F}_h^{n+1}-\mathbf{F}_h^{n}\right) \colon \boldsymbol{\Theta}_h+\Delta t \frac{\beta_h^{n+1}}{\sqrt{\int_{\Omega}(j(\phi_h^n,\mathbf{F}_h^n))+k}}\int_{\Omega}\left({\mathbf{v}_h^{n}} \cdot \nabla \right)\mathbf{F}_h^n\colon \boldsymbol{\Theta}_h\\
 \displaystyle -\Delta t\frac{\beta_h^{n+1}}{\sqrt{\int_{\Omega}(j(\phi_h^n,\mathbf{F}_h^n))+k}}\int_{\Omega}\left(\nabla {\mathbf{v}_h^{n}}\right)\mathbf{F}_h^n\colon \boldsymbol{\Theta}_h+\Delta t\gamma \int_{\Omega}\mathbf{M}_h^{n+1}\colon \boldsymbol{\Theta}_h=0,\\ \\
 \displaystyle \int_{\Omega}\mathbf{M}_h^{n+1}\colon \boldsymbol{\Pi}_h=\frac{\beta_h^{n+1}}{\sqrt{\int_{\Omega}(j(\phi_h^n,\mathbf{F}_h^n))+k}}\int_{\Omega}\partial_{\mathbf{F}}w(\phi_h^n,\mathbf{F}_h^n)\colon \boldsymbol{\Pi}_h+\lambda \int_{\Omega}\nabla \mathbf{F}_h^{n+1} \mathbin{\tensorm} \nabla \boldsymbol{\Pi}_h,\\ \\
 \displaystyle \int_{\Omega}\left(\phi_h^{n+1}-\phi_h^n\right)\xi_h+\frac{\beta_h^{n+1}}{\sqrt{\int_{\Omega}(j(\phi_h^n,\mathbf{F}_h^n))+k}}\Delta t\int_{\Omega}\left(\mathbf{v}_h^{n} \cdot \nabla \phi_h^n \right) \xi_h+\Delta t\int_{\Omega}b\left(\phi_h^n\right)\nabla \mu_h^{n+1} \cdot \nabla \xi_h =0,\\ \\
\displaystyle \int_{\Omega}\mu_h^{n+1} \chi_h=\int_{\Omega}\nabla \phi_h^{n+1} \cdot \nabla \chi_h+\frac{\beta_h^{n+1}}{\sqrt{\int_{\Omega}(j(\phi_h^n,\mathbf{F}_h^n))+k}}\int_{\Omega}\partial_{\phi}j(\phi_h^n,\mathbf{F}_h^n)\chi_h,\\ \\
\displaystyle \left(\beta_h^{n+1}-\beta_h^n\right)=\frac{1}{2\sqrt{\int_{\Omega}(j(\phi_h^n,\mathbf{F}_h^n))+k}}\biggl[\left(\int_{\Omega}\partial_{\phi}j(\phi_h^n,\mathbf{F}_h^n)(\phi_h^{n+1}-\phi_h^n)\right)\\
\displaystyle+\left(\int_{\Omega}\partial_{\mathbf{F}}w(\phi_h^n,\mathbf{F}_h^n)\colon(\mathbf{F}_h^{n+1}-\mathbf{F}_h^{n})\right)+\Delta t\int_{\Omega}\left({\mathbf{v}_h^{n}} \cdot \nabla \right)\mathbf{F}_h^n\colon \mathbf{M}_h^{n+1}-\Delta t\int_{\Omega}\left(\nabla \mathbf{F}_h^n\right)^T\mathbf{M}_h^{n}\cdot \mathbf{v}_h^{n+1}\\
\displaystyle -\Delta t \int_{\Omega}\left(\nabla {\mathbf{v}_h^{n}}\right)\mathbf{F}_h^n\colon \mathbf{M}_h^{n+1}+\Delta t\int_{\Omega}\left(\mathbf{M}_h^{n}(\mathbf{F}_h^n)^T\right)\colon \nabla\mathbf{v}_h^{n+1}+\Delta t\int_{\Omega}\left(\mathbf{v}_h^{n} \cdot \nabla \phi_h^n \right) \mu_h^{n+1}\\
\displaystyle -\Delta t \int_{\Omega}\mu_h^{n} \nabla \phi_h^n \cdot \mathbf{v}_h^{n+1}\biggr].
\end{cases}
\end{equation}
We observe that \eqref{eqn:galerkinfe2} is a linear coupled system. 
The existence and uniqueness of a solution to \eqref{eqn:galerkinfe2} and its gradient stability are given in the following lemma.
\newpage
\begin{lem}
\label{lem:gradstab2}
For all $n=0,\dots,N-1$, given $\mathbf{v}_h^n\in V_h$, $\phi_h^n,\,\mu_h^n\in S_h$, $\mathbf{F}_h^n,\,\mathbf{M}_h^n\in X_h$, $\beta_h^n\in \mathbb{R}_+$, with $\phi_h^0=P_h^S(\phi_0)$, $\mathbf{F}_h^0=P_h^X(\mathbf{F}_0)$, $\beta_h^0=\sqrt{\int_{\Omega}(j(\phi_h^0,\mathbf{F}_h^0))+k}$, $\mathbf{v}_h^0=\boldsymbol{0}$, $\mu_h^0=0$, $\mathbf{M}_h^0=\boldsymbol{0}$, there exists a unique solution $(\mathbf{v}_h^{n+1},s_h^{n+1},\mathbf{F}_h^{n+1},\mathbf{M}_h^{n+1},\phi_h^{n+1},\mu_h^{n+1},\beta_h^{n+1})$ of system \eqref{eqn:galerkinfe2}, which satisfies the following stability bound:
\begin{align}
\label{eqn:stabbound2}
&\max_{n=0\to N-1}\left(\frac{1}{2}||\nabla \phi_h^{n+1}||_{L^2(\Omega;\mathbb{R}^3)}^2+\frac{\lambda}{2}||\nabla \mathbf{F}_h^{n+1}||_{L^2(\Omega;\mathbb{R}^{3\times 3 \times 3})}^2+\left(\beta_h^{n+1}\right)^2\right)\\
& \notag +\frac{(\Delta t)^2}{2}\sum_{n=0}^{N-1}\left|\left|\frac{\nabla (\phi_h^{n+1}-\phi_h^n)}{(\Delta t)^2}\right|\right|_{L^2(\Omega;\mathbb{R}^3)}^2+\frac{\lambda(\Delta t)^2}{2}\sum_{n=0}^{N-1}\left|\left|\frac{\nabla (\mathbf{F}_h^{n+1}-\mathbf{F}_h^{n})}{(\Delta t)^2}\right|\right|_{L^2(\Omega;\mathbb{R}^{3\times 3 \times 3})}^2\\
& \notag+(\Delta t)^2\sum_{n=0}^{N-1}\left|\left|\frac{(\beta_h^{n+1}-\beta_h^n)}{(\Delta t)^2}\right|\right|_{L^2(\Omega)}^2+ \Delta t\nu \sum_{n=0}^{N-1}||\nabla \mathbf{v}_h^{n+1}||_{L^2(\Omega;\mathbb{R}^{3\times 3})}^2+\Delta t\gamma \sum_{n=0}^{N-1}|| \mathbf{M}_h^{n+1}||_{L^2(\Omega;\mathbb{R}^{3\times 3})}^2\\
& \notag +\Delta t \sum_{n=0}^{N-1}\int_{\Omega}b\left(\phi_h^n\right)\nabla \mu_h^{n+1} \cdot \nabla \mu_h^{n+1}\leq C(\phi_h^0,\mathbf{F}_h^0,\beta_h^0)+C.
\end{align}
Moreover, the function 
\begin{equation}
\label{eqn:lyapunov2}
\left(\beta_h^{n+1}\right)^2+\frac{1}{2}||\nabla \phi_h^{n+1}||_{L^2(\Omega;\mathbb{R}^3)}^2+\frac{\lambda}{2}||\nabla \mathbf{F}_h^{n+1}||_{L^2(\Omega;\mathbb{R}^{3\times 3 \times 3})}^2
\end{equation}
is a decreasing (Lyapunov) function for the discrete solution.
\end{lem}
\begin{pf}
As in the proof of Lemma \ref{lem:gradstab1}, we start by considering a reduced expression of System \eqref{eqn:galerkinfe2} in which $\mathbf{v}_h^n, \mathbf{v}_h^{n+1}, \mathbf{u}_h \in V_{h,\text{div}}$, i.e. the following:
given $\mathbf{v}_h^{n}\in V_{h,\text{div}}$, $\phi_h^n,\mu_h^n\in S_h, \beta_h^n \in \mathbb{R}, \mathbf{F}_h^n,\mathbf{M}_h^n\in X_h$,
find $(\mathbf{v}_h^{n+1},\mathbf{F}_h^{n+1},\mathbf{M}_h^{n+1},\phi_h^{n+1},\mu_h^{n+1})\in V_{h,\text{div}}\times X_h\times X_h\times S_h\times S_h$ and $\beta_h^{n+1}\in \mathbb{R}$ such that,
for all $(\mathbf{u}_h,\boldsymbol{\Theta}_h,\boldsymbol{\Pi}_h,\xi_h,\chi_h)\in V_{h,\text{div}} \times X_h\times X_h\times S_h\times S_h$,
\begin{equation}
\label{eqn:galerkinfe2div}
\begin{cases}
\displaystyle \nu \int_{\Omega}\nabla \mathbf{v}_h^{n+1}\colon \nabla \mathbf{u}_h =\frac{\beta_h^{n+1}}{\sqrt{\int_{\Omega}(j(\phi_h^n,\mathbf{F}_h^n))+k}}\int_{\Omega}\mu_h^{n} \nabla \phi_h^n \cdot \mathbf{u}_h\\
\displaystyle -\frac{\beta_h^{n+1}}{\sqrt{\int_{\Omega}(j(\phi_h^n,\mathbf{F}_h^n))+k}}\int_{\Omega}\left(\mathbf{M}_h^{n}(\mathbf{F}_h^n)^T\right)\colon \nabla\mathbf{u}_h+\frac{\beta_h^{n+1}}{\sqrt{\int_{\Omega}(j(\phi_h^n,\mathbf{F}_h^n))+k}}\int_{\Omega}\left(\nabla \mathbf{F}_h^n\odot \mathbf{M}_h^{n}\right)\cdot \mathbf{u}_h,\\ \\
 \displaystyle \int_{\Omega} \left(\mathbf{F}_h^{n+1}-\mathbf{F}_h^{n}\right) \colon \boldsymbol{\Theta}_h+\Delta t \frac{\beta_h^{n+1}}{\sqrt{\int_{\Omega}(j(\phi_h^n,\mathbf{F}_h^n))+k}}\int_{\Omega}\left({\mathbf{v}_h^{n}} \cdot \nabla \right)\mathbf{F}_h^n\colon \boldsymbol{\Theta}_h\\
 \displaystyle -\Delta t\frac{\beta_h^{n+1}}{\sqrt{\int_{\Omega}(j(\phi_h^n,\mathbf{F}_h^n))+k}}\int_{\Omega}\left(\nabla {\mathbf{v}_h^{n}}\right)\mathbf{F}_h^n\colon \boldsymbol{\Theta}_h+\Delta t\gamma \int_{\Omega}\mathbf{M}_h^{n+1}\colon \boldsymbol{\Theta}_h=0,\\ \\
 \displaystyle \int_{\Omega}\mathbf{M}_h^{n+1}\colon \boldsymbol{\Pi}_h=\frac{\beta_h^{n+1}}{\sqrt{\int_{\Omega}(j(\phi_h^n,\mathbf{F}_h^n))+k}}\int_{\Omega}\partial_{\mathbf{F}}w(\phi_h^n,\mathbf{F}_h^n)\colon \boldsymbol{\Pi}_h+\lambda \int_{\Omega}\nabla \mathbf{F}_h^{n+1} \mathbin{\tensorm} \nabla \boldsymbol{\Pi}_h,\\ \\
 \displaystyle \int_{\Omega}\left(\phi_h^{n+1}-\phi_h^n\right)\xi_h+\frac{\beta_h^{n+1}}{\sqrt{\int_{\Omega}(j(\phi_h^n,\mathbf{F}_h^n))+k}}\Delta t\int_{\Omega}\left(\mathbf{v}_h^{n} \cdot \nabla \phi_h^n \right) \xi_h+\Delta t\int_{\Omega}b\left(\phi_h^n\right)\nabla \mu_h^{n+1} \cdot \nabla \xi_h =0,\\ \\
\displaystyle \int_{\Omega}\mu_h^{n+1} \chi_h=\int_{\Omega}\nabla \phi_h^{n+1} \cdot \nabla \chi_h+\frac{\beta_h^{n+1}}{\sqrt{\int_{\Omega}(j(\phi_h^n,\mathbf{F}_h^n))+k}}\int_{\Omega}\partial_{\phi}j(\phi_h^n,\mathbf{F}_h^n)\chi_h,\\ \\
\displaystyle \left(\beta_h^{n+1}-\beta_h^n\right)=\frac{1}{2\sqrt{\int_{\Omega}(j(\phi_h^n,\mathbf{F}_h^n))+k}}\biggl[\left(\int_{\Omega}\partial_{\phi}j(\phi_h^n,\mathbf{F}_h^n)(\phi_h^{n+1}-\phi_h^n)\right)\\
\displaystyle+\left(\int_{\Omega}\partial_{\mathbf{F}}w(\phi_h^n,\mathbf{F}_h^n)\colon(\mathbf{F}_h^{n+1}-\mathbf{F}_h^{n})\right)+\Delta t\int_{\Omega}\left({\mathbf{v}_h^{n}} \cdot \nabla \right)\mathbf{F}_h^n\colon \mathbf{M}_h^{n+1}-\Delta t\int_{\Omega}\left(\nabla \mathbf{F}_h^n\right)^T\mathbf{M}_h^{n}\cdot \mathbf{v}_h^{n+1}\\
\displaystyle -\Delta t \int_{\Omega}\left(\nabla {\mathbf{v}_h^{n}}\right)\mathbf{F}_h^n\colon \mathbf{M}_h^{n+1}+\Delta t\int_{\Omega}\left(\mathbf{M}_h^{n}(\mathbf{F}_h^n)^T\right)\colon \nabla\mathbf{v}_h^{n+1}+\Delta t\int_{\Omega}\left(\mathbf{v}_h^{n} \cdot \nabla \phi_h^n \right) \mu_h^{n+1}\\
\displaystyle -\Delta t \int_{\Omega}\mu_h^{n} \nabla \phi_h^n \cdot \mathbf{v}_h^{n+1}\biggr].
\end{cases}
\end{equation}
The existence of a solution to System \eqref{eqn:galerkinfe2div}, which is a finite dimensional algebraic system with the same number of equations and unknowns, is guaranteed by its linearity. The solution is also unique. Indeed, let us consider two solutions $(\mathbf{v}_1,\mathbf{F}_1,\mathbf{M}_1,\phi_1,\mu_1,\beta_1)$ and $(\mathbf{v}_2,\mathbf{F}_2,\mathbf{M}_2,\phi_2,\mu_2,\beta_2)$ of System \eqref{eqn:galerkinfe2div}, satisfying the same initial and boundary conditions, and let us define $\bar{\mathbf{v}}=\mathbf{v}_1-\mathbf{v}_2$, $\bar{\mathbf{F}}=\mathbf{F}_1-\mathbf{F}_2$, $\bar{\mathbf{M}}=\mathbf{M}_1-\mathbf{M}_2$, $\bar{\phi}=\phi_1-\phi_2$, $\bar{\mu}=\mu_1-\mu_2$, $\bar{\beta}=\beta_1-\beta_2$. Taking the difference of the equations in \eqref{eqn:galerkinfe2div} for the variables $(\mathbf{v}_1,\mathbf{F}_1,\mathbf{M}_1,\phi_1,\mu_1,\beta_1)$ and $(\mathbf{v}_2,\mathbf{F}_2,\mathbf{M}_2,\phi_2,\mu_2,\beta_2)$, we obtain
\begin{equation}
\label{eqn:galerkinfe2div2}
\begin{cases}
\displaystyle \nu \int_{\Omega}\nabla \bar{\mathbf{v}}\colon \nabla \mathbf{u}_h =\frac{\bar{\beta}}{\sqrt{\int_{\Omega}(j(\phi_h^n,\mathbf{F}_h^n))+k}}\int_{\Omega}\mu_h^{n} \nabla \phi_h^n \cdot \mathbf{u}_h\\
\displaystyle -\frac{\bar{\beta}}{\sqrt{\int_{\Omega}(j(\phi_h^n,\mathbf{F}_h^n))+k}}\int_{\Omega}\left(\mathbf{M}_h^{n}(\mathbf{F}_h^n)^T\right)\colon \nabla\mathbf{u}_h+\frac{\bar{\beta}}{\sqrt{\int_{\Omega}(j(\phi_h^n,\mathbf{F}_h^n))+k}}\int_{\Omega}\left(\nabla \mathbf{F}_h^n\odot \mathbf{M}_h^{n}\right)\cdot \mathbf{u}_h,\\ \\
 \displaystyle \int_{\Omega} \bar{\mathbf{F}} \colon \boldsymbol{\Theta}_h+\Delta t \frac{\bar{\beta}}{\sqrt{\int_{\Omega}(j(\phi_h^n,\mathbf{F}_h^n))+k}}\int_{\Omega}\left({\mathbf{v}_h^{n}} \cdot \nabla \right)\mathbf{F}_h^n\colon \boldsymbol{\Theta}_h\\
 \displaystyle -\Delta t\frac{\bar{\beta}}{\sqrt{\int_{\Omega}(j(\phi_h^n,\mathbf{F}_h^n))+k}}\int_{\Omega}\left(\nabla {\mathbf{v}_h^{n}}\right)\mathbf{F}_h^n\colon \boldsymbol{\Theta}_h+\Delta t\gamma \int_{\Omega}\bar{\mathbf{M}}\colon \boldsymbol{\Theta}_h=0,\\ \\
 \displaystyle \int_{\Omega}\bar{\mathbf{M}}\colon \boldsymbol{\Pi}_h=\frac{\bar{\beta}}{\sqrt{\int_{\Omega}(j(\phi_h^n,\mathbf{F}_h^n))+k}}\int_{\Omega}\partial_{\mathbf{F}}w(\phi_h^n,\mathbf{F}_h^n)\colon \boldsymbol{\Pi}_h+\lambda \int_{\Omega}\nabla \bar{\mathbf{F}} \mathbin{\tensorm} \nabla \boldsymbol{\Pi}_h,\\ \\
 \displaystyle \int_{\Omega}\bar{\phi}\,\xi_h+\frac{\bar{\beta}}{\sqrt{\int_{\Omega}(j(\phi_h^n,\mathbf{F}_h^n))+k}}\Delta t\int_{\Omega}\left(\mathbf{v}_h^{n} \cdot \nabla \phi_h^n \right) \xi_h+\Delta t\int_{\Omega}b\left(\phi_h^n\right)\nabla \bar{\mu} \cdot \nabla \xi_h =0,\\ \\
\displaystyle \int_{\Omega}\bar{\mu} \chi_h=\int_{\Omega}\nabla \bar{\phi} \cdot \nabla \chi_h+\frac{\bar{\beta}}{\sqrt{\int_{\Omega}(j(\phi_h^n,\mathbf{F}_h^n))+k}}\int_{\Omega}\partial_{\phi}j(\phi_h^n,\mathbf{F}_h^n)\chi_h,\\ \\
\displaystyle \bar{\beta}=\frac{1}{2\sqrt{\int_{\Omega}(j(\phi_h^n,\mathbf{F}_h^n))+k}}\biggl[\left(\int_{\Omega}\partial_{\phi}j(\phi_h^n,\mathbf{F}_h^n)\bar{\phi}\right)\\
\displaystyle+\int_{\Omega}\partial_{\mathbf{F}}w(\phi_h^n,\mathbf{F}_h^n)\bar{\mathbf{F}}+\Delta t\int_{\Omega}\left({\mathbf{v}_h^{n}} \cdot \nabla \right)\mathbf{F}_h^n\colon \bar{\mathbf{M}}-\Delta t\int_{\Omega}\left(\nabla \mathbf{F}_h^n\right)^T\mathbf{M}_h^{n}\cdot \bar{\mathbf{v}}\\
\displaystyle -\Delta t \int_{\Omega}\left(\nabla {\mathbf{v}_h^{n}}\right)\mathbf{F}_h^n\colon \bar{\mathbf{M}}+\Delta t\int_{\Omega}\left(\mathbf{M}_h^{n}(\mathbf{F}_h^n)^T\right)\colon \nabla\bar{\mathbf{v}}+\Delta t\int_{\Omega}\left(\mathbf{v}_h^{n} \cdot \nabla \phi_h^n \right) \bar{\mu}\\
\displaystyle -\Delta t \int_{\Omega}\mu_h^{n} \nabla \phi_h^n \cdot \bar{\mathbf{v}}\biggr].
\end{cases}
\end{equation}
Taking $\mathbf{u}_h=\Delta t \bar{\mathbf{v}}$, $\boldsymbol{\Theta}_h=\bar{\mathbf{M}}$, $\boldsymbol{\Pi}_h=-\bar{\mathbf{F}}$, $\xi_h=\bar{\mu}$, $\chi_h=-\bar{\phi}$ in \eqref{eqn:galerkinfe2div2} and multiplying \eqref{eqn:galerkinfe2div2}$_6$ by $2\bar{\beta}$, summing all the equations and using Assumption \textbf{A1}, we obtain that
\begin{align}
\label{eqn:galerkinfe2div3}  
& \displaystyle \Delta t\nu ||\nabla \bar{\mathbf{v}}||_{L^2(\Omega;\mathbb{R}^{3\times 3})}^2+\Delta t\gamma || \bar{\mathbf{M}}||_{L^2(\Omega;\mathbb{R}^{3\times 3})}^2+\lambda||\nabla \bar{\mathbf{F}}||_{L^2(\Omega;\mathbb{R}^{3\times 3 \times 3})}^2\\
& \displaystyle \notag+\Delta t b_0||\nabla \bar{\mu}||_{L^2(\Omega;\mathbb{R}^{3})}^2+||\nabla \bar{\phi}||_{L^2(\Omega;\mathbb{R}^{3})}^2+2\bar{\beta}^2\leq 0.
\end{align}
Taking moreover $\boldsymbol{\Theta}_h=\mathbf{e}_{l,r}$, $l,r=1,2,3$ and $\xi_h=\chi_h=1$
in \eqref{eqn:galerkinfe2div2}, using moreover the fact that $\mathbf{v}_{h}^{n}\in V_{h,\text{div}}$, we obtain that
\begin{equation}
\label{eqn:galerkinfe2div4}  
\left|\int_{\Omega}\bar{\mathbf{F}}\right|\leq C||\bar{\mathbf{M}}||_{L^2(\Omega;\mathbb{R}^{3\times 3})}, \; \int_{\Omega}\bar{\phi}=0, \; \left|\int_{\Omega}\bar{\mu}\right|\leq C|\bar{\beta}|.
\end{equation}
From \eqref{eqn:galerkinfe2div3} and \eqref{eqn:galerkinfe2div4} we finally conclude that $\bar{\mathbf{v}}=\mathbf{0}$, $\bar{\mathbf{F}}=\mathbf{0}$, $\bar{\mathbf{M}}=\mathbf{0}$, $\bar{\phi}=0$, $\bar{\mu}=0$, $\bar{\beta}=0$, hence the solution of System \eqref{eqn:galerkinfe2div} is unique. 
Since \eqref{eqn:galerkinfe2div}$_1$ defines a linear functional from $V_h$ to $\mathbb{R}$ which vanishes on $V_{h,\text{div}}$, we conclude by standard arguments (see e.g. \cite[p.22]{gr}) that there exists a unique $s_h^{n+1}\in R_h$ such that 
 \[
 (\mathbf{v}_{h}^{n+1},s_h^{n+1},\mathbf{F}_{h}^{n+1},\mathbf{M}_{h}^{n+1},\phi_{h}^{n+1},\mu_{h}^{n+1},\beta_h^{n+1})\in V_h\times R_h\times X_h \times X_h \times S_h \times X_h \times \mathbb{R}
 \]
 is the unique solution of System \eqref{eqn:galerkinfe2}.
 
Let us now take $\mathbf{u}_h=\Delta t \mathbf{v}_h^{n+1}$, $p_h=\Delta t s_h^{n+1}$, $\boldsymbol{\Theta}_h=\mathbf{M}_h^{n+1}$, $\boldsymbol{\Pi}_h=-(\mathbf{F}_h^{n+1}-\mathbf{F}_h^{n})$, $\xi_h=\mu_h^{n+1}$, $\chi_h=-(\phi_h^{n+1}-\phi_h^n)$ in \eqref{eqn:galerkinfe2} and let us multiply \eqref{eqn:galerkinfe2}$_7$ by $2\beta_h^{n+1}$. Summing all the equations, and using the identity $a\star (a-b)=\frac{1}{2}|a|^2+\frac{1}{2}|a-b|^2-\frac{1}{2}|b|^2$, we obtain
\begin{align}
\label{eqn:discreteenergy2}
&\displaystyle\Delta t\nu ||\nabla \mathbf{v}_h^{n+1}||_{L^2(\Omega;\mathbb{R}^{3\times 3})}^2+\Delta t\gamma || \mathbf{M}_h^{n+1}||_{L^2(\Omega;\mathbb{R}^{3\times 3})}^2+\frac{\lambda}{2}||\nabla \mathbf{F}_h^{n+1}||_{L^2(\Omega;\mathbb{R}^{3\times 3 \times 3})}^2\\
&\notag \displaystyle+\frac{\lambda}{2}\left|\left|\nabla \left(\mathbf{F}_h^{n+1}-\mathbf{F}_h^{n}\right)\right|\right|_{L^2(\Omega;\mathbb{R}^{3\times 3 \times 3})}^2+
\Delta t \int_{\Omega}b\left(\phi_h^n\right)\nabla \mu_h^{n+1}\cdot \nabla \mu_h^{n+1}\\
&\notag +\frac{1}{2}||\nabla \phi_h^{n+1}||_{L^2(\Omega;\mathbb{R}^{3})}^2+\frac{1}{2}\left|\left|\nabla \left(\phi_h^{n+1}-\phi_h^{n}\right)\right|\right|_{L^2(\Omega;\mathbb{R}^{3})}^2+ \left(\beta_h^{n+1}\right)^2+\left(\beta_h^{n+1}-\beta_h^{n}\right)^2\\
 &\notag \displaystyle =\frac{\lambda}{2}||\nabla \mathbf{F}_h^{n}||_{L^2(\Omega;\mathbb{R}^{3\times 3 \times 3})}^2+\frac{1}{2}||\nabla \phi_h^{n}||_{L^2(\Omega;\mathbb{R}^{3})}^2+ \left(\beta_h^{n}\right)^2,
\end{align}
which gives that \eqref{eqn:lyapunov2} is a decreasing (Lyapunov) function for the discrete solutions. Summing \eqref{eqn:discreteenergy2} from $n=0\to m$, for $m=0\to N-1$, we finally get \eqref{eqn:stabbound2}.
\end{pf}
\\
\noindent
Following the ideas reported \cite{lishen}, starting from the coupled system \eqref{eqn:galerkinfe2} we obtain a decoupled system by introducing the following ansatz for its solution:
\begin{align}
\label{eqn:ansatzdecoupled}
& \mathbf{v}_h^{n+1}=\mathbf{v}_{h,0}^{n+1}+A_h^{n+1}\mathbf{v}_{h,1}^{n+1},\\
& \notag s_h^{n+1}=s_{h,0}^{n+1}+A_h^{n+1}s_{h,1}^{n+1},\\
& \notag \mathbf{F}_h^{n+1}=\mathbf{F}_{h,0}^{n+1}+A_h^{n+1}\mathbf{F}_{h,1}^{n+1},\\
& \notag \mathbf{M}_h^{n+1}=\mathbf{M}_{h,0}^{n+1}+A_h^{n+1}\mathbf{M}_{h,1}^{n+1},\\
& \notag \phi_h^{n+1}=\phi_{h,0}^{n+1}+A_h^{n+1}\phi_{h,1}^{n+1},\\
& \notag \mu_h^{n+1}=\mu_{h,0}^{n+1}+A_h^{n+1}\mu_{h,1}^{n+1},\\
\end{align}
where 
\[A_h^{n+1}:=\frac{\beta_h^{n+1}}{\sqrt{\int_{\Omega}(j(\phi_h^n,\mathbf{F}_h^n))+k}}.\] 
Inserting the expansions \eqref{eqn:ansatzdecoupled} in \eqref{eqn:galerkinfe2} and equating the terms of the same order in $A_h^{n+1}$, we obtain that the variables $\mathbf{v}_{h,i}^{n+1},s_{h,i}^{n+1},\mathbf{F}_{h,i}^{n+1},\mathbf{F}_{h,i}^{n+1},\phi_{h,i}^{n+1},\mu_{h,i}^{n+1}$, $i=0,1$, satisfy the following decoupled systems.

\noindent
\textbf{Stokes subsystem:}
\begin{equation}
\label{eqn:subS}
\begin{cases}
\displaystyle \nu \int_{\Omega}\nabla \mathbf{v}_{h,0}^{n+1}\colon \nabla \mathbf{u}_h-\int_{\Omega}s_{h,0}^{n+1}\diver \mathbf{u}_h =0,\\
\displaystyle \int_{\Omega}p_h\diver \mathbf{v}_{h,0}^{n+1}=0\\ \\
\displaystyle \nu \int_{\Omega}\nabla \mathbf{v}_{h,1}^{n+1}\colon \nabla \mathbf{u}_h-\int_{\Omega}s_h^{n+1}\diver \mathbf{u}_h =\int_{\Omega}\mu_h^{n} \nabla \phi_h^n \cdot \mathbf{u}_h-\int_{\Omega}\left(\mathbf{M}_h^{n}(\mathbf{F}_h^n)^T\right)\colon \nabla\mathbf{u}_h\\
\displaystyle +\int_{\Omega}\left(\nabla \mathbf{F}_h^n\odot \mathbf{M}_h^{n}\right)\cdot \mathbf{u}_h,\\
\displaystyle \int_{\Omega}p_h\diver \mathbf{v}_{h,1}^{n+1}=0.
\end{cases}
\end{equation}
Since $\mathbf{v}_{h,0}^{n+1}|_{\partial \Omega}=0$, \eqref{eqn:subS}$_{1,2}$ imply that $\mathbf{v}_{h,0}^{n+1}\equiv 0$.

\noindent
\textbf{Elasticity subsystem:}
\begin{equation}
\label{eqn:subEL}
\begin{cases}
\displaystyle \int_{\Omega} \left(\mathbf{F}_{h,0}^{n+1}-\mathbf{F}_h^{n}\right) \colon \boldsymbol{\Theta}_h+\Delta t\gamma \int_{\Omega}\mathbf{M}_{h,0}^{n+1}\colon \boldsymbol{\Theta}_h=0,\\
 \displaystyle \int_{\Omega}\mathbf{M}_{h,0}^{n+1}\colon \boldsymbol{\Pi}_h=\lambda \int_{\Omega}\nabla \mathbf{F}_{h,0}^{n+1} \mathbin{\tensorm} \nabla \boldsymbol{\Pi}_h,\\ \\
 \displaystyle \int_{\Omega} \mathbf{F}_{h,1}^{n+1} \colon \boldsymbol{\Theta}_h+\Delta t \int_{\Omega}\left({\mathbf{v}_h^{n}} \cdot \nabla \right)\mathbf{F}_h^n\colon \boldsymbol{\Theta}_h-\Delta t\int_{\Omega}\left(\nabla {\mathbf{v}_h^{n}}\right)\mathbf{F}_h^n\colon \boldsymbol{\Theta}_h+\Delta t\gamma \int_{\Omega}\mathbf{M}_{h,1}^{n+1}\colon \boldsymbol{\Theta}_h=0,\\
 \displaystyle \int_{\Omega}\mathbf{M}_{h,1}^{n+1}\colon \boldsymbol{\Pi}_h=\int_{\Omega}\partial_{\mathbf{F}}w(\phi_h^n,\mathbf{F}_h^n)\colon \boldsymbol{\Pi}_h+\lambda \int_{\Omega}\nabla \mathbf{F}_{h,1}^{n+1} \mathbin{\tensorm} \nabla \boldsymbol{\Pi}_h.
\end{cases}
\end{equation}

\noindent
\textbf{Cahn--Hilliard subsystem:}
\begin{equation}
\label{eqn:subCH}
\begin{cases}
 \displaystyle \int_{\Omega}\left(\phi_{h,0}^{n+1}-\phi_h^n\right)\xi_h+\Delta t\int_{\Omega}b\left(\phi_h^n\right)\nabla \mu_{h,0}^{n+1} \cdot \nabla \xi_h =0,\\
\displaystyle \int_{\Omega}\mu_{h,0}^{n+1} \chi_h=\int_{\Omega}\nabla \phi_{h,0}^{n+1} \cdot \nabla \chi_h,\\ \\
 \displaystyle \int_{\Omega}\phi_{h,1}^{n+1}\xi_h+\Delta t\int_{\Omega}\left(\mathbf{v}_h^{n} \cdot \nabla \phi_h^n \right) \xi_h+\Delta t\int_{\Omega}b\left(\phi_{h}^n\right)\nabla \mu_h^{n+1} \cdot \nabla \xi_h =0,\\
\displaystyle \int_{\Omega}\mu_{h,1}^{n+1} \chi_h=\int_{\Omega}\nabla \phi_{h,1}^{n+1} \cdot \nabla \chi_h+\int_{\Omega}\partial_{\phi}j(\phi_h^n,\mathbf{F}_h^n)\chi_h.
\end{cases}
\end{equation}
Moreover,
\begin{align}
\label{eqn:an}
& \displaystyle \biggl[\frac{\sqrt{\int_{\Omega}(j(\phi_h^n,\mathbf{F}_h^n))+k}}{\Delta t}-\frac{1}{2\sqrt{\int_{\Omega}(j(\phi_h^n,\mathbf{F}_h^n))+k}}\biggl\{\left(\int_{\Omega}\partial_{\phi}j(\phi_h^n,\mathbf{F}_h^n)\frac{\phi_{h,1}^{n+1}}{\Delta t}\right)\\
& \displaystyle \notag + \left(\int_{\Omega}\partial_{\mathbf{F}}w(\phi_h^n,\mathbf{F}_h^n)\colon \frac{\mathbf{F}_{h,1}^{n+1}}{\Delta t}\right)+\int_{\Omega}\left({\mathbf{v}_h^{n}} \cdot \nabla \right)\mathbf{F}_h^n\colon \mathbf{M}_{h,1}^{n+1}-\int_{\Omega}\left(\nabla \mathbf{F}_h^n\right)^T\mathbf{M}_h^{n}\cdot \mathbf{v}_{h,1}^{n+1}\\
& \displaystyle -\int_{\Omega}\left(\nabla {\mathbf{v}_h^{n}}\right)\mathbf{F}_h^n\colon \mathbf{M}_{h_1}^{n+1}+\int_{\Omega}\left(\mathbf{M}_h^{n}(\mathbf{F}_h^n)^T\right)\colon \nabla\mathbf{v}_{h,1}^{n+1}+\int_{\Omega}\left(\mathbf{v}_h^{n} \cdot \nabla \phi_h^n \right) \mu_{h,1}^{n+1}\\
& \displaystyle -\int_{\Omega}\mu_h^{n} \nabla \phi_h^n \cdot \mathbf{v}_{h,1}^{n+1}\biggr\}\biggr]A_h^{n+1}\\
& \displaystyle =\biggl[\frac{\beta_h^n}{\Delta t}+\frac{1}{2\sqrt{\int_{\Omega}(j(\phi_h^n,\mathbf{F}_h^n))+k}}\biggl\{\left(\int_{\Omega}\partial_{\phi}j(\phi_h^n,\mathbf{F}_h^n)\frac{\phi_{h,0}^{n}-\phi_h^n}{\Delta t}\right)\\
& \displaystyle \notag + \left(\int_{\Omega}\partial_{\mathbf{F}}w(\phi_h^n,\mathbf{F}_h^n)\colon \frac{\mathbf{F}_{h,0}^{n+1}-\mathbf{F}_{h}^{n}}{\Delta t}\right)+\int_{\Omega}\left({\mathbf{v}_h^{n}} \cdot \nabla \right)\mathbf{F}_h^n\colon \mathbf{M}_{h,0}^{n+1}-\int_{\Omega}\left(\nabla {\mathbf{v}_h^{n}}\right)\mathbf{F}_h^n\colon \mathbf{M}_{h_0}^{n+1}\\
& \displaystyle +\int_{\Omega}\left(\mathbf{v}_h^{n} \cdot \nabla \phi_h^n \right) \mu_{h,0}^{n+1}\biggr\}\biggr].
\end{align}
It's easy to show, by taking $\mathbf{u}_h=\mathbf{v}_{h,1}^{n+1}$ and $p_h=s_{h,1}^{n+1}$ in \eqref{eqn:subS}$_{3,4}$, $\boldsymbol{\Theta}_h=\mathbf{M}_{h,1}^{n+1}$, $\boldsymbol{\Pi}_h=\frac{\mathbf{F}_{h,1}^{n+1}}{\Delta t}$ in \eqref{eqn:subEL}$_{3,4}$ and $\xi_h=\mu_{h,1}^{n+1}$, $\chi_h=\frac{\phi_{h,1}^{n+1}}{\Delta t}$ in \eqref{eqn:subCH}$_{3,4}$, that the term in the square brackets which multiplies $A_h^{n+1}$ in \eqref{eqn:an} is strictly positive, hence $A_h^{n+1}$ is uniquely determined. 

\section{Results}
In this section we report the results of numerical simulations in two space dimensions for different test cases. 
We consider the following form of the free energy density:
\begin{equation}
\label{eqn:psinum}
\psi(\phi)=\frac{\beta}{4\alpha}\phi^2(1-\phi)^2,
\end{equation}
\begin{equation}
\label{eqn:wnum}
w(\phi,\mathbf{F})=\frac{\zeta}{2}\left(\mathbf{F}^T\mathbf{F}-\mathbf{H}(\phi)\right)\colon \left(\mathbf{F}^T\mathbf{F}-\mathbf{H}(\phi)\right),
\end{equation}
where $\beta,\alpha,\zeta>0$ are model parameters, and
\[
 \mathbf{H}(\phi)=\tilde{\mathbf{F}}(\phi)^T\tilde{\mathbf{F}}(\phi),
 \]
 with
 \[
 \tilde{\mathbf{F}}(\phi)=\begin{pmatrix}
 1 && a\phi\\
 0 && 1
\end{pmatrix}
\]
 and $a>0$. In particular, \eqref{eqn:psinum} is the standard Cahn--Hilliard smooth double-well potential, with stable minima at $\phi\equiv 0$ and $\phi\equiv 1$, where the parameter $\beta$ is proportional to the surface tension and the parameter $\alpha$ represents the interface thickness in the Cahn--Hilliard surface term. Consequently, the term proportional to $|\nabla \phi|^2$ in \eqref{eqn:free} is multiplied by a factor $\epsilon=\beta \alpha$, ($\epsilon, \alpha, \beta$ were taken equal to one in the previous analysis for ease of notation). Moreover, \eqref{eqn:wnum} represents an elastic energy density of shape memory alloy type (see e.g. \cite[Eq. (3.20a)]{rubiceksm}), where the pure phases associated to the stable minima of the phase field potential are characterized by different elastic properties. Indeed, we have that
 \[
 \partial_{\phi}\left(\psi(\phi)+w(\phi,\mathbf{F})\right)=\psi^{\prime}(\phi)-\zeta\mathbf{H}^{\prime}(\phi)\colon\left(\mathbf{F}^T\mathbf{F}-\mathbf{H}(\phi)\right),
 \]
 and
 \[
 \partial_{\mathbf{F}}w(\phi,\mathbf{F})=2\zeta\mathbf{F}\left(\mathbf{F}^T\mathbf{F}-\mathbf{H}(\phi)\right).
 \]
 Hence, the global minima for $\psi(\phi)+w(\phi,\mathbf{F})$, at which it takes its minimum value equal to zero, are attained at $\phi\equiv 0$, $\phi\equiv 1$ and $\mathbf{F}\equiv \mathbf{R}\tilde{\mathbf{F}}(\phi)$, where $\mathbf{R}\in SO(\mathbb{R}^{2\times 2})$ is a generic rotation matrix. We observe that \eqref{eqn:wnum} is frame indifferent but not isotropic.
 
 In the following, we will consider two test cases to investigate the phase separation dynamics and the coarsening dynamics associated to \eqref{eqn:psinum} and \eqref{eqn:wnum}, starting from different initial configurations. In \textit{Test Case 1} we will consider the phase separation dynamics starting from different initial conditions for $\phi$. In particular, we will consider for $\phi_0$ two different configurations given by a small uniformly distributed random perturbation around the values $\phi_0=0.3$ or $\phi_0=0.7$, leading to different topologies of the stationary states, determined both from the metastability of $\psi(\cdot)$ and from the elasticity dynamics described by $w(\phi,\mathbf{F})$. In both cases, we will take as initial condition for the deformation gradient $\mathbf{F}_0=\mathbf{I}$.

The values of the parameters for \textit{Test Case 1} are taken as $\nu=1$, $\lambda=0.001$, $\beta=0.1$, $\alpha=0.002$, $\zeta=10$, $a=0.5$. Moreover, we take a constant mobility $b(\phi)\equiv 1$ and we consider a domain $\Omega=[0,1]\times [0,1]$, with a uniform triangulation of dimension $64\times 64$, and $\Delta t=0.001 \epsilon$. In \textit{Test Case 1} we will vary the value of the parameter $\gamma$, considering $\gamma=1$ or $\gamma=0.001$, in order to observe the phase separation dynamics under different degrees of diffusive regularization of the transport equation for the deformation gradient. 

In \textit{Test Case 2} we will consider the merging and coarsening dynamics of isolated circular subregions of a pure phase immersed in a bath of the opposite pure phase. In particular, we will consider two different configurations with two and four initial circular subregions placed symmetrically with respect to the centre of the domain. The values of the parameters for \textit{Test Case 2} are taken as $\nu=1$, $\lambda=0.001$, $\beta=0.1$, $\alpha=0.02$, $a=0.5$ and $\gamma=0.001$, $b(\phi)\equiv 1$. We will vary the value of the elastic modulus, taking $\zeta=1$ and $\zeta=10$, in order to observe the effects of higher stiffness on the merging and coarsening dynamics. We consider a domain $\Omega=[0,1]\times [0,1]$, with a uniform triangulation of dimension $64\times 64$, furtherly refined in a neighborhood of the support set of $\phi_0$. The time step is taken as $\Delta t=0.00001 \epsilon$. 
 \begin{rem}
 \label{rem:gradstabnum}
 The expression \eqref{eqn:wnum} can be rewritten as
 \begin{align}
\label{eqn:wnum2}
w(\phi,\mathbf{F})=\frac{\zeta}{2}\mathbf{F}^T\mathbf{F}\colon \mathbf{F}^T\mathbf{F}-\zeta \mathbf{H}(\phi)\colon \mathbf{F}^T\mathbf{F}+\frac{\zeta}{2}\mathbf{H}(\phi)\colon \mathbf{H}(\phi),
\end{align}
which is of the particular form \eqref{eqn:wassum2} and satisfies Assumption $\bf{A2_{h}}$ in the general case of Remark \ref{rem:a2hgen} when $\phi$ is bounded, which is always practically verified by the numerical solutions.
 \end{rem}
The fully coupled system \eqref{eqn:galerkinfe1} is solved through the iteration method \eqref{eqn:galerkinfe1dec}, while the solution of \eqref{eqn:galerkinfe2} is obtained by solving the independent subsystems \eqref{eqn:subS}, \eqref{eqn:subEL}, \eqref{eqn:subCH} and \eqref{eqn:an} and using \eqref{eqn:ansatzdecoupled}.
\subsection{Gradient stability of \eqref{eqn:galerkinfe1} and \eqref{eqn:galerkinfe2}}
In this section we show the gradient stability of \eqref{eqn:galerkinfe1} and \eqref{eqn:galerkinfe2}. As representative results, we report the solutions of \eqref{eqn:galerkinfe1} and \eqref{eqn:galerkinfe2}, with the energy density defined in \eqref{eqn:psinum} and \eqref{eqn:wnum}, in the case $\gamma=1$, $\phi_0=0.3\pm 0.5\iota$ and $\mathbf{F}_0=\mathbf{I}$, where $\iota$ is a random perturbation uniformly distributed in the interval $[0,1]$. 

In Figure \ref{fig:1} we report the plots of the Lyapunov functionals \eqref{eqn:lyapunov1} and  \eqref{eqn:lyapunov2} (subtracting to it the value of $k$), which we call $L_{CS}$ and $L_{DSAV}$ respectively, versus time, together with the functionals
\[
E_{CH}:=\int_{\Omega}\left(\psi(\phi)+\frac{\epsilon}{2}|\nabla \phi|^2\right)
\]
and 
\[
E_{EL}:=\int_{\Omega}\left(w(\phi,\mathbf{F})+\frac{\lambda}{2}|\nabla \mathbf{F}|^2\right).
\]

\begin{figure}[h!]
\includegraphics[width=0.9\linewidth]
{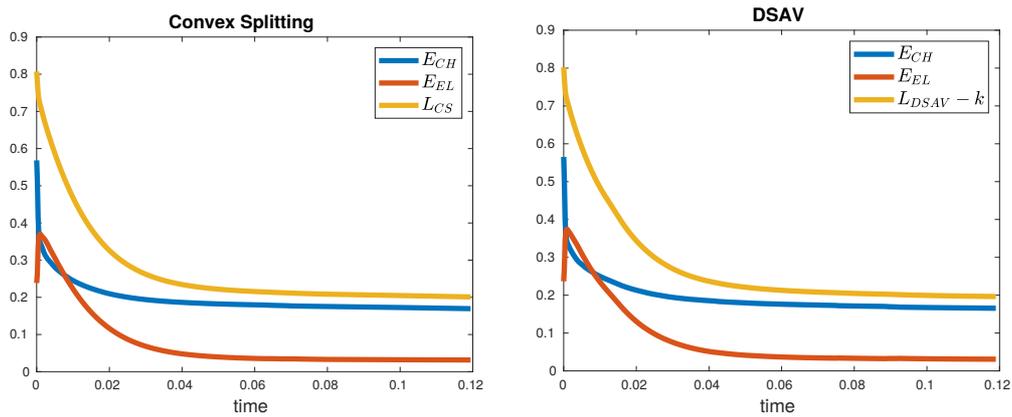}
\centering
\caption{Plot of the Lyapunov functionals \eqref{eqn:lyapunov1} ($L_{CS}$) and \eqref{eqn:lyapunov2} ($L_{DSAV}$) and of the functionals $E_{CH}$ and $E_{EL}$ vs time for the schemes \eqref{eqn:galerkinfe1} (Convex Splitting, left panels) and \eqref{eqn:galerkinfe2} (DSAV, right panels).}
\label{fig:1}
\end{figure}
We observe from Figure \ref{fig:1} that both systems \eqref{eqn:galerkinfe1} and \eqref{eqn:galerkinfe2} are gradient stable. The plots show that the two systems behave almost identically over the whole time span until the attainment of the stationary state, with the convex splitting scheme showing a steeper decrease of the Lyapunov functional at early times than the DSAV scheme. 
In Figure \ref{fig:2} we compare the plots of $\phi$ and $|\mathbf{F}|$ for the two schemes at initial and at late time points, observing that the numerical solutions for $\phi$ and $\mathbf{F}$ of \eqref{eqn:galerkinfe1} and \eqref{eqn:galerkinfe2} are statistically and topologically similar throughout the whole dynamics. Note that, since the initial condition $\phi_0$ is random, the numerical simulations for the two schemes do not start from the same initial conditions. In Figure \ref{fig:3} we also compare the line plots, along a vertical line, of $|\mathbf{F}|$ at time $t=0.01$, i.e. at an early time where the two schemes show slight differences in the Lyapunov functional decrease, observing higher numerical oscillations for the scheme \eqref{eqn:galerkinfe2}.

\begin{figure}[ht!]
\includegraphics[width=0.9\linewidth]
{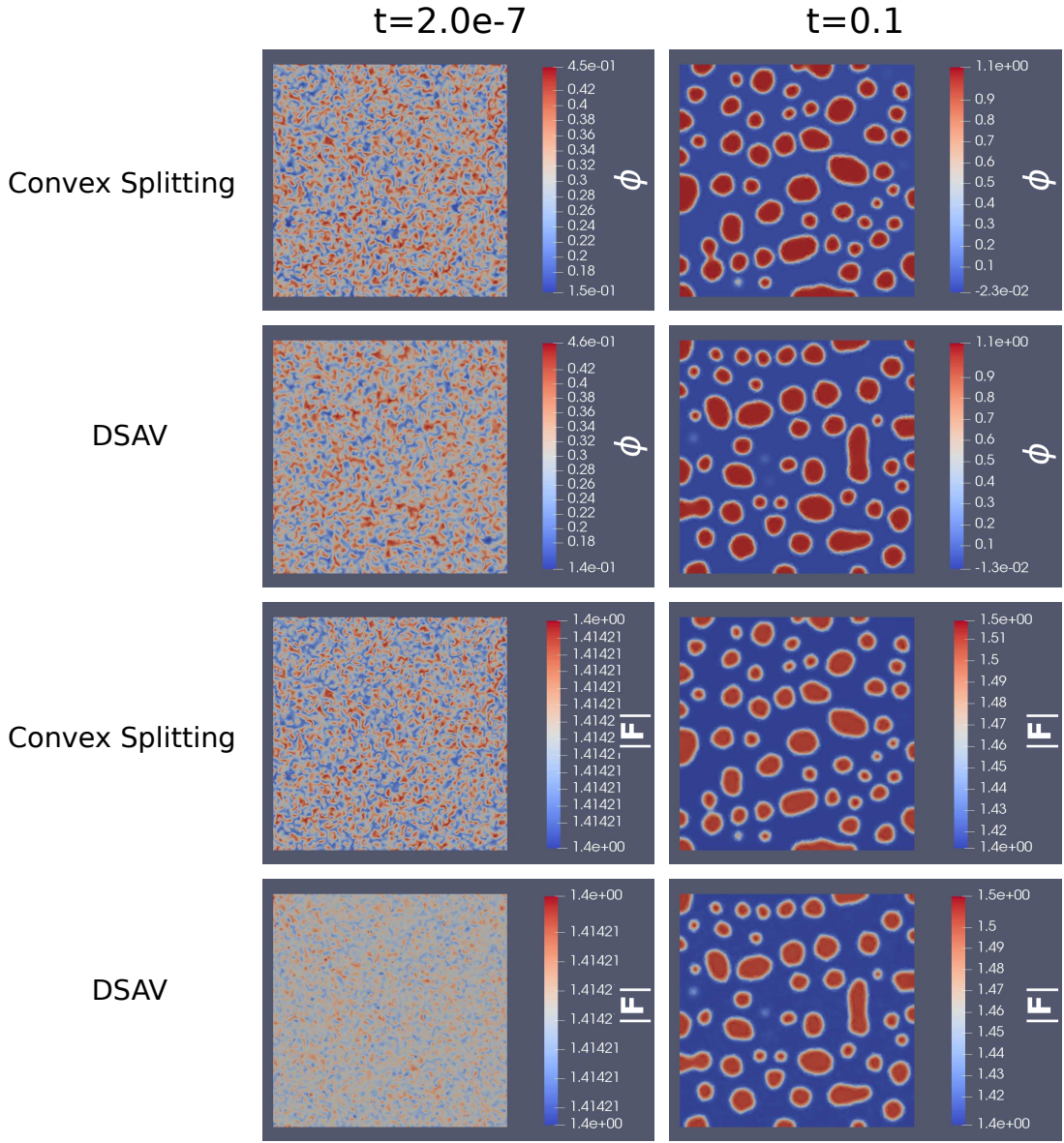}
\centering
\caption{Plots of $\phi$ and $|\mathbf{F}|$ at different time points for the schemes \eqref{eqn:galerkinfe1} (Convex Splitting) and \eqref{eqn:galerkinfe2} (DSAV).}
\label{fig:2}
\end{figure}

\begin{figure}[ht!]
\includegraphics[width=0.9\linewidth]
{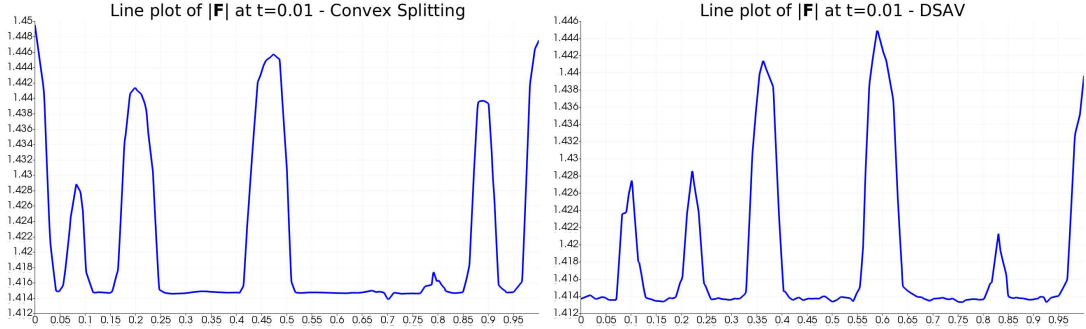}
\centering
\caption{Line plots of $|\mathbf{F}|$, along a vertical line, at time $t=0.01$, for the schemes \eqref{eqn:galerkinfe1} (Convex Splitting) and \eqref{eqn:galerkinfe2} (DSAV).}
\label{fig:3}
\end{figure}
The computational time to solve \eqref{eqn:galerkinfe1} for $150000$ time steps is $\sim 1.5\mathrm{e}{+6}\, s$, while computational time to solve \eqref{eqn:galerkinfe2} for $150000$ time steps is $\sim 3.4\mathrm{e}{+5}\, s$. Hence the time needed to solve \eqref{eqn:galerkinfe1} is almost one order of magnitude greater than the computational time needed to solve \eqref{eqn:galerkinfe2}. 

In the following, the numerical results for \textit{Test Case 1} and \textit{Test Case 2} are obtained as solutions of \eqref{eqn:galerkinfe1}.
\subsection{Test case 1 -- Phase separation}
We first consider the initial conditions $\phi_0=0.3\pm 0.5\iota$, where $\iota$ is a random perturbation uniformly distributed in the interval $[0,1]$, and $\mathbf{F}_0=\mathbf{I}$. In Figure \ref{fig:4} we show the numerical results at different time points throughout the phase separation dynamics, up to late times at which we can observe the coarsening dynamics of the separated domain subregions, for both the cases with $\gamma=1$ and $\gamma=0.001$.

\begin{figure}[ht!]
\includegraphics[width=0.9\linewidth]
{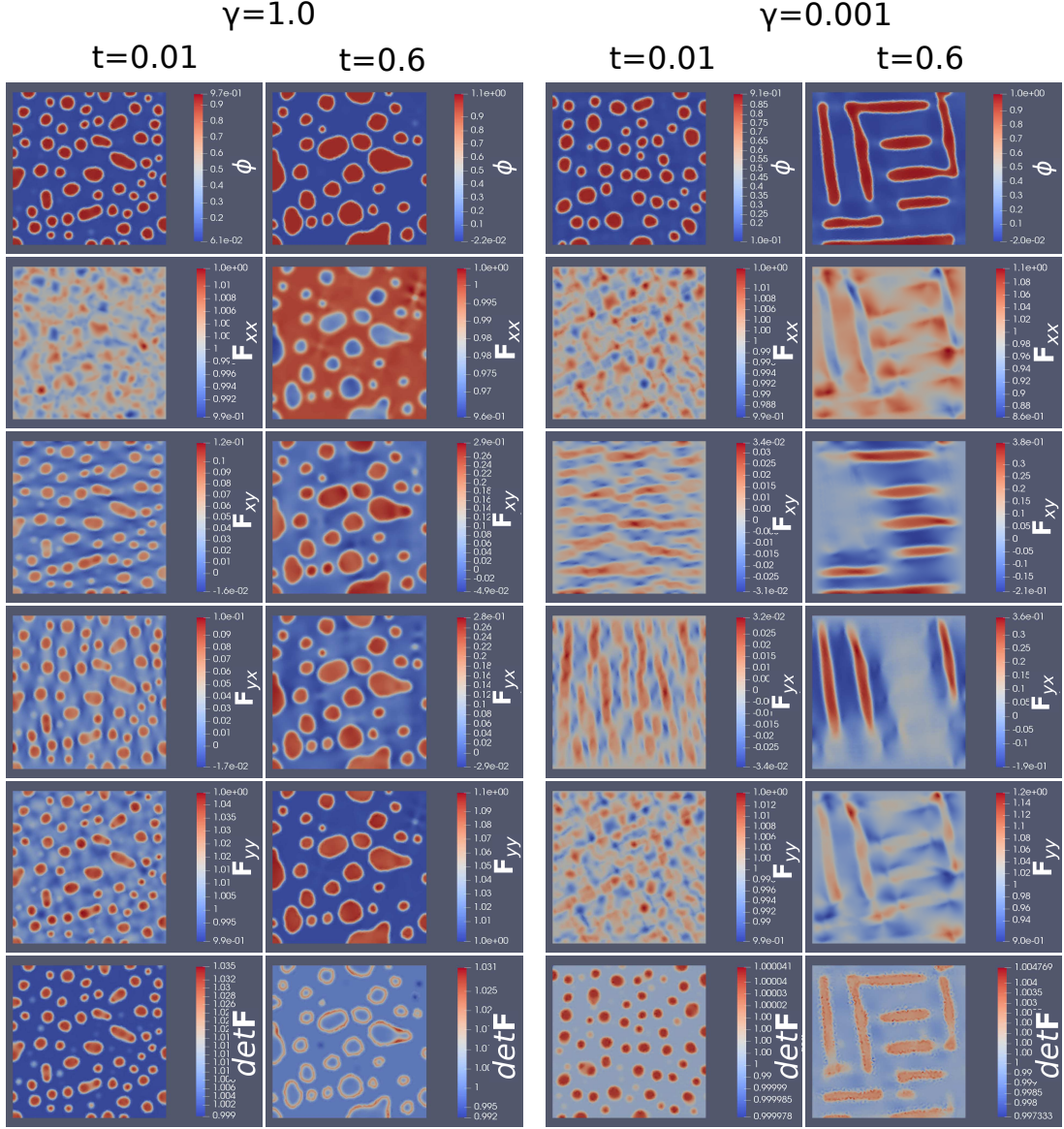}
\centering
\caption{Plot of $\phi$ (first row) and $\mathbf{F}_{xx}$ (second row), $\mathbf{F}_{xy}$ (third row), $\mathbf{F}_{yx}$ (fourth row), $\mathbf{F}_{yy}$ (fifth row), $\text{det}\mathbf{F}$ (sixth row) at different time points in the case $\phi_0=0.3\pm 0.5\iota$, for $\gamma=1$ (first and second columns) and $\gamma=0.001$ (third and fourth columns).}
\label{fig:4}
\end{figure}
We observe from Figure \ref{fig:4} that the phase separation dynamics for the $\phi$ variable consists in the formation of elongated circular clusters with $\phi\sim 1$ along orthogonal directions determined by the separation dynamics for the $\mathbf{F}$ variable, immersed in a bath with $\phi\sim 0$. At late times, these clusters grow and merge, forming strips oriented in orthogonal directions. The variable $\mathbf{F}$ assumes at late times values $\mathbf{F}\sim \mathbf{R}\mathbf{F}(\phi)$ in the regions with $\phi\sim 1$, with off-diagonal components oriented along strips, and $\mathbf{F}\sim \mathbf{I}$ in the regions with $\phi\sim 0$. For instance, checking the values taken by the components of the deformation gradient in one cluster with $\phi\sim 1$ along the second column of Figure \ref{fig:4}, we observe the values $\mathbf{F}_{xx}\sim 0.98$, $\mathbf{F}_{xy}\sim 0.29$, $\mathbf{F}_{yx}\sim 0.20$, $\mathbf{F}_{yy}\sim 1.08$, which corresponds to
\[
\begin{pmatrix}
 \cos(\theta) && -\sin(\theta)\\
 \sin(\theta) && \cos(\theta)
\end{pmatrix}
\begin{pmatrix}
 1 && 0.5\\
 0 && 1
\end{pmatrix}
\]
with $\theta=0.2$ radiants.
In the case $\gamma=1$ we also observe higher deviations from the values $\text{det}\mathbf{F}=1$, concentrated at the boundary regions of the clusters with $\phi\sim 1$, then in the case $\gamma=0.001$, where $\text{det}\mathbf{F}\sim 1$ over the whole domain.

We then consider the initial conditions $\phi_0=0.7\pm 0.2\iota$ and $\mathbf{F}_0=\mathbf{I}$. In Figure \ref{fig:5} we show the numerical results at different time points, for both the cases with $\gamma=1$ and $\gamma=0.001$.
\begin{figure}[ht!]
\includegraphics[width=0.9\linewidth]
{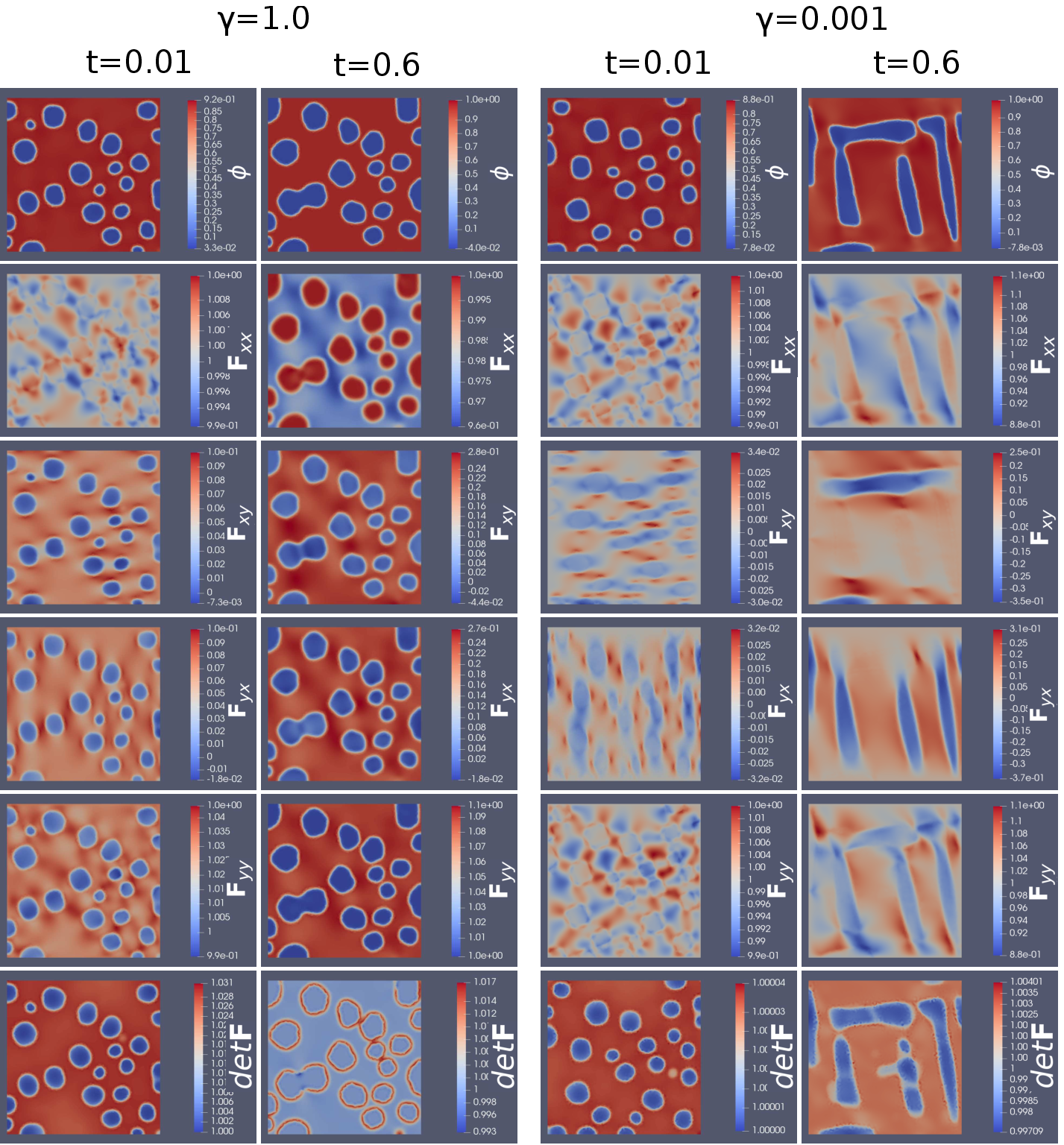}
\centering
\caption{Plot of $\phi$ (first row) and $\mathbf{F}_{xx}$ (second row), $\mathbf{F}_{xy}$ (third row), $\mathbf{F}_{yx}$ (fourth row), $\mathbf{F}_{yy}$ (fifth row), $\text{det}\mathbf{F}$ (sixth row) at different time points in the case $\phi_0=0.7\pm 0.2\iota$, for $\gamma=1$ (first and second columns) and $\gamma=0.001$ (third and fourth columns).}
\label{fig:5}
\end{figure}
We observe from Figure \ref{fig:5} that the phase separation dynamics for the $\phi$ variable consists in the formation of elongated circular clusters with $\phi\sim 0$ along orthogonal directions determined by the separation dynamics for the $\mathbf{F}$ variable, immersed in a bath with $\phi\sim 1$. As in the case reported in Figure \ref{fig:4}, at late times the clusters grow and merge, forming strips oriented in orthogonal directions and $\mathbf{F}$ assumes values $\mathbf{F}\sim \mathbf{R}\mathbf{F}(\phi)$ in the regions with $\phi\sim 1$, with off-diagonal components oriented along strips, and $\mathbf{F}\sim \mathbf{I}$ in the regions with $\phi\sim 0$.
 As in Figure \ref{fig:4}, in the case $\gamma=1$ we observe higher deviations from the values $\text{det}\mathbf{F}=1$, concentrated at the boundary regions of the clusters with $\phi\sim 0$, then in the case $\gamma=0.001$, where $\text{det}\mathbf{F}\sim 1$ over the whole domain.

\subsection{Test case 2 -- Coarsening}
We start by considering the initial conditions $\phi_0=1.0(\chi_{B1}+\chi_{B2})$, where $\chi_{B1}$ and $\chi_{B2}$ are the characteristic functions of two circular regions placed symmetrically along the $x$ direction, and
 \[
\mathbf{F}_0=\begin{pmatrix}
 1 && a\phi_0\\
 0 && 1
\end{pmatrix}.
\] 
In Figure \ref{fig:6} we show the numerical results at different time points, both for the cases $\zeta=1$ and $\zeta=10$.
\begin{figure}[ht!]
\includegraphics[width=0.9\linewidth]
{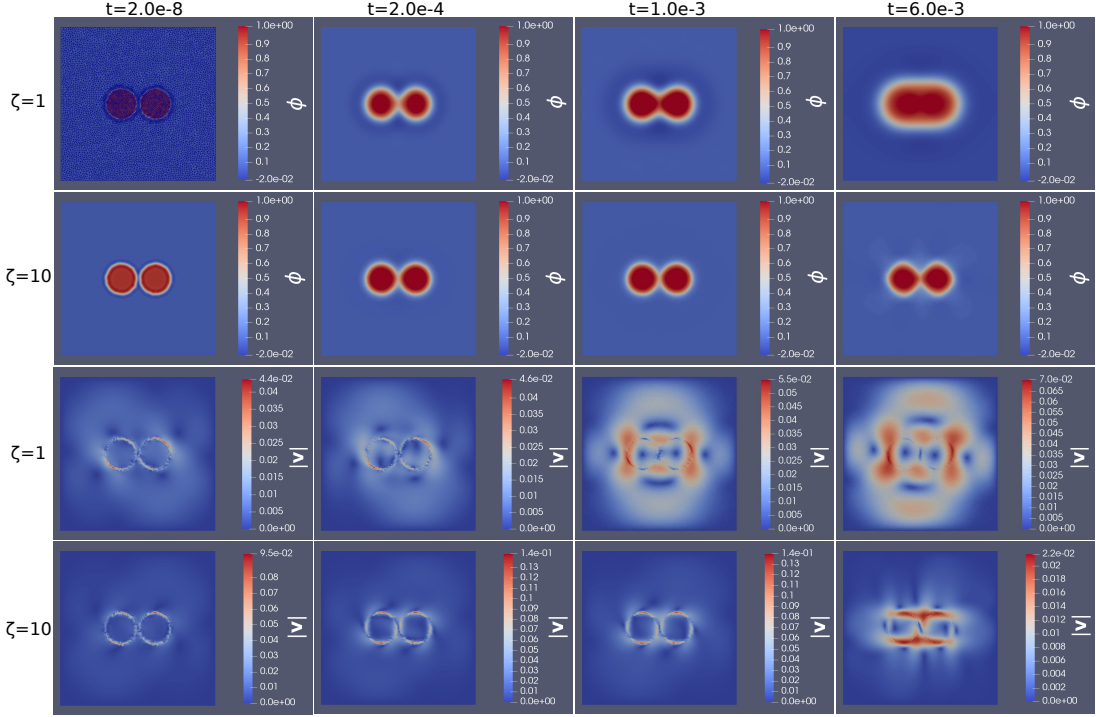}
\centering
\caption{Plot of $\phi$ (first and second rows) and $|\mathbf{v}|$ (third and fourth rows) at different time points, for $\zeta=1$ (first and third rows) and $\zeta=10$ (second and fourth rows), in the case $\phi_0=1.0(\chi_{B1}+\chi_{B2})$.}
\label{fig:6}
\end{figure}
We observe from Figure \ref{fig:6} that, for a low value of the elastic modulus $\zeta=1$, the initial circular clusters evolve, with interacting tips initially oriented along the bisecting directions of the plane and advected by the velocity field, finally merging into the equilibrium shape of an ellipse with the major axis oriented along the $x$ axis. Increasing the elastic modulus to $\zeta=10$, the circular clusters interact only weakly and do not merge during the observed time window, deforming to a rectangular shape. The results in Figure \ref{fig:6} may be compared to the results shown in \cite[Figures 6-7]{leo} for the merging of two circular precipitates in a linear isotropic inhomogeneous elastic medium. While in the latter case with the considered elastic parameters the equilibrium shape is circular, in the present case the anisotropy associated to the pure phase $\phi\equiv 1$ drives an elliptic or rectangular equilibrium shape. 

We finally consider the initial conditions $\phi_0=1.0(\chi_{B1}+\chi_{B2}+\chi_{B3}+\chi_{B4})$, where $\chi_{B1}$, $\chi_{B2}$, $\chi_{B3}$ and $\chi_{B4}$ are the characteristic functions of four circular regions placed symmetrically with respect to the center of the domain, and
 \[
\mathbf{F}_0=\begin{pmatrix}
 1 && a\phi_0\\
 0 && 1
\end{pmatrix}.
\] 
In Figure \ref{fig:7} we show the numerical results at different time points, both for the cases $\zeta=1$ and $\zeta=10$.
\begin{figure}[ht!]
\includegraphics[width=0.9\linewidth]
{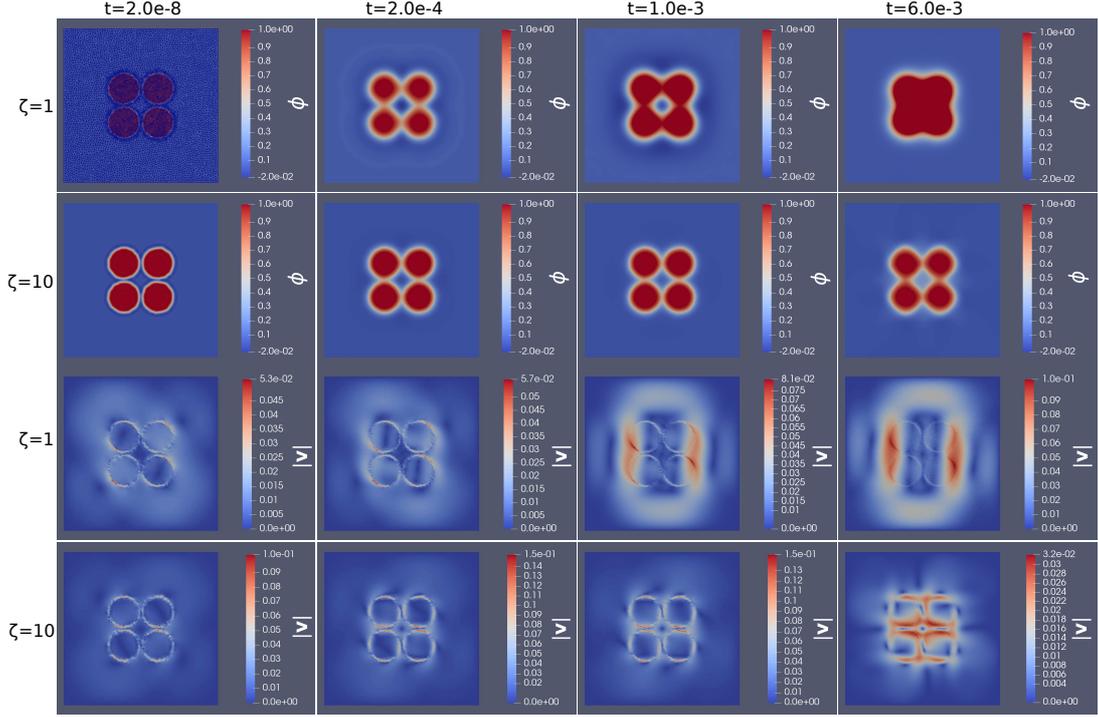}
\centering
\caption{Plot of $\phi$ (first and second rows) and $|\mathbf{v}|$ (third and fourth rows) at different time points, for $\zeta=1$ (first and third rows) and $\zeta=10$ (second and fourth rows), in the case $\phi_0=1.0(\chi_{B1}+\chi_{B2}+\chi_{B3}+\chi_{B4})$.}
\label{fig:7}
\end{figure}
We observe from Figure \ref{fig:7} that, for $\zeta=1$, the initial circular clusters evolve, with interacting tips initially oriented along the bisecting directions of the plane and advected by the velocity field, similarly to the case in Figure \ref{fig:6}. The clusters merge both horizontally and vertically, initially surrounding a circular region with $\phi\equiv 0$ which dissolves at late times. The final equilibrium shape is given by a square. This is different from the results with linear isotropic inhomogeneous elasticity and periodic boundary conditions reported in \cite[Figure 11]{leo}, where the final equilibrium configuration is given by a cross.  Also, in the case with $\zeta=10$ the circular clusters interact only weakly and do not merge during the observed time window, deforming into a rectangular shape. 

We highlight the fact that the numerical results shown in \cite[Figures 6-7-11]{leo} concern interacting soft precipitates in a linear isotropic inhomogeneous elastic medium. In the case with hard precipitates, i.e. when the elastic modulus associated to the phase $\phi\equiv 1$ is higher than the elastic modulus associated to the phase $\phi\equiv 0$, the numerical results in \cite{leo} show a repulsion between the precipitates, which do not interact. This is comparable to what happens in our numerical simulations in the case of an high value of the elastic modulus.

\section{Conclusion}
In this paper we have proposed  a new Cahn--Hilliard phase field model coupled to nonlinear incompressible finite viscoelasticity, where a new kind of diffusive regularization, of Allen--Cahn type, is introduced in the transport equation for the deformation gradient, together with a regularizing interface term depending on the gradient of the deformation gradient in the free energy density of the system. The designed regularization, which preserves the dissipative structure of the equations, helps in enhancing the space and time regularity of the deformation gradient. The resulting transport equation for the deformation gradient with Allen--Cahn type regularization was expressed in a dual mixed formulation, introducing a dual variable of the deformation gradient which enters also in the expression for the Cauchy stress tensor. Through a Galerkin approximation of the model equations and the introduction of truncated problems, where the polynomial growth of the elastic energy density is truncated to degree $4$, we proved existence of a global in time weak solution in three space dimensions and for elastic energy densities which are coupled to the phase field variable and which possibly degenerate for some values of the phase field variable. Then, thanks to an iterative argument based on elliptic regularity bootstrap steps applied to the Allen--Cahn transport equation for the deformation gradient, we extended the existence result, passing to the limit for the truncation parameter tending to infinity, to the case of a Cahn--Hilliard potential and an elastic energy density with maximum allowed polynomial growths. In three space dimensions, we found maximum admissible polynomial growth degrees of $10$ for the Cahn--Hilliard potential and $6$ for the elastic energy density, with the degree of the Cahn--Hilliard potential depending on the degree of the elastic energy density. We also proposed two kind of unconditionally energy stable and efficient finite element approximations of the model, based on convex splitting ideas and on the use of a scalar auxiliary variable, proving the existence and stability of discrete solutions. We finally showed numerical results for different test cases with shape memory alloy type free energies, characterized by different elastic properties of the pure phases of the phase field variable, which verify the gradient stability properties of the proposed schemes and show qualitatively how the topology of stationary states depends on both the phase separation and the elasticity dynamics.
Future developments of the present work will investigate the cases with singular phase field potential and compressible elasticity, together with the generalization of the proposed model to models for biomathematical applications.}

\section{Acknowledgements}
This research was supported by the Italian Ministry of Education,
University and Research (MIUR): Dipartimenti di Eccellenza Program (2018--2022)
-- Dept.~of Mathematics ``F.~Casorati'', University of Pavia.
In addition, AA, PC and ER  gratefully mention some other support
from the MIUR-PRIN Grant 2020F3NCPX ``Mathematics for industry 4.0 (Math4I4)'' and
their affiliation to the GNAMPA (Gruppo Nazionale per l'Analisi Matematica,
la Probabilit\`a e le loro Applicazioni) of INdAM (Isti\-tuto
Nazionale di Alta Matematica). 

%
\bibliographystyle{plain}
\bibliography{biblio_CHVE} 
\end{document}